\pdfoutput=1
\documentclass[10pt]{article}
\usepackage{cancel}
\usepackage{color}
\usepackage{graphicx}
\usepackage[colorlinks=true, pdfstartview=FitV, linkcolor=blue, citecolor=blue,urlcolor=blue]{hyperref}

\usepackage{framed,color}

\usepackage{amsbsy}
\usepackage{amssymb,amsmath}
 \usepackage{amsfonts}
\usepackage{graphicx}

\usepackage{cancel}
\usepackage{mathtools}
\usepackage{amsmath}
\usepackage{amsthm}
\usepackage{amsfonts}
\usepackage{amssymb}
\usepackage{mathrsfs}
\usepackage{graphicx}
\usepackage[T1]{fontenc}
\usepackage[OT2,T1]{fontenc}
\usepackage[latin1]{inputenc}
 \usepackage{lmodern}
\usepackage[all,cmtip]{xy}

\usepackage{bbm}
\usepackage[vcentermath]{youngtab}
\usepackage{ stmaryrd }
\usepackage{calc}

\usepackage{setspace}

\usepackage{epigraph}
\usepackage{mathdots}
\usepackage{mathabx}
\usepackage{braket}
\usepackage{longtable}
\usepackage{multirow}
\usepackage{float}
\usepackage{tikz-cd}

\def\eqref#1{(\ref{#1})}
\newtheorem{theorem}{Theorem}[section]
\newtheorem{example}{Example}[section]
\newtheorem{exercise}{Exercise}[section]

\newtheorem{lemma}{Lemma}[section]
\newtheorem{remark}{Remark}[section]

\newtheorem{proposition}{Proposition}[section]
\newtheorem{corollary}{Corollary}[section]
\newtheorem{definition}{Definition}[section]

\def\bre{\begin{remark}}
\def\ere{\end{remark}}
\def\bth{\begin{theorem}}
\def\eth{\end{theorem}}
\def\bcr{\begin{corollary}}
\def\ecr{\end{corollary}}
\def\bex{\begin{example}\small}
\def\eex{\end{example}}
\def\bexr{\begin{exercise}\small}
\def\eexr{\end{exercise}}
\def\ble{\begin{lemma}}
\def\ele{\end{lemma}}

\def\bde{\begin{definition}}
\def\ede{\end{definition}}
\def\bpr{\begin{proposition}}
\def\epr{\end{proposition}}
\def\be{\begin{equation}}
\def\ee{\end{equation}}
\def\bea{\begin{eqnarray}}
\def\eea{\end{eqnarray}}
\def\beas{\begin{eqnarray*}}
\def\eeas{\end{eqnarray*}}

\def\vPsi{\vec{\Psi}}

\newcommand{\minitab}[2][l]{\begin{tabular}#1 #2\end{tabular}}

\newlanguage\fakelanguage
\newcommand\cyr{\fontencoding{OT2}\fontfamily{wncyr}\selectfont
   \language\fakelanguage}
\DeclareTextFontCommand{\textcyr}{\cyr}

\numberwithin{equation}{section}

\numberwithin{theorem}{section}

\numberwithin{proposition}{section}

\numberwithin{definition}{section}

\numberwithin{remark}{section}

\numberwithin{lemma}{section}

\numberwithin{corollary}{section}

\date{January 19, 2021}
\begin{document}
\baselineskip=14pt

\vspace{0.2cm}
\begin{center}
\begin{Large}
\fontsize{17pt}{27pt}
\selectfont

\textbf{The sixth Painlev\'e equation as isomonodromy deformation of an irregular system:  monodromy data, coalescing eigenvalues, locally holomorphic transcendents and Frobenius manifolds}
\end{Large}
\\
\bigskip
\begin{large} {Gabriele Degano${}^1$, Davide Guzzetti${}^{2,3}$ }\end{large}
\bigskip
\\{${}^1$ Grupo de F\'isica Matem\'atica, Universidade de Lisboa, Campo Grande, Edif\'icio C6, 1749-016 Lisboa  -- Portugal. E-mail:  gabriele.degano.gd@gmail.com}
\\{${}^2$ SISSA,  Via Bonomea, 265,  34136 Trieste -- Italy.  E-mail: guzzetti@sissa.it}
\\{${}^3$ Davide Guzzetti's ORCID ID: 0000-0002-6103-6563}
%\bigskip

\end{center}
%\begin{spacing}{1.5}
%\begin{small}
{\bf Abstract:}  We consider a 3-dimensional Pfaffian system, whose $z$-component is a differential system with  irregular singularity at infinity  and Fuchsian  at zero. In the first part of the paper, we  prove that  its Frobenius integrability is equivalent to the sixth Painlev\'e equation PVI. The coefficients of the system will be explicitly written in terms of the solutions of PVI.  In this way, we remake a result of  \cite{Hard,MazzoIrr}. 
We then express in terms of the Stokes matrices of the $3\times 3$  irregular system the  monodromy invariants $p_{jk}=\hbox{\rm Tr}(\mathcal{M}_j\mathcal{M}_k)$ of the 2-dimensional isomonodromic Fuchsian system with four singularities, traditionally associated to PVI \cite{Fuchs,JM} and used to solve the non-linear connection problem.    
 Several years after  \cite{Hard,MazzoIrr},  the authors of  \cite{CDG} showed that  the computation of the monodromy data of a class of 
 irregular systems may be facilitated  in case of coalescing eigenvalues. This coalescence corresponds to the critical points
  (fixed singularities) of PVI. In the second part of the paper, we  classify  the branches of PVI transcendents  holomorphic at a critical point 
   such that    the analyticity  and semisimplicity properties described in  \cite{CDG} are satisfied, and we compute the associated  
   Stokes matrices and the invariants $p_{jk}$. Finally, we compute the monodromy data parametrizing the chamber of a 
   3-dim Dubrovin-Frobenius manifold associated with a transcendent holomorphic at $x=0$.

\vspace{0.5cm}

\noindent
\begin{small}
\end{small}
%\vskip 15pt

%\tableofcontents

%{\color{blue} ESEMPIO DI USO DEL COLORE} 
%{\color{red} ESEMPIO DI USO DEL COLORE} 
%{\color{green} ESEMPIO DI USO DEL COLORE} 
%{\color{magenta} ESEMPIO DI USO DEL COLORE}

%\setlength{\belowdisplayskip}{4.6pt} \setlength{\belowdisplayshortskip}{4.6pt}
%\setlength{\abovedisplayskip}{4.64pt} \setlength{\abovedisplayshortskip}{4.64pt}

\section{Introduction}
\vskip 0.2 cm 
The sixth Painlev\'e equation, hereafter denoted by PVI, is the non-linear ODE
$$
\frac{d^2y }{ dx^2}=\frac{1}{ 2}\left[ 
\frac{1}{ y}+\frac{1}{ y-1}+\frac{1}{ y-x}
\right]
           \left(\frac{dy}{ dx}\right)^2
-\left[
\frac{1}{ x}+\frac{1}{ x-1}+\frac{1}{ y-x}
\right]\frac{dy }{ dx}
$$
$$
+
\frac{y(y-1)(y-x)}{ x^2 (x-1)^2}
\left[
\alpha+\beta \frac{x}{ y^2} + \gamma \frac{x-1}{ (y-1)^2} +\delta
\frac{x(x-1)}{ (y-x)^2}
\right]
,
$$ 
where the  coefficients can be parameterized  by four constants $\theta_1,\theta_2,\theta_3,\theta_\infty$, with $\theta_\infty\neq 0$, as follows
\be
\label{11giugno2021-2}
2\beta=-\theta_1^2,~~~2\delta=1-\theta_2^2,~~~2\gamma=\theta_3^2,~~~2\alpha=(\theta_\infty-1)^2.
\ee
PVI has three fixed singularities at $x=0,1,\infty$,   called {\it critical points}, and the behaviour of its solutions  at these points is called {\it critical}. Since these singularities are equivalent by the symmetries of PVI \cite{Okamoto}, as far as the local analysis is concerned it suffices to study the point $x=0$.\footnote{We can use the following symmetries
\begin{align*}
&1)&&\theta_3^{\prime}=\theta_1,&&\theta_1^{\prime}=\theta_3;&
\theta_2^{\prime}=\theta_2,&&\theta_\infty^{\prime}=\theta_\infty;   
&&y^{\prime}(x^\prime)=1-y(x), &&x=1-x^\prime. 
\\
&2)&& \theta_1^{\prime}=\theta_\infty-1,&&\theta_\infty^{\prime}=\theta_1+1;&
\theta_3^{\prime}=\theta_3,&&\theta_2^{\prime}=\theta_2;&&
y^{\prime}(x^\prime)={1\over y(x)},&&x={1\over x^{\prime}}.
\\
&3)&& \theta_2^\prime=\theta_3,&&\theta_3^\prime=\theta_2;&
\theta_1^\prime=\theta_1,&&\theta_\infty^\prime=\theta_\infty; &&
y^\prime(x^\prime)={1\over x}y(x), && x={1\over x^\prime}.
\end{align*}
}
 The solutions of PVI are  called {\it transcendents}, because they generically are highly transcendental functions not expressible in terms of classical functions by means of Umemura's admissible operations \cite{Um1,Um2,Um3,Watanabe}. 
The critical behaviours of PVI trascendents have been obtained in  \cite{Jimbo,Sh82-1,Sh82-2,Sh82-3,Sh87,guz2001-1,guz2002,guz2006,guz2008,guz2011,guz2012,guz2012-1,Sh15} and classified and tabulated in  \cite{guz2012} (see also \cite{guz2015}), including critical behaviours given by  Taylor expansions. The latter   have also been  studied in \cite{kaneko,kaneko07,kaneko09}, where  convergence has been proved. See also  \cite{Bruno4,Bruno5,Bruno7} for expansions of PVI transendents, and  \cite{GG15, GG16, GG17,GG20-1,GG20-2,Sh15} for convergence issues.

PVI is equivalent to  the isomonodromy deformation equations, the {\it Schlesinger equations},   of a
 $2\times 2$ isomonodromic Fuchsian system 
\be
\label{23luglio2021-4}
\frac{d \Phi}{d\lambda} =\sum_{k=1}^3 \frac{\mathcal{A}_k(u)}{\lambda-u_k} \Phi
\ee
with $u=(u_1,u_2,u_3)$ and 
$$ \hbox{eigenvalues of $\mathcal{A}_k$}= \pm \frac{\theta_k}{2},
\quad \quad 
 \sum_{k=1}^3  \mathcal{A}_k=
   \begin{pmatrix}
  - \theta_\infty/2& 0
     \\
      0 & \theta_\infty/2
   \end{pmatrix}.
   $$
   This means that  solutions   $\mathcal{A}_1(u),\mathcal{A}_2(u),\mathcal{A}_3(u)$ of the Schlesinger equations, up to the equivalence relation $   \mathcal{A}_k\longmapsto \hbox{\rm diag}(\varepsilon_1,\varepsilon_3)^{-1}\cdot\mathcal{A}_k\cdot \hbox{\rm diag}(\varepsilon_1,\varepsilon_3)$, $\varepsilon_1\varepsilon_3\neq 0$,  are in one-to-one correspondence with solutions of PVI, with $x=(u_2-u_1)/(u_3-u_1)$. The correspondence is  given by the
    explicit formulae in  appendix C of \cite{JM} (in  \cite{JM}, $x$ is called $t$ and the choice $(u_1,u_2,u_3)=(0,t,1)$ is made).
      
   Given a branch of a  PVI transcendent, defined in the $x$-plane with cuts, such as $|\arg x|<\pi$, $|\arg(1-x)|<\pi$,  the non-linear connection problem is to express the (one or two) integration constants parametrizing the critical behaviour of the branch at a critical point, in terms of the integration constants expressing the critical behaviour of the branch at another critical point.
  Since Jimbo's work \cite{Jimbo} 
    (see also the review  \cite{guz2015}) the strategy to solve this problem   
    has been to explicitly express   the  integration constants  parameterizing the three  critical behaviours   at $x=0, 1 , \infty$ respectively   in terms of {\it the same} traces
   \be
\label{27luglio2021-1}
p_{jk}=\hbox{tr}(\mathcal{M}_j\mathcal{M}_k),\quad 1\leq j\neq k \leq 3, 
\ee
where $p_{jk}=p_{kj}$ and   $\mathcal{M}_1,\mathcal{M}_2,\mathcal{M}_3\in SL(2,\mathbb{C})$ are the monodromy matrices  of a fundamental matrix solution
  $\Phi(\lambda,u)$ of the system \eqref{23luglio2021-4},   whose $\mathcal{A}_1,\mathcal{A}_2,\mathcal{A}_3$ are in one-to-one correspondence (up to the above  equivalence) with the transcendent. Here,  $\Phi(\lambda,u)$ is defined for $u$ varying in a sufficiently small    simply  connected domain  of $\mathbb{C}^3\backslash \bigcup_{j\neq k }\{u_j=u_k\}$,  and for $\lambda$  in  a plane 
   with brach-cuts issuing from $u_1,u_2,u_3$ towards infinity,   in a direction that 
   we will call  $3\pi/2-\tau^{*}\in\mathbb{R}$, satisfying $3\pi/2-\tau^{*}\neq \arg(u_j-u_k)$ mod $\pi$,  $1 \leq j\neq k\leq 3$. 
         
      The $p_{jk}$, together with the $p_j:=\hbox{tr}(\mathcal{M}_j)=2\cos(\pi \theta_j)$, 
      $j=1,2,3$, and\footnote{The sequence  $j_1\neq j_2\neq j_3$  depends on the ordering of the basis of loops for the monodromy matrices. More precisely, $ j_1\prec j_2\prec j_3$, where the ordering is defined  in  \eqref{16aprile2021-3}.} $p_\infty:=\hbox{tr}(\mathcal{M}_{j_3}\mathcal{M}_{j_2}\mathcal{M}_{j_1})=2\cos(\pi\theta_\infty)$ generate the ring 
      of {invariant functions} of $\boldsymbol{\mathcal{M}}:=SL(2,\mathbb{C})^3/SL(2,\mathbb{C})$,  the space of 
      conjugacy classes  of  triples $\mathcal{M}_1,\mathcal{M}_2,\mathcal{M}_3$.  A conjugacy class  belonging to a ``big'' 
      open subset\footnote{ Let $(i,j,k)$ be a cyclic permutation of $(j_1,j_2,j_3)$. 
       According to \cite{Iwa}, the big open is the complement in $\boldsymbol{\mathcal{M}}$ of the set where the following six algebraic equations are satisfied:
      \[
\begin{cases}
 (p_{jk}^2-4)(p_{jk}^2+p_i^2+p_\infty^2-p_{jk}p_ip_\infty-4)=0, \\[2ex]
 (p_{jk}^2-4)(p_{jk}^2+p_j^2+p_k^2-p_{jk}p_jp_k-4)=0,
\end{cases} \quad \hbox{ for all $i=1,2,3$}
\]
For example, the above polynomials vanish if  $\mathcal{M}_1,\mathcal{M}_2,\mathcal{M}_3$ generate a reducible group.  Also the triple $(p_{12}, p_{23}, p_{31})=(\pm 2, \pm 2,\pm2)$,  is a common root of the above polynomials.  
      }
       of $\boldsymbol{\mathcal{M}}$ can be explicitly parameterized  by  the $p_{jk}$, $p_j$, $p_\infty$, 
      according to Tables 1 and 2 of \cite{Iwa} (a generalization of the formulae of Tables 1 and 2 of \cite{Iwa} for the Garnier system with two times is given in \cite{MC}). 
 There  is a one-to-one correspondence between monodromy  data and branches  of Painlev\'e transcendents, if  none of the $\mathcal{M}_j$, $j=1,2,3$ and  
$\mathcal{M}_{j_3}\mathcal{M}_{j_2}\mathcal{M}_{j_1}$ is  equal to  $\pm I$, where $I$ is the identity matrix (see \cite{guz2008} and also \cite{DM1}). In this case,    the integration constants expressing the critical behaviours  can be univocally written in terms of the $p_{jk}$, provided   the triple  $\mathcal{M}_1,\mathcal{M}_2,\mathcal{M}_3$ is in the big open  of $\boldsymbol{\mathcal{M}}$ (an example of a point not in the big open is a triple generating a reducible group, another example is when $(p_{12}, p_{23}, p_{31})=(\pm 2, \pm 2,\pm2)$;  see also   Section \ref{13agosto2021-1}).

The analytic continuation of a branch is obtained by an action of the braid group on the monodromy data \cite{DM,guz2002,Iwa,Boalch1}.

\vskip 0.2 cm 

It was  shown in   \cite{Hard} that PVI is also the isomonodromy condition for the $3\times 3$ irregular system \eqref{23luglio2021-2} below. In  \cite{Hard} a  class of integrable Pfaffian systems (studied in   \cite{JMMS})  are described  in the loop algebra framework of \cite{AHP} and,  using duality of moment maps \cite{AHP1}, dual isomonodromic systems are obtained. In particular, the "dual" to  \eqref{23luglio2021-4} turns out to be\footnote{ The notations in  in \cite{Hard} are different from ours.} \be
\label{23luglio2021-2}
\frac{dY}{dz}=\left(U+\frac{V(u)}{z}\right) Y,\quad \quad U=\hbox{\rm diag}(u_1,u_2,u_3),
\ee
with a suitable matrix $V(u)$. Here, $z=0$ is a Fuchsian singularity, and $z=\infty$ is a singularity of the second kind with Poincar\'e rank 1.  
     By another approach based on the Laplace transform of   \cite{BJL4}, it was equally proved that the isomonodromy deformation equations of \eqref{23luglio2021-2}  reduce to PVI, in \cite{Dub1,Dub2} for  $\theta_1=\theta_2=\theta_3=0$,   in \cite{MazzoIrr}  for the general case. This was later shown also in \cite{Boalch1}.  In particular, \cite{MazzoIrr} gives   the one-to-one correspondence between solutions of PVI and  matrices  $V$ satisfying the isomonodromy conditions  (up to the equivalence relation $V\longmapsto 
\hbox{\rm diag}(k_1^0,k_2^0,k_3^0)\cdot V \cdot \hbox{\rm diag}(k_1^0,k_2^0,k_3^0)^{-1}$, $ k_1^0k_2^0k_3^0\neq 0$) by means of formulae $V=V(x,y(x))$  (up to  typing  misprints) which  make completely explicit the results of  section 3.c. of  \cite{Hard}.

Though  it was suggested in  \cite{MazzoIrr}  that the point of view based on \eqref {23luglio2021-2} may be useful, for example for the study of symmetries, an this was further investigated in \cite{Boalch1}, in the literature the isomondormic approach to  PVI has  mainly been based on the Fuchsian system \eqref{23luglio2021-4}, even in the most recent developments \cite{GL1,ILP1,ILT1,ILST1,ILTy1}.
 
   \vskip 0.2 cm 
 In this paper, we  show some advantages and uses  of  \eqref {23luglio2021-2}. First of all, we  show that \eqref{23luglio2021-2} can be employed  to solve  the non-linear connection problem. Indeed, we  prove an explicit {\it  formula expressing  the 2-dimensional  invariants $p_{jk}$ in terms of the 3-dimensional Stokes matrices of \eqref{23luglio2021-2}}.  This  problem was not considered in \cite{MazzoIrr} and -- to our knowledge --   in the literature.\footnote{What we can find in the literature is theorem I (theorem 1) at page 389 of \cite{guz2016}, where the relation is given between Stokes matrices of the irregular system of dimension $n$ and traces of products of monodromy matrices of certain selected solutions of the Laplace-transformed Fuchsian system {\it of the same dimension}. Moreover, in theorem 2 of \cite{Boalch1}, and in \cite{Dub1,Dub2} in case of Frobenius manifolds, a Killing-Coxeter identity is given relating a  product of Stokes matrices and a product of  monodromy matrices (pseudo-reflections) of the Laplace-transformed Fuchsian system of the same dimension.}  
  
   Then, we  demonstrate  that \eqref {23luglio2021-2}  allows to implement in case of PVI some results of \cite{CDG} on isomonodromy deformations with coalescing eigenvalues,  naturally yielding a  {\it classification} of   certain transcendents admitting a holomorphic  branch at a critical point, whose associated matrix $V(u)$ satisfies the analyticity properties of \cite{CDG}. 
  We will also  show   that  \eqref{23luglio2021-2} is advantageous for the computation of the monodromy data parametrizing  the above classified  transcendents, because  the properties  \cite{CDG} make the computation feasible in terms of classical special functions. Several examples of these explicit computations will be given, including an application to Dubrovin-Frobenius manifolds.

 %%%%%%%%%%%%%%%%%%%%%%%%%%
 \subsubsection*{Results} 
 
{\bf (1)} Firstly,  we    remake the results of  \cite{MazzoIrr} and of section 3.c. of \cite{Hard}, by means of a Laplace transform of the type of \cite{BJL4,guz2016,guz2021}  with deformations parameters, relating the two Pfaffian systems \eqref{19febbraio2021-11} and \eqref{15ottobre-bloccoGreenPass-1} introduced in Section \ref{19giugno2021-6}.   {\bf  Theorem \ref{20gennaio2021-2}}  states  the equivalence between PVI and  the Frobenius integrability  of  \eqref{19febbraio2021-11}. More precisely, we show that $V(u)$ can always  be  factorized as
$$
V(u)=(u_3-u_1)^{\hbox{\rm diag}(\theta_1,\theta_2,\theta_3)}~ \Omega(x)~ (u_3-u_1)^{-\hbox{\rm diag}(\theta_1,\theta_2,\theta_3)},\quad\quad x=\frac{u_2-u_1}{u_3-u_1}.
$$
The Frobenius integrability of  \eqref{19febbraio2021-11} will be  expressed by a non-linear differential system
\be
\label{23luglio2021-1}
\frac{d \Omega}{d x}
=\left[\Omega,\widehat{\Omega}_2 \right]
,\quad\quad
 \widehat{\Omega}_2 :=
\begin{pmatrix}
0 & -\frac{\Omega_{12}}{x} &0
\\ 
-\frac{\Omega_{21}}{x} & 0 & \frac{\Omega_{23}}{1-x}
\\ 
0 &\frac{\Omega_{32}}{1-x} & 0
\end{pmatrix},
\ee
  that we prove to be equivalent to PVI, in the sense that there is a one-to-one correspondence between solutions $y(x)$ and  equivalence classes of solutions $\Omega(x)$  of \eqref{23luglio2021-1} satisfying 
$$
\hbox{\rm eigenvalues of }\Omega(x)= 0, \frac{\theta_\infty-\theta_1-\theta_2-\theta_3}{2},  \frac{-\theta_\infty-\theta_1-\theta_2-\theta_3}{2},\quad \quad \hbox{\rm diag}(\Omega)=-\hbox{\rm diag}(\theta_1,\theta_2,\theta_3).
$$ 
The equivalence classes are defined by the relation  $\Omega
\longmapsto   K^0\cdot \Omega \cdot (K^0)^{-1}$, where $K^0$ is a constant diagonal matrix with non zero eigenvalues.  
   We write the explicit  formulae expressing  a matrix $\Omega(x)$ (up to $\Omega \longmapsto   K^0\cdot \Omega \cdot (K^0)^{-1}$) in terms of a Painlev\'e transcendent $y(x)$, that will be  used for actual computations in the second part of the paper.  They  correct some typing misprints in the analogous formulae of \cite{MazzoIrr}, making them usable.

\vskip 0.2 cm 

{\bf (2)} The second result  is  {\bf Theorem \ref{16aprile2021-5}}, which gives the invariants $p_{jk}$ in terms of the entries of the Stokes matrices of \eqref{23luglio2021-2}. In the first part of the theorem, we show that 
\be
\label{23luglio2021-3}
p_{jk}
   =
   \left\{
  \begin{array}{cc}
 2\cos\pi(\theta_j-\theta_k)-e^{i\pi(\theta_j-\theta_k)} (\mathbb{S}_1)_{jk}(\mathbb{S}_{2}^{-1})_{kj}, & j\prec k,
  \\
  \\
 2\cos\pi(\theta_j-\theta_k)-e^{i\pi(\theta_k-\theta_j)}(\mathbb{S}_1)_{kj}(\mathbb{S}_{2}^{-1})_{jk},
  , & j\succ k.
  \end{array}
  \right.,
\ee
where the ordering relation 
$  j\prec k $ means $ \Re(e^{i\tau^{*}} (u_j-u_k))<0$. 
Formula \eqref{23luglio2021-3} is new in the literature, as far as we are informed. Here, $\mathbb{S}_1$ and $\mathbb{S}_2$ are two successive Stokes matrices of \eqref{23luglio2021-2}, defined by $Y_{j+1}=Y_j \mathbb{S}_j$, $j=1,2$, where $Y_1(z,u),Y_2(z,u),Y_3(z,u)$ are the unique fundamental matrix solutions with canonical asymptotic behaviour given by the formal solution 
\be \label{3agosto2021-5}
Y_F(z,u)=(I+\sum_{k\geq 1} F_k(u)z^{-k})z^{\hbox{\rm diag}(V)} \exp\{zU\},
\ee  
in respectively three successive overlapping sectors  (for small $\varepsilon>0$)
\begin{align}
\label{27luglio2021-5}
&\mathcal{S}_1:\quad (\tau^{*}-\pi)-\varepsilon<\arg z <  \tau^{*}+\varepsilon,   
 \\
 \label{27luglio2021-6}
& \mathcal{S}_2:\quad \tau^{*}-\varepsilon<\arg z <  (\tau^{*}+\pi)+\varepsilon, 
  \\
  \label{27luglio2021-7}
& \mathcal{S}_3:\quad (\tau^{*}+\pi)-\varepsilon<\arg z <  (\tau^{*}+2\pi)+\varepsilon .
 \end{align}
 The intersections  $\mathcal{S}_1\cap \mathcal{S}_2$ and  $\mathcal{S}_2\cap \mathcal{S}_3$ do not contain Stokes rays of $U(u)$, as $u$ varies in a small subset of $\mathbb{C}^3\backslash \bigcup_{j\neq k }\{u_j=u_k\}$ as mentioned above. 

In the second part of Theorem \ref{16aprile2021-5} we establish the same formulae also in case $u$ varies in a sufficiently small polydisc $\mathbb{D}(u^c)$ centered at a {\it coalescence point} $u^c$, such that $u_j^c=u_k^c$ for some $j\neq k$, under the asumption that $V(u)$ is holomorphic in $\mathbb{D}(u^c)$ and satisfies the vanishing conditions 
$$ 
V_{jk}(u)\to 0 \hbox{ holomorphically, for $u_j-u_k\to 0$ in $\mathbb{D}(u^c)$}.
$$ 
 We only need to consider the case when  two out of the three $u_1,u_2,u_3$ may coalesce at $u^c$, because if  $u_1^c=u_2^c=u_3^c$ the results are trivial. Moreover, at the critical points $x=0,1,\infty$   only two out of the three $u_1,u_2,u_3$  coalesce, because $x=(u_2-u_1)/(u_3-u_1)$. 
We show that in  the coalescent case formulae \eqref{23luglio2021-3} hold for $j\neq k$ such that $u_j^c\neq u_k^c$, while\footnote{In this case, $(\mathbb{S}_1)_{jk}= (\mathbb{S}_1)_{kj}= (\mathbb{S}_2)_{jk}= (\mathbb{S}_2)_{kj}=0 $ for $j\neq k$ such that $u_j^c= u_k^c$.} 
$$ 
p_{jk}=  2\cos\pi(\theta_j-\theta_k)\quad \hbox{ for $j\neq k$ such that $u_j^c= u_k^c$}.
$$
 Here, the $p_{jk}$ are still defined for a fundamental matrix solution  $\Phi(\lambda,u)$ with $u$  varying in a small subset of $\mathbb{D}(u^c)$ in whose interior $u_1,u_2,u_3$ are pairwise distinct. In this case, $\tau^*$ need to be chosen  so to satisfy $3\pi/2-\tau^*\neq \arg(u_j^c-u_k^c)$ mod $\pi$, for $j\neq k$ such that $u_j^c\neq u_k^c$. 

\bre {\rm
In the paper, $\tau^*$ will be called $\tau^{(0)}$ in case there are not coalescences, and $\tau$ in case of coalescences. 
}\ere

The proof of Theorem \ref{16aprile2021-5} is based   on the  isomonodromic Laplace
 transform \cite{guz2021}, relating an $n$ dimensional isomonodromic system of the type \eqref{23luglio2021-2} to an isomonodromic  Fuchsian system of {\it the same} dimension, with poles at $\lambda=u_k$, $k=1,...,n$. The latter admits certain
 {\it selected} vector solutions, whose monodromy can be written in terms of certain   {\it connection coefficients}. In \cite{guz2021}, the connection coefficients are expressed as linear functions of the entries of the 
  Stokes matrices of the $n$-dimensional  \eqref{23luglio2021-2} and vice-versa, including the coalescent case.  
 While in \cite{guz2021} both the irregular system \eqref{23luglio2021-2} and the Fuchsian system have the same dimension,   in our case the irregular system \eqref{23luglio2021-2}  has $n=3$ , while the Fuchsian system \eqref{23luglio2021-4} is $2$-dimensional. By a dimensional reduction,   we will express the monodromy matrices of the $2$-dimensional  system \eqref{23luglio2021-4} in terms of the connection coefficients of the $3$-dimensional Fuchsian system associated with \eqref{23luglio2021-2} by Laplace transform.  
 %A technical difficulty will be that the selected solutions of the $3$-dimensional Fuchsian system may be not  linearly independent. 

\vskip 0.2 cm

{\bf (3)} In order to introduce the third result, we need to recall the extension of the isomonodromy deformation theory given in \cite{CDG} (see also \cite{Eretico1,Guzz-SIGMA,Guzz-Spring18,Guzz-Degru20}).   In theorem 1.1. of \cite{CDG} an irregular $n\times n$ system of type \eqref{23luglio2021-2} is considered, with $U=\hbox{\rm diag}(u_1,...,u_n)$ and $V(u)$ holomorphic in a sufficiently small polydisc $\mathbb{D}(u^c)$  centered at a coalescence point $u^c$. Consequently, the polydisc contains a coalescence  locus where some components of $u=(u_1,...,u_n)$ merge.  The polydisc is so small that $u^c$ is the most coalescent point.  If the following vanishing conditions hold in $\mathbb{D}(u^c)$
\be
\label{5agosto2021-1} 
V_{jk}(u)\longrightarrow 0 \quad \hbox{ for $j\neq k$, whenever  $u_j-u_k\to 0$},
\ee
then the matrix coefficients  $F_k(u)$ of the  $n$-dimensional formal solution  \eqref{3agosto2021-5} 
 are holomorphic in $\mathbb{D}(u^c)$. Also  the canonical solutions $Y_1(z,u),Y_2(z,u),Y_3(z,u)$ are well defined and holomorphic in $\mathbb{D}(u^c)$, having the asymptotic expansion  $Y_F(z,u)$ in sectors $\mathcal{S}_1,\mathcal{S}_2,\mathcal{S}_3$ as in  \eqref{27luglio2021-5}-\eqref{27luglio2021-7}, 
with $\tau^*$ such that $3\pi/2-\tau^*\neq \arg(u_j^c-u_k^c)$ mod $\pi$, for $j\neq k$ such that $u_j^c\neq u_k^c$. Moreover, if \eqref{5agosto2021-1} holds, the monodromy data of \eqref{23luglio2021-2} are well defined and  constant  on the whole  $\mathbb{D}(u^c)$: they are the Stokes matrices $\mathbb{S}_1$, $\mathbb{S}_2$, the formal monodromy exponent $\hbox{\rm diag}(V)$, the exponents of a fundamental matrix solution $Y^{(0)}(z,u)$ with Levelt form at $z=0$ (such as  $\hat{\mu}$ and $R$ described in the example of  Section \ref{dfmanifolds}),  and the central connection matrix $ \mathcal{C}^{(0)}$ such that $Y_1(z,u)=Y^{(0)}(z,u) \mathcal{C}^{(0)}$.   The  Stokes matrices satisfy
$$ 
(\mathbb{S}_1)_{jk}= (\mathbb{S}_1)_{kj}= (\mathbb{S}_2)_{jk}= (\mathbb{S}_2)_{kj}=0 \quad \hbox{ for $j\neq k$ such that } u_j^c=u_k^c. 
$$ 
As a consequence of theorem 1.1. of \cite{CDG}, in order to compute the monodromy  data of system  \eqref{23luglio2021-2}, it   suffices to compute the data of 
\be
\label{27luglio2021-2}
\frac{dY}{dz}=\left(U(u^c)+\frac{V(u^c)}{z}\right) Y,
\ee
which is simpler than \eqref{23luglio2021-2} at generic $u$.\footnote{In particular, for the computation of the Stokes matrices, one has to consider the formal solution of \eqref{27luglio2021-2}  at infinity given by  $Y_F(z,u^c)$, and the corresponding unique canonical solutions  of \eqref{27luglio2021-2} on $\mathcal{S}_1,\mathcal{S}_2,\mathcal{S}_3$ (see  Remark \ref{28luglio2021-1} in the paper).}

Coming back to our  third result, in case $n=3$ and $V(u)$ associated  to   a Painlev\'e transcendent,  we classify in Section \ref{27luglio2021-3} all the branches  $y(x)$ which are  holomorphic 
at $x=0$ and  such that  $V(u)$  satisfies \eqref{5agosto2021-1}  and theorem 1.1 
of \cite{CDG}, for  the coalescence  is $u_2-u_1\to 0$, corresponding  to  $x\to 0$.  
 They form a sub-class of the Taylor expansions  tabulated in \cite{guz2012}. 
 For such transcendents, we show that system \eqref{27luglio2021-2} can be solved in terms of classical special functions, being reducible, depending on the case,  to either the confluent hypergeometric equation or the generalized hypergeometric equation of type $(p,q)=(2,2)$. Our classification and tabulation given in Section \ref{27luglio2021-3} reports for each transcendent   the corresponding  classical special functions.  In Section \ref{20giugno2021-1}, we compute the Stokes matrices and the  invariants $p_{jk}$ for a selection of the classified transcendents.

\vskip 0.2 cm

As a special and important case, in Section  \ref{dfmanifolds}, we compute all the monodromy data (the Stokes matrices, 
the Levelt exponents and the central connection matrix)  in a  {\it chamber} of a 3-dimensional Dubrovin-Frobenius manifold  associated 
with a branch $y(x)$  analytic at $x=0$. In case of a Dubrovin-Frobenius manifold, a transcendent analytic at a one of 
the critical points always satisfies theorem 1.1 of \cite{CDG}, so falling in our classification.  The knowledge of the 
 monodromy data  is very important for the analytic theory of a Dubrovin-Frobenius manifold, since they 
 parameterize the chambers into which the manifold is split (chambers are very similar to local charts,  see  
 \cite{CDG1} for the definition and discussion). By a Riemann-Hilbert boundary value problem, these data allow to  
 construct the manifold structure in the chamber \cite{Dub1,Dub2}. This is a local approach, but the passage from local to global is possible by means of an action of the braid group  \cite{Dub1,Dub2,CDG1}.  See   also \cite{Cotti20} for the analytic Riemann-Hilbert  problem at a coalescent point, and \cite{sabbah} for an algebraic viewpoint.

\vskip 0.2 cm 
In conclusion, the $3\times3$ isomonodromic representation \eqref{23luglio2021-2} of PVI, together with the results of \cite{CDG} on non-generic isomonodromic deformations, implemented in \cite{guz2021}  through an isomonodromic Laplace transform, allow us to: 

\begin{itemize}
\item[-] relate the monodromy invariants $p_{jk}$ to the Stokes matrices, including  the semisimple coalescent case;
\item[-] classify a class of solutions of PVI which are analytic at a critical point, and explicitly compute their monodromy data using classical special functions. 
\item[-] explicitly compute the local monodromy data of 3-dim  Dubrovin-Frobenius manifolds  locally parametrized by Painlev\'e transcendents having  an analytic branch at a critical point.
 
\end{itemize}

%%%%%%%%%%%%%%%%%%%%%%%%%%
\section{ PVI as isomonodromy condition of  an irregular system}
\label{19giugno2021-6}

Theorem \ref{20gennaio2021-2}  below establishes the equivalence between PVI and  the isomonodromy deformation 
equations of a certain 3-dimensional system \eqref{23luglio2021-2}, with explicit formulae.  This result is equivalent to that of \cite{MazzoIrr}, here formulated in the language of two  Pfaffian systems 
(integrable deformations) related by Laplace transform.  
      Consider a $3$-dimensional Frobenius integrable Pfaffian system  
\be
\label{19febbraio2021-11}
 dY=\omega(z,u)Y,\quad\quad \omega(z,u)=\left(U+\frac{V(u)}{z}\right)dz +\sum_{k=1}^3 (zE_k+V_k(u))du_k,
\ee
where $U=\hbox{\rm diag}(u_1,u_2,u_3)$ and $u=(u_1,u_2,u_3)$ varies in a domain of $\mathbb{C}^3$.  
The set of ``diagonals''  
  \be
  \label{19giugno2021-1}
  \Delta_{\mathbb{C}^3}:=\bigcup_{i\neq j}\{u\in\mathbb{C}^3~|~u_i-u_j=0\}
  \ee
  is called {\bf coalescence locus}. 
  In order to work in the analytic domain, we assume that $V$ is holomorphic on a polydisc $\mathbb{D}$, that we can choose in two ways.
  
   \underline{Case 1}. $\mathbb{D}=\mathbb{D}(u^0)$ centered at $u^0$, such that $\mathbb{D}(u^0)\cap  \Delta_{\mathbb{C}^3}=\emptyset$.
   
   \underline{Case 2}. $\mathbb{D}=\mathbb{D}(u^c)$,  such that $  \mathbb{D}(u^c)\cap  \Delta_{\mathbb{C}^3}\neq \emptyset$, 
  with center at $u^c\in  \Delta_{\mathbb{C}^3}$. We assume that $u^c$ is the most coalescent point, namely  if $u_j^c\neq u_k^c$ for some $j\neq k$,  then $u_j\neq u_k$ for all  $u\in \mathbb{D}(u^c)$.  
   There are two possibilities:  either  the case with  two distinct eigenvalues $\lambda_1\neq \lambda_2$, namely  
   $$\lambda_1:=u_i^c=u_j^c \hbox{ for some $1\leq i\neq j\leq 3$,  and  $\lambda_2:=u_k^c\neq u_i^c$ for $k\in\{1,2,3\}\backslash\{i,j\}$},
   $$
    or  the case $u_1^c=u_2^c=u_3^c$. The latter will not be considered here.  We also assume that 
\be
\label{19giugno2021-3}
   \lim_{u_i-u_j\to 0} V_{ij}(u)=0 \hbox{  holomorphically when $u_i-u_j\to 0$ in } \mathbb{D}(u^c),
\ee
so that theorem 1.1 of \cite{CDG} holds.

 \ble
\label{23gennaio2021-7}
  The integrability condition  $d\omega=\omega\wedge \omega$ of \eqref{19febbraio2021-11}   holds  on a domain of $\mathbb{C}^3$ if and only if $V$ and  $V_k$ satisfy on that domain the system\footnote{ The notation in \eqref{20gennaio2021-6} means that $V_k$ is the matrix whose entries $(i,j)$  are the ratios $V_{ij}(\delta_{ik}-\delta_{jk})/(u_i-u_j)$. 
}
  
  \begin{align}
\label{20gennaio2021-6}
&
V_k(u)=\left(
 \frac{V_{ij}(\delta_{ik}-\delta_{jk})}{u_i-u_j}
 \right)_{i,j=1}^3
,\quad 1\leq k  \leq 3.
\\
\label{20gennaio2021-7}
&
\partial_k V=[V,V_k],\quad k=1,2,3.
\end{align}
Moreover, $d\omega=\omega\wedge \omega$ implies the integrability condition  $ \partial_iV_j-\partial_jV_i=V_iV_j-V_jV_i$ of \eqref{20gennaio2021-7} and of the Pfaffian system 
\be
\label{19giugno20201-2}
dG=\Bigl(\sum_{j=1}^3 V_j(u) du_j \Bigr)~G.
\ee
A Jordan form of $V$ is constant  and is given by $G^{-1}VG$ for a suitable  fundamental solution of \eqref{19giugno20201-2}. 
\ele
\begin{proof}
The proof of \eqref{20gennaio2021-6}-\eqref{20gennaio2021-7} and  $ \partial_iV_j-\partial_jV_i=V_iV_j-V_jV_i$ is a computation. The last statement follows from the fact that \eqref{20gennaio2021-7} and \eqref{19giugno20201-2} imply $\partial_j (G^{-1} V G)=0$ for every fundamental solution of \eqref{19giugno20201-2}. For more details, see Sec. 2 of \cite{Guzz-SIGMA}. 
\end{proof}

\bre
\label{19giugno20201-5}
{\rm 
In case $V$ is analytic on either  $\mathbb{D}=\mathbb{D}(u^0)$, or on $\mathbb{D}=\mathbb{D}(u^c)$ with vanishing conditions \eqref{19giugno2021-3}, then a fundamental solution $G(u)$ of \eqref{19giugno20201-2},  such that $G^{-1}VG=J$ is a Jordan form, is holomorphic and holomorphically invertible on $\mathbb{D}$.  For details, we refer to 
 \cite{CDG, Guzz-SIGMA,guz2021} (see also \cite{Eretico1}).
 }
  \ere

\bpr  
\label{23gennaio2021-5} 
Let
$$\Theta=\hbox{\rm diag}(\theta_1,\theta_2,\theta_3):=-\hbox{\rm diag}V(u),\quad\quad x:= \frac{u_2-u_1}{u_3-u_1}.
$$
Then, $\theta_1,\theta_2,\theta_3$ are constant along the solutions of \eqref{20gennaio2021-6}-\eqref{20gennaio2021-7}. Every solution of \eqref{20gennaio2021-6}-\eqref{20gennaio2021-7} admits  the factorization
\be
\label{20gennaio2021-1}
 V(u)=(u_3-u_1)^\Theta~ \Omega(x)~ (u_3-u_1)^{-\Theta}, 
\quad\quad 
V_k(u)= (u_3-u_1)^\Theta ~\Omega_k(u)~ (u_3-u_1)^{-\Theta},
\ee
 where $\Omega(x)$ depends only on $x$, and 
$$
\Omega_1
=
\begin{pmatrix}
0 & \frac{\Omega_{12}(x)}{u_1-u_2} & \frac{\Omega_{13}(x)}{u_1-u_3}
\\ \noalign{\medskip}
\frac{\Omega_{21}(x)}{u_1-u_2} & 0 & 0 
\\ \noalign{\medskip}
\frac{\Omega_{31}(x)}{u_1-u_3} & 0 & 0
\end{pmatrix},
\quad
\Omega_2
=
\begin{pmatrix}
0 & \frac{\Omega_{12}(x)}{u_2-u_1} & 0
\\ \noalign{\medskip}
\frac{\Omega_{21}(x)}{u_2-u_1} & 0 & \frac{\Omega_{23}(x)}{u_2-u_3}
\\ \noalign{\medskip}
0 & \frac{\Omega_{32}(x)}{u_2-u_3} & 0
\end{pmatrix},
\quad
\Omega_3
=
\begin{pmatrix}
0 & 0& \frac{\Omega_{13}(x)}{u_3-u_1}
\\ \noalign{\medskip}
0 & 0 & \frac{\Omega_{23}(x)}{u_3-u_2}
\\ \noalign{\medskip}
\frac{\Omega_{31}(x)}{u_3-u_1} & \frac{\Omega_{32}(x)}{u_3-u_2} & 0
\end{pmatrix}
$$
\epr
\begin{proof} Substituting \eqref{20gennaio2021-6} into \eqref{20gennaio2021-7} we explicitly see that 
$$ 
\partial_k V_{jj}=0 \quad \quad \forall~ j=1,2,3,\quad  \forall~ k=1,2,3.
$$
Notice that $\sum_k V_k=0$. This, and the conditions \eqref{20gennaio2021-6}-\eqref{20gennaio2021-7}, with $\hbox{\rm diag}V(u)=-\Theta$,  imply
$$ 
\sum_{k=1}^3 \partial_k V=0
,\quad\quad 
\sum_{k=1}^3 u_k\partial_k V_{ij}=(\theta_i-\theta_j)V_{ij}.
$$
Thus,  \eqref{20gennaio2021-1} for $V$  follows from  Lemma \ref{21gennaio2021-3} in the  Appendix. Then,   \eqref{20gennaio2021-1}  for $V_k$ follows from  \eqref{20gennaio2021-6}.
\end{proof}

\bpr
\label{23gennaio2021-6}
The three equations $\partial_k V=[V,V_k]$, where $V_k$ is \eqref{20gennaio2021-6},   are equivalent to 
\be
\label{20gennaio2021-26}
\frac{d \Omega}{d x}
=\bigl[\Omega,\widehat{\Omega}_2\bigr],
\ee
where 
$$ 
 \widehat{\Omega}_2(x) :=
\frac{1}{1-x}
\begin{pmatrix}
0 & 0&0
\\ 
0 & 0 & \Omega_{23}(x)
\\ 
0 & \Omega_{32}(x) & 0
\end{pmatrix}- \frac{1}{x}
 \begin{pmatrix}
0 & \Omega_{12}(x) & 0
\\ 
\Omega_{21}(x) & 0 & 0 
\\ 
0& 0 & 0
\end{pmatrix}.
$$

\epr

\begin{proof}  We substitute the factorization \eqref{20gennaio2021-1}  into the three equations  $\partial_k V=[V,V_k]$, $k=1,2,3$,  and then use the chain rule 
$$
\frac{\partial \Omega}{\partial u_j}=\frac{\partial x}{\partial u_j}~\frac{d \Omega}{dx},\quad j=1,2,3.
$$ Since $x=(u_2-u_1)/(u_3-u_1)$, then we obtain the three equations
$$
\frac{d\Omega}{dx}= \frac{(u_2-u_1) [\Omega_1,\Omega]+[\Theta,\Omega]}{x(x-1)},
\quad
\frac{d\Omega}{dx}= \frac{(u_2-u_1) [\Omega_2,\Omega]-[\Theta,\Omega]}{x},
\quad
\frac{d\Omega}{dx}= \frac{(u_2-u_1) [\Omega_3,\Omega]}{-x^2}.
$$
It is a tedious computation to check entry by entry that the above equations are equivalent to \eqref{20gennaio2021-26}.
\end{proof}

The integrability of the Pfaffian system \eqref{19febbraio2021-11} is equivalent to the strong isomonodromy of  its $z$-component 
\be
\label{systempainleve-bis}
\frac{dY}{dz}=\left(U+\frac{V}{z}\right) Y.
\ee
This means that the essential monodormy data (Stokes matrices, monodromy exponents, central connection matrix, see \cite{CDG,Guzz-SIGMA}) are independent of $u$. 
The factorization of $V$ in Proposition \ref{23gennaio2021-5} implies that  {\it in order to compute the monodromy data of \eqref{systempainleve-bis}, it suffices to  assume $u_1=0$, $u_2=x$, $u_3=1$}. With this setting, the $z$-component of  \eqref{19febbraio2021-11} becomes (see also Remark \ref{par-reduction})
\begin{equation}
\label{systempainleve}
\frac{dY}{dz}=\left(U(x)+\frac{\Omega(x)}{z}\right)Y,
\quad\quad
U(x)=
\begin{pmatrix}
0 & 0 & 0\\
0 & x & 0\\
0 & 0 & 1
\end{pmatrix}.
\end{equation}

\vskip 0.2 cm 
 
In the sequel, we will be interested  in the solutions of \eqref{20gennaio2021-6}-\eqref{20gennaio2021-7} satisfying: 
\begin{align}
\label{17giugno2021-4}
&\hbox{\rm diag}V(u)=\hbox{\rm diag}(-\theta_1,-\theta_2,-\theta_3), 
\\
\label{19gennaio2021-12}
&\hbox{\rm $V$  has distinct  eigenvalues}= 0, ~\frac{\theta_\infty - \theta_1 - \theta_2 - \theta_3}{2},~\frac{-\theta_\infty - \theta_1 -\theta_2 - \theta_3}{2}.
\end{align}
$V$ is then diagonalizable. Here, $\theta_\infty$ is just introduced to give a name to the eigenvalues. By Lemma \ref{23gennaio2021-7}, the eigenvalues of $V$ are constant, so that $\theta_\infty$ is a constant. The eigenvalues are distinct if and only if
\be
\label{22giugno2021-2} 
\theta_\infty \neq 0,\pm ( \theta_1 +\theta_2 + \theta_3).
\ee
  As explained in  Remark \ref{22giugno2021-1},  condition \eqref{22giugno2021-2} can always be fulfilled  for a matrix $V$ associated to a  transcendent by Theorem \ref{20gennaio2021-2}.   The following statement is straightforward.

\ble
\label{17giugno2021-10}
If $V$ is a solution of \eqref{20gennaio2021-6}-\eqref{20gennaio2021-7} with constraints \eqref{17giugno2021-4}-\eqref{19gennaio2021-12},  then all the matrices 
$$V^\prime=K^0\cdot  V \cdot (K^0)^{-1}, \quad K^0:=\hbox{\rm diag} (k_1^0,k_2^0,k_3^0),\quad  k_1^0,k_2^0,k_3^0\in\mathbb{C}\backslash\{0\},
$$ 
are solutions with the same constraints.
There is no loss of generality in taking $k_3^0=1$. 
\ele

The main result of the section is 

\bth
\label{20gennaio2021-2}

   The integrability condition     \eqref{20gennaio2021-6}-\eqref{20gennaio2021-7}, with the constraints
      \eqref{17giugno2021-4}-\eqref{19gennaio2021-12},  is equivalent to PVI with 
      coefficients \eqref{11giugno2021-2} given in terms of the parameters $\theta_1,\theta_2,\theta_3,\theta_\infty$ of \eqref{17giugno2021-4}-\eqref{19gennaio2021-12}. 
  Equivalently, the non-linear system \eqref{20gennaio2021-26} with $\Omega$ satisfying the same constraints   \eqref{17giugno2021-4}-\eqref{19gennaio2021-12}  is equivalent to PVI. 
  
  There  is a one-to-one correspondence between transcendents   $y(x)$ and 
  equivalence classes 
 $$
 \Bigl\{K^0\cdot  V \cdot (K^0)^{-1},~K^0=\hbox{\rm diag}  (k_1^0,k_2^0,1),\quad (k_1^0,k_2^0)\in\mathbb{C}^2\backslash\{0,0\}
 \Bigr\}
 $$ 
  of solutions of \eqref{20gennaio2021-7}, 
 or the corresponding  classes $\{K^0\cdot \Omega  \cdot (K^0)^{-1}\}$ of solutions of \eqref{20gennaio2021-26}. The following explicit formulae hold.  \begin{align*}
\Omega_{12}=& \frac{k_1(x)}{k_2(x)}\cdot \dfrac{(x^2 - x)\dfrac{dy}{dx} + (\theta_\infty - 1)y^2 + \Bigl( 
\theta_2 - \theta_1 + 
 1-(\theta_\infty +\theta_2)x\Bigr)y + \theta_1x}{2(x - 1)y},
 \\
 \Omega_{21}=&\frac{k_2(x)}{k_1(x)}\cdot \frac{(x^2 - x)\dfrac{dy}{dx} + (\theta_\infty - 1)y^2 + \Bigl(  \theta_1 - \theta_2 + 1
 -(\theta_\infty - \theta_2)x\Bigr)y - \theta_1x}{2(x - y)},
 \\
 \Omega_{13}=&
 k_1(x)\cdot \frac{(x-x^2)\dfrac{dy}{dx} + (1-\theta_\infty )y^2 + \Bigl((\theta_1 - \theta_3)x + \theta_\infty + \theta_3 - 1\Bigr)y - \theta_1x}{2(x - 1)y},
 \\
 \Omega_{31}=&\frac{1}{k_1(x)}\cdot  
 \frac{(x-x^2)\dfrac{dy}{dx} + ( 1-\theta_\infty)y^2 + \Bigl(( \theta_3-\theta_1)x + \theta_\infty - \theta_3 - 1\Bigr)y + \theta_1x}{2x(y - 1)},
 \\
 \Omega_{23}=& k_2(x)\cdot \frac{(x-x^2 )\dfrac{dy}{dx} +  ( 1-\theta_\infty)y^2 + \Bigl((\theta_\infty - \theta_2)x + \theta_\infty + \theta_3 - 1\Bigr)y - x(\theta_\infty - \theta_2 + \theta_3)}{2(x - y)},
 \\
 \Omega_{32}=&\frac{1}{k_2(x)}\cdot 
 \frac{(x - x^2)\dfrac{dy}{dx} + (  1-\theta_\infty)y^2 + ((\theta_\infty + \theta_2)x + \theta_\infty- \theta_3 - 1)y - x(\theta_\infty + \theta_2 - \theta_3)}{2x(1 - y)},
\end{align*}
and $\hbox{\rm diag} ~\Omega= \hbox{\rm diag} ~V ={\rm diag}(-\theta_1,-\theta_2,-\theta_3)$. 
The functions $k_j(x)$ are obtained by the quadratures 
\be
\label{18gennaio2021-5}
k_j(x)=k_j^0\exp\left\{L_j(x)\right\},\quad  k_j^0\in\mathbb{C}\backslash\{0\},\quad  L_j(x)=\int^x l_j(\xi)d\xi, \quad\quad j=1,2,
\ee
with 

\begin{align*}
l_1(x):=& \frac{x(1-x)\dfrac{dy}{dx} + (\theta_2 -\theta_1 - \theta_3 + 1)y^2 + ((\theta_1 + \theta_3)x + \theta_1 - \theta_2 - 1)y - \theta_1 x}{2x(x - 1)(y - 1)y}
,
\\
\noalign{\medskip}
l_2(x):= &\frac{1}{2x(1-x)(1-y)(x - y)}\left(-x(x - 1)^2\dfrac{dy}{d x} + \Bigl((\theta_1 - 3\theta_2 + \theta_3 + 1)x - \theta_1 + \theta_2 + \theta_3 - 1\Bigr)y^2 + \right.
\\
&
\left.
+ \Bigl(( 2\theta_2-\theta_1  - \theta_3)x^2 + (3\theta_2 - 3\theta_3 - 1)x + \theta_1 - \theta_2 + 1\Bigr)y + \Bigl((\theta_1 - 2\theta_2 + 2\theta_3)x - \theta_1\Bigr)x\right).
\end{align*}
%Conversely, the relation $\Omega\longmapsto y$ is
%\be
%\label{15giugno2021-1}
%y(x)= \frac{xR(x)}{x(1+R(x))-1},\quad R(x):=-\frac{\theta_1+k_1(x)\Omega_{31}(x)}{\theta_3+\Omega_{13}(x)/k_1(x)},
%\ee
%where ....COMPLETARE: VA SPIEGATO CHI SONO k1 e k2 indipendentemente da y
\eth

The proof will be given  in the analytic case $\mathbb{D}=\mathbb{D}(u^0)$, or the  case  $\mathbb{D}=\mathbb{D}(u^c)$ with  conditions \eqref{19giugno2021-3}, but  it is based on  linear algebra  and calculus of derivatives.  Thus, the formulae hold at every point $(u_1,u_2,u_3) $ in a domain of $\mathbb{C}^3$  where the calculations  make sense. They   emend editing  misprints of \cite{MazzoIrr}.  

\bre
\label{22giugno2021-1} 
{\rm 
The condition \eqref{22giugno2021-2} is not restrictive. Indeed, given PVI with coefficients $\alpha,\beta,\gamma,\delta$, the  $\{\theta_\nu\}_{\mu=1,2,3,\infty}$ are defined by  \eqref{11giugno2021-2}, so the changes $\theta_\infty\mapsto 2-\theta_\infty$, $\theta_1\mapsto-\theta_1$, $\theta_2\mapsto-\theta_2$, $\theta_3\mapsto-\theta_3$  do not change the specific equation PVI under consideration.  Therefore, one has to choose the sign of the square roots defining the  $\{\theta_\nu\}_{\mu=1,2,3,\infty}$  in such a way to avoid $\theta_\infty = 0,\pm ( \theta_1 +\theta_2 + \theta_3)$. }
\ere

\bre
\label{11giugno2021-3}
{\rm
In case $\theta_1=\theta_2=\theta_3=0$, then from Theorem \ref{20gennaio2021-2} we compute
$$ 
k_1(x)=k_1^0\frac{\sqrt{y}\sqrt{x-1}}{\sqrt{y-1}\sqrt{x}},\quad\quad 
k_2(x)=k_2^0\frac{\sqrt{y-x}}{\sqrt{y-1}\sqrt{x}}.
$$
If we choose $k_1^0=\pm \sqrt{-1}$, $k_2^0=\pm \sqrt{-1}$, then $$ 
V^T=-V,\quad\Omega^T=-\Omega,
$$
and the matrix $V$  is associated with a  Dubrovin-Frobenius manifold \cite{Dub2}.  The  expressions of Theorem \ref{20gennaio2021-2}  reduce to the formulae  in section 4 of \cite{guzz2001} (see page 269 there for the relation  $y\mapsto V$ and (48) at page 270, for the relation $
V\mapsto y$). These  formulae were  later  used in \cite{guz2011} and in section  22 of   \cite{CDG}.
}
\ere

\subsection{Proof of Theorem \ref{20gennaio2021-2}}

We assume that {\it $V$ is diagonalizable with pairwise distinct eigenvalues,  and that (at least)  one eigenvalue is zero}. Let 
$\mu_1,0,\mu_3$ denote the constant eigenvalues and let
$$ 
\widehat{\mu}:= \hbox{\rm diag}( \mu_1, 0 ,\mu_3)=G(u)^{-1} V(u) G(u),
$$
where $G$ is a fundamental solution of \eqref{19giugno20201-2} (here  Remark \ref{19giugno20201-5} applies).  

\begin{remark}
\label{17giugno2021-1}
{\rm
$G$ is determined up to $G\mapsto G\cdot  \hbox{\rm diag}(\varepsilon_1,\varepsilon_2,\varepsilon_3)$, for $\varepsilon_1,\varepsilon_2,\varepsilon_3\in \mathbb{C}\backslash\{0\}$.  
}
\end{remark}
Notice that  \eqref{19febbraio2021-11}  is integrable if and only if  
\be
\label{20gennaio2021-4}
dY=\omega(z,u)Y,\quad\quad \omega(z,u)=\left(U+\frac{\widetilde{V}}{z}\right)dz +\sum_{k=1}^3 (zE_k+V_k)du_k,\quad\quad  \widetilde{V}:=V-I
\ee
is integrable, because  in both cases the integrability is system   \eqref{20gennaio2021-6}-\eqref{20gennaio2021-7}.

\subsubsection{Step 1: equivalence  between  \eqref{20gennaio2021-6}-\eqref{20gennaio2021-7} and the Schlesinger system \eqref{20gennaio2021-8}}

We consider a  Fuchsian system 
$
(U-\lambda)\frac{d\Psi}{d\lambda} =(I+\widetilde{V}(u))\Psi, 
$
namely, 
\be
\label{20gennaio2021-3}
 \frac{d\Psi}{d\lambda} =\sum_{k=1}^3 \frac{B_k(u)}{\lambda-u_k} \Psi,
 \ee
 where 
 $ B_k=-E_k(\widetilde{V}+I)=-E_k V$. 
 Explicitly, 
 $$ 
 B_1=\begin{pmatrix}
 \theta_1 & -V_{12} & -V_{13}
 \\
 0 & 0 & 0
 \\
 0 & 0 & 0
 \end{pmatrix},
 \quad
 B_2=\begin{pmatrix}
  0 & 0 & 0\\
   -V_{21} &\theta_2 & -V_{23}
 \\
 0 & 0 & 0
 \end{pmatrix},
 \quad
B_3=\begin{pmatrix}
0 & 0 & 0
 \\
 0 & 0 & 0
 \\
  -V_{31} & -V_{32} & \theta_3
 \end{pmatrix}.
 $$
For every  fixed $u$, it is known \cite{BJL4,guz2016} that the relation between \eqref{20gennaio2021-3} and the $z$-component
\be
\label{8maggio2021-1}
 \frac{dY}{dz}=\left(U+\frac{\widetilde{V}}{z}\right) Y 
\ee
 of \eqref{20gennaio2021-4} is given by the Laplace transform
 $$
\vec{Y}(z,u)=\int_\gamma e^{\lambda z}\vec{\Psi}(\lambda,u)d\lambda 
,
$$
   associating  to  a column vector solution $\vec{\Psi}(\lambda,u)$ of \eqref{20gennaio2021-3} a vector solution $\vec{Y}(z,u)$ of  \eqref{8maggio2021-1}, for a suitable path $\gamma$ such that $e^{\lambda z}(\lambda-U)\vPsi(\lambda)\Bigl|_\gamma=0$. 
It is shown in \cite{guz2021}, by elementary computations,  that
the integrability $d\omega=\omega\wedge \omega$   of  system \eqref{20gennaio2021-4} is   equivalent to the integrability $dP=P\wedge P$ of the Pfaffian system 
\be
\label{15ottobre-bloccoGreenPass-1}
d\Psi=P(\lambda,u) \Psi, 
\quad 
\quad
P(\lambda,u)=\sum_{k=1}^3 \left(\frac{B_k(u)}{\lambda-u_k} d(\lambda-u_k) + V_k(u)du_k\right).
\ee
The above is an example of a {\it non-normalized Schlesinger deformation} \cite{Bo}. Explicitly,  $dP=P\wedge P$ is the  non-normalized Schlesinger system 
\begin{equation}
\label{20gennaio2021-8}
\partial_i B_k= \frac{[B_i,B_k]}{u_i-u_k} +[V_i,B_k],\quad i\neq k;
\quad\quad \partial_i B_i=-\sum_{k\neq i} \frac{[B_i,B_k]}{u_i-u_k} +[V_i,B_i]. 
\end{equation}
Therefore,    \eqref{20gennaio2021-6}-\eqref{20gennaio2021-7}  is  equivalent to  \eqref{20gennaio2021-8}.  The equivalence is proved by linear algebra  and calculus of derivatives, so it holds  in  every open domain of $\mathbb{C}^3$  where the derivatives  make sense. 
The equivalence of $d\omega=\omega\wedge \omega$ and $dP=P\wedge P$ is also stated in \cite{MazzoIrr,Boalch1,Boalch2}.  

For further reading on the Laplace transform see  \cite{BJL5,DubDual,Loday,REMY, Scha0, Scha1,Scha2, Scha, Monica}.
   
   \subsubsection{Step 2. From \eqref{20gennaio2021-8} to \eqref{20gennaio2021-12}}
   
   We reduce  \eqref{20gennaio2021-3}  to a $2\times 2$ Fuchsian system. The  associated $2\times 2$ normalized Schlesinger system induced by \eqref{20gennaio2021-8}  is sufficient to find the solutions of  \eqref{20gennaio2021-6}-\eqref{20gennaio2021-7}, up to the freedom $V\mapsto K^0\cdot  V \cdot (K^0)^{-1}$ of Lemma \ref{17giugno2021-10}.   
We start by rewriting \eqref{20gennaio2021-3} as 
\be
\label{20gennaio2021-5}
\frac{dX}{dz}=\sum_{k=1}^3\frac{\widetilde{B}_k }{\lambda-u_k}X ;\quad\quad\quad
\widetilde{B}_k :=-G^{-1}E_kG\widehat{\mu},~\quad X=G^{-1} \Psi.
\ee
Explicitly, 
\be
\label{20gennaio2021-16}
\widetilde{B}_1=\frac{1}{\det G}\begin{pmatrix}
-(G_{22} G_{33}-G_{23} G_{32}) b_1 & 0  & -(G_{22} G_{33}-G_{23} G_{32}) d_1
\\
(G_{21} G_{33}-G_{23} G_{31}) b_1 &   0 &   (G_{21} G_{33}-G_{23} G_{31})d_1 
\\
-(G_{21} G_{32}-G_{22} G_{31}) b_1 & 0  &  -(G_{21} G_{32}-G_{22} G_{31}) d_1
\end{pmatrix}
\ee
\be
\label{20gennaio2021-17}
\widetilde{B}_2=\frac{1}{\det G}\begin{pmatrix}
(G_{12} G_{33}-G_{13} G_{32}) b_2 & 0 & (G_{12} G_{33}-G_{13} G_{32}) d_2
\\
-(G_{11} G_{33}-G_{31} G_{13}) b_2 &0&-(G_{11} G_{33}-G_{31} G_{13}) d_2
\\
(G_{11} G_{32}-G_{31} G_{12}) b_2  &0 &(G_{11} G_{32}-G_{31} G_{12}) d_2
\end{pmatrix}
\ee
\be
\label{20gennaio2021-18}
\widetilde{B}_3=\frac{1}{\det G}\begin{pmatrix}
-(G_{12} G_{23}-G_{13} G_{22}) b_3&0 &-(G_{12} G_{23}-G_{13} G_{22}) d_3
\\
(G_{11 }G_{23}-G_{21} G_{13}) b_3&0 &(G_{11} G_{23}-G_{21}G_{13}) d_3
\\
-(G_{11 }G_{22}-G_{21} G_{12}) b_3&0&-(G_{11} G_{22}-G_{21} G_{12}) d_3
\end{pmatrix},
\ee
where 
\be
\label{20gennaio2021-19}
b_1=G_{11}\mu_1,\quad  b_2= G_{21} \mu_1, \quad b_3= G_{31} \mu_1,
\quad
d_1=G_{13}\mu_3,\quad  d_2= G_{23} \mu_3, \quad d_3= G_{33} \mu_3,
\ee
By construction,
$$ 
\hbox{\rm Tr} \widetilde{B}_k=\hbox{\rm Tr} {B}_k=\theta_k;\quad \quad \hbox{eigenvalues of $B_k$ and $\widetilde{B}_k$  } = \theta_k,~0,~0.
$$
The non-normalized Schlesinger system \eqref{20gennaio2021-8} is equivalent to the normalized one 
\begin{equation}
\label{20gennaio2021-10}
\partial_i \widetilde{B}_k= \frac{[\widetilde{B}_i,\widetilde{B}_k]}{u_i-u_k},\quad i\neq k ;
\quad\quad
\partial_i \widetilde{B}_i=-\sum_{k\neq i} \frac{[\widetilde{B}_i,\widetilde{B}_k]}{u_i-u_k} .
\end{equation}
Let a vector solution of  \eqref{20gennaio2021-5}   be denoted by 
$$\vec{X}:=
 \begin{pmatrix}
 X_1\\ X_2\\ X_3
 \end{pmatrix}, 
 \quad 
 \quad \hbox{ and let }\quad  
 \boldsymbol{X}:=
 \begin{pmatrix}
 X_1\\ X_3
 \end{pmatrix}.
 $$
Then,  $\boldsymbol{X}$ satisfies 
 \be
 \label{20gennaio2021-22}
 \frac{d \boldsymbol{X}}{d\lambda} =\sum_{k=1}^3 \frac{A_k}{\lambda-u_k}\boldsymbol{X}, \quad
  \quad
   A_k:=\begin{pmatrix}
   (\widetilde{B}_k)_{11} &  (\widetilde{B}_k)_{13} 
   \\
    (\widetilde{B}_k)_{31}  &  (\widetilde{B}_k)_{33} 
   \end{pmatrix},
   \ee
  while $X_2$ is obtained by a quadrature
   $$ 
   X_2(\lambda)=\sum_{k=1}^3 \int^\lambda \frac{(B_k)_{21}X_1(s)+(B_k)_{23}X_3(s)}{s-u_k}ds.
   $$
In the integration, $\lambda$ lies  in a plane with branch cuts issuing from $u_1,u_2,u_3$ (see Section \ref{25gennaio2021-1} for details). 
   By construction 
   \be
   \label{13giugno2021-1} 
   A_1+A_2+A_3= -\begin{pmatrix} \mu_1 & 0 \\ 0 & \mu_3\end{pmatrix}, 
   \quad \quad \hbox{\rm Tr} A_k=\theta_k,\quad \quad \hbox{Eigenvalues of $A_k$} = \theta_k,0.
\ee
The Schlesinger system \eqref{20gennaio2021-10} implies 
 \be
\label{20gennaio2021-12}
\partial_i A_k= \frac{[A_i,A_k]}{u_i-u_k} ,\quad i\neq k;\quad\quad
\partial_i A_i=-\sum_{k\neq i} \frac{[A_i,A_k]}{u_i-u_k} .
\ee

The following straightforward lemma  is  the analogue of Lemma \ref{17giugno2021-10}.

\ble
\label{17giugno2021-14-bis}
If $\{A_k\}_{k=1,2,3}$ is a solution of \eqref{20gennaio2021-12} with constraints \eqref{13giugno2021-1}, then for every $\varepsilon_1,\varepsilon_3\in\mathbb{C}\backslash\{0\}$   the matrices
$
\hbox{\rm diag}(\varepsilon_1,\varepsilon_3)^{-1}\cdot A_k\cdot \hbox{\rm diag}(\varepsilon_1,\varepsilon_3) 
$ 
are a solution with the same constraints \eqref{13giugno2021-1}. 
\ele

\ble
\label{17giugno2021-14}
There is a one-to-one  correspondence between  equivalent classes 
\be
\label{17giugno2021-12}
\{\mathcal{E}^{-1}A_1 \mathcal{E},~\mathcal{E}^{-1}A_2 \mathcal{E},~\mathcal{E}^{-1}A_3 
\mathcal{E},\quad\mathcal{E}:=\hbox{\rm diag}(\varepsilon_1,\varepsilon_3),
\quad
 \varepsilon_1,\varepsilon_3\in\mathbb{C}\backslash\{0\}\}
\ee
of solutions of  \eqref{20gennaio2021-12}, satisfying the constraint \eqref{13giugno2021-1}, and  equivalence classes  
 \be
 \label{18giugno2021-1}
 \{K^0\cdot  V \cdot (K^0)^{-1},\quad K^0=\hbox{\rm diag}  (k_1^0,k_2^0,k_3^0),\quad k_1^0,k_2^0,k_3^0\in\mathbb{C}\backslash\{0\} \}
 \ee
   of solutions of \eqref{20gennaio2021-6}-\eqref{20gennaio2021-7} satisfying the constraints \eqref{17giugno2021-4}-\eqref{19gennaio2021-12}.

\ele

\begin{proof}  We propose a proof slightly different from \cite{MazzoIrr}. To every solution $V$ of \eqref{20gennaio2021-6}-\eqref{20gennaio2021-7} satisfying the constraints \eqref{17giugno2021-4}-\eqref{19gennaio2021-12},  we associate $A_k$ in \eqref{20gennaio2021-22}, using  \eqref{20gennaio2021-16}-\eqref{20gennaio2021-19}, so that 
\be
\label{17giugno2021-8}
A_k =
\begin{pmatrix}
a_kb_k & a_k d_k
\\
c_k b_k & c_k d_k
\end{pmatrix},\quad k=1,2,3,
\ee
where
\be
\label{17giugno2021-6}
a_1= \frac{-(G_{22} G_{33}-G_{23} G_{32})}{\det G},\quad a_2=\frac{G_{12} G_{33}-G_{13} G_{32}}{\det G} ,
\quad a_3=\frac{-(G_{12} G_{23}-G_{13} G_{22})}{\det G}.
\ee
\be
\label{17giugno2021-7}
c_1=\frac{-(G_{21} G_{32}-G_{22} G_{31})}{\det G},\quad c_2=\frac{G_{11} G_{32}-G_{31} G_{12}}{\det G} ,
\quad c_3=\frac{-(G_{11 }G_{22}-G_{21} G_{12}) }{\det G}
\ee
The above expressions determine  a class \eqref{17giugno2021-12} in terms of $V$. Indeed, V determines $G$ up to $G\mapsto G\cdot \hbox{\rm diag}(\varepsilon_1,\varepsilon_2,\varepsilon_3)$, so that there is the freedom   
$$\det G\mapsto \varepsilon_1\varepsilon_2\varepsilon_3\det G, \quad
a_k \mapsto \varepsilon_1^{-1}a_k, \quad  b_k\mapsto  \varepsilon_1b_k,
\quad 
 c_k\mapsto \varepsilon_3 c_k, \quad d_k\mapsto \varepsilon_3 d_k,
 $$
which determines the $A_k$ up to 
\be
\label{17giugno2021-2} 
A_k \longmapsto 
\hbox{\rm diag}(\varepsilon_1,\varepsilon_3)^{-1}\cdot A_k\cdot \hbox{\rm diag}(\varepsilon_1,\varepsilon_3) 
\ee
Moreover, if we change $V\mapsto V^\prime =K^0\cdot  V \cdot (K^0)^{-1}$ as in  Lemma \ref{17giugno2021-10}, to  $G$ such that $G^{-1}VG=\widehat{\mu}$ there corresponds 
$$ 
G^\prime=\hbox{\rm diag} (k_1^0,k_2^0,k_3^0)\cdot G,\quad \quad (G^\prime)^{-1}V^\prime  G^\prime=\widehat{\mu}.
$$
Therefore, from \eqref{20gennaio2021-19} and \eqref{17giugno2021-6}-\eqref{17giugno2021-7}
$$
a_j^\prime= a_j/k_j^0, \quad b_j^\prime=k_j^0 b_j, \quad c^\prime_j=c_j/k_j^0, \quad d_j^\prime=k_j^0 d_j
.
$$   It follows that  $A_k$ in \eqref{17giugno2021-8} is invariant under the change $V\mapsto V^\prime$, namely it is associated to the equivalence class of $V$. 

Conversely, consider   a solution of   \eqref{20gennaio2021-12} with constraint \eqref{13giugno2021-1}, which can be written as   in \eqref{17giugno2021-8}.  We find the entries of $V=G\widehat{\mu}G^{-1}$ starting from $a_j$, $b_j$, $c_j$, $d_j$ ($1\leq j \leq 3$)  as follows
\be
\label{20gennaio2021-23}
V= \begin{pmatrix}
-\theta_1 & -a_2b_1-c_2d_1 & -a_3 b_1 -c_3 d_1
\\
-a_1b_2-c_1d_2 & -\theta_2 & -a_3b_2-c_3d_2
\\
-a_1b_3-c_1d_3 & -a_2b_3-c_2d_3 & -\theta_3
\end{pmatrix}
.
\ee
The above is invariant under the map  $
a_k \mapsto \varepsilon_1^{-1}a_k$, $ b_k\mapsto  \varepsilon_1b_k$, $ c_k\mapsto \varepsilon_3 c_k$, $d_k\mapsto \varepsilon_3 d_k $ and the map \eqref{17giugno2021-2}. 
Moreover,  $A_j$   determines $a_j$, $b_j$, $c_j$, $d_j$ up to $a_j\mapsto  a_j/k_j^0$, $b_j \mapsto k_j^0 b_j$, $c_j \mapsto c_j/k_j^0$, $d_j \mapsto k_j^0 d_j$, $j=1,2,3$,  so that to each $\{A_1,A_2,A_3\}$ we can only associate  $V$ up to the freedom $V\mapsto K^0\cdot  V \cdot (K^0)^{-1}$.  \end{proof}

%%%%%%%
\subsubsection{Step 3. From \eqref{20gennaio2021-12} to  \eqref{23gennaio2021-1}} 
The third step is the equivalence between the Schlesinger equations \eqref{20gennaio2021-12}   and the equations \eqref{23gennaio2021-1} below  with  only  one independent variable $x$. The  solutions of \eqref{23gennaio2021-1}  in terms of scalar functions of $x$, and the  corresponding  $\Omega(x)$, will be given in Proposition \ref{13giugno2021-3}. 
We start with the gauge 
$$
\boldsymbol{X}=\prod_{j=1}^3 (\lambda-u_j)^{\theta_j/2} \Phi
$$
applied to \eqref{20gennaio2021-22}. It  yields 
\be
\label{20gennaio2021-21}
 \frac{d \Phi}{d\lambda} =\sum_{k=1}^3 \frac{\mathcal{A}_k}{\lambda-u_k} \Phi, \quad
  \quad
   \mathcal{A}_k:=A_k-\frac{\theta_k}{2}I.
   \ee
   Let $ \mathcal{A}_\infty:=-(\mathcal{A}_1+\mathcal{A}_2+\mathcal{A}_3)$. Then, the equalities \eqref{13giugno2021-1} become
   \be
   \label{20gennaio2021-14}
   \mathcal{A}_\infty=
   \begin{pmatrix}
   \mu_1+ \sum_{j=1}^2{\theta_j/2} & 0
     \\
      0 &   \mu_3+ \sum_{j=1}^2{\theta_j/2}    
      \end{pmatrix},\quad
 \quad \hbox{eigenvalues of $\mathcal{A}_k$}= \pm \frac{\theta_k}{2}.
   \ee
In particular,    $\hbox{Tr}\mathcal{A}_k=0$.   Since $  \hbox{Tr}  \mathcal{A}_\infty=0$, we define $\theta_\infty$  by setting 
 $$
   \mathcal{A}_\infty=:
   \begin{pmatrix}
    \theta_\infty/2& 0
     \\
      0 & -\theta_\infty/2
   \end{pmatrix},
$$
so that 
\be
   \label{5febbraio2021-1} 
  \mu_1= \frac{\theta_\infty - \theta_1 - \theta_2 - \theta_3}{2},\quad \mu_3=\frac{-\theta_\infty - \theta_1 -\theta_2 - \theta_3}{2}.
   \ee
   
   \ble
Let $z_k$ and $v_k$ be defined by $z_k+\theta_k=a_kb_k$ and $v_kz_k=-a_kd_k$, for $k=1,2,3$. 
Then 
\be
\label{22gennaio2021-3}
\mathcal{A}_k= \begin{pmatrix}
z_k+\dfrac{\theta_k}{2} & -v_k z_k
\\
\noalign{\medskip}
\dfrac{z_k+\theta_k}{v_k} & -z_k-\dfrac{\theta_k}{2}
\end{pmatrix};
\quad\quad
V_{ij}= \left(z_j-\frac{v_jz_j}{v_iz_i}(z_i+\theta_i)\right)\frac{d_i}{d_j}, \quad 1\leq i\neq j\leq 3.
\ee
 \ele

     \begin{proof}
   The conditions  \eqref{20gennaio2021-14} imply that $c_k=(\theta_k-a_kb_k)/d_k$,  so that the matrix $\mathcal{A}_k$ from  \eqref{17giugno2021-8}  and  \eqref{20gennaio2021-21} is 
   \be
   \label{20gennaio2021-25}
   \mathcal{A}_k= \begin{pmatrix}
   a_kb_k -{\theta_k/2} & a_kd_k
   \\
   \noalign{\medskip}
   {(\theta_k-a_kb_k)b_k/d_k}	& {\theta_k/2}-a_kb_k
   \end{pmatrix}, 
\ee
 We substitute into \eqref{20gennaio2021-23} and  obtain 
 $$
   V_{ij}= -\left( a_jb_i+\frac{d_i}{d_j}(\theta_j-a_jb_j)\right),\quad 1\leq i \neq j \leq 3.
$$
   Then, \eqref{22gennaio2021-3} follows by an elementary  computation.
\end{proof}
 
 \bre
 {\rm 
 Matrices \eqref{20gennaio2021-25} and \eqref{20gennaio2021-23} (with $c_k=(\theta_k-a_kb_k)/d_k$) are the analogous of respectively (3.55) and (3.62) in \cite{Hard}. With the coordinates  $v_k$, $z_k$, (3.55) and (3.62) of \cite{Hard} correspond to  \eqref{22gennaio2021-3}. 
 }
 \ere
 
 It follows from \eqref{20gennaio2021-12} that also the matrices $\mathcal{A}_k$ satisfy the Schlesinger equations 
  \be
\label{21gennaio2021-1}
\partial_i \mathcal{A}_k= \frac{[\mathcal{A}_i,\mathcal{A}_k]}{u_i-u_k} ,\quad i\neq k;
\quad\quad
\partial_i \mathcal{A}_i=-\sum_{k\neq i} \frac{[\mathcal{A}_i,\mathcal{A}_k]}{u_i-u_k} ,
\quad \hbox{with constraints \eqref{20gennaio2021-14}}
\ee
\ble
 Every solution of \eqref{21gennaio2021-1} admits the factorization 
\be
\label{22gennaio2021-1}
\mathcal{A}_k(u)= (u_3-u_1)^{-\mathcal{A}_\infty} \widetilde{\mathcal{A}}_k(x) ~(u_3-u_1)^{\mathcal{A}_\infty},\quad k=1,2,3,\quad \quad x:=\frac{u_2-u_1}{u_3-u_1}.
\ee
where $\widetilde{\mathcal{A}}_k(x)$ solves the Schlesinger equations 
\be
  \label{23gennaio2021-1}
  \frac{d\widetilde{\mathcal{A}}_1}{dx}=\frac{[\widetilde{\mathcal{A}}_2,\widetilde{\mathcal{A}}_1]}{x}, \quad
  \quad
   \frac{d\widetilde{\mathcal{A}}_3}{dx}=\frac{[\widetilde{\mathcal{A}}_2,\widetilde{\mathcal{A}}_3]}{x-1},
  \quad
  \quad
   \frac{d\widetilde{\mathcal{A}}_2}{dx}=-\frac{[\widetilde{\mathcal{A}}_2,\widetilde{\mathcal{A}}_1]}{x}-\frac{[\widetilde{\mathcal{A}}_2,\widetilde{\mathcal{A}}_3]}{x-1}.
  \ee
with constraints
\be
\label{13giugno2021-2} 
\hbox{\rm  eigenvalues of  $\widetilde{\mathcal{A}}_k$ }=\pm \frac{\theta_k}{2},
\quad\quad
 \widetilde{\mathcal{A}}_1+\widetilde{\mathcal{A}}_2+\widetilde{\mathcal{A}}_3=-\mathcal{A}_\infty=\hbox{\rm diag}(-\theta_\infty/2,\theta_\infty/2),
\ee

\ele
\begin{proof}
The equations \eqref{21gennaio2021-1} imply that\footnote{ They also imply  
$ 
\partial_i(\sum_{k=1}^3\mathcal{A}_k)=0
$, for every $i=1,2,3$, 
in accordance with the fact that $\sum_{k=1}^3\mathcal{A}_k=-\mathcal{A}_\infty$ constant.} 
$$
\sum_{i=1}^3 \partial_i \mathcal{A}_k=0 ;\quad\quad  \sum_{i=1}^3 u_i\partial_i \mathcal{A}_k =\bigl[ \mathcal{A}_k, \mathcal{A}_\infty\bigr],\quad k=1,2,3.
$$
In particular, 
$$ 
 \sum_{i=1}^3 u_i\partial_i (\mathcal{A}_k)_{12}=-\theta_\infty (\mathcal{A}_k)_{12},
 \quad 
 \quad
 \sum_{i=1}^3 u_i\partial_i (\mathcal{A}_k)_{21}=\theta_\infty (\mathcal{A}_k)_{21},
 \quad
 \quad 
 \sum_{i=1}^3 u_i\partial_i (\mathcal{A}_k)_{jj}=0 .
 $$
 Therefore, form Lemma \ref{21gennaio2021-3}, 
 $$  (\mathcal{A}_k(u))_{12}=(u_3-u_1)^{-\theta_\infty} (\widetilde{\mathcal{A}}_k(x))_{12}
 ,\quad(\mathcal{A}_k(u))_{21}=(u_3-u_1)^{\theta_\infty} (\widetilde{\mathcal{A}}_k(x))_{21}
 ,\quad(\mathcal{A}_k(u))_{jj}=(\widetilde{\mathcal{A}}_k(x))_{jj} .
 $$
 This proves \eqref{22gennaio2021-1}.  It follows from  \eqref{22gennaio2021-1} 
  and   the chain rule 
    that  \eqref{21gennaio2021-1} is equivalent to \eqref{23gennaio2021-1} with the given constraints.
 \end{proof}

 \bpr 
 \label{13giugno2021-3}
  Let $ \widetilde{\mathcal{A}}_k(x)$ be the matrix in the factorization \eqref{22gennaio2021-1}. Then, 
 \be
 \label{22gennaio2021-8}
 \widetilde{\mathcal{A}}_k(x)= \begin{pmatrix}
z_k(x)+\dfrac{\theta_k}{2} &~ -\widetilde{v}_k(x) z_k(x)
\\
\noalign{\medskip}
\dfrac{z_k(x)+\theta_k}{\widetilde{v}_k(x)} &~ -z_k(x)-\dfrac{\theta_k}{2}
\end{pmatrix},\quad k=1,2,3,
\ee
holds for some functions  $\widetilde{v}_k(x)$, where $z_k=z_k(x)$ is the same as in \eqref{22gennaio2021-3} and depends only on $x$.  
 Moreover, the matrix $\Omega(x)$ in the factorization \eqref{20gennaio2021-1} can be written as 
 \be
 \label{22gennaio2021-10}
 \Omega(x)=K(x) \boldsymbol{Z}(x) K(x)^{-1},\quad \quad K(x)=\hbox{\rm diag}(k_1(x),k_2(x),1)
 \ee
 for some scalar functions $k_1(x)$, $k_2(x)$, where 
 \begin{align}
 \label{22gennaio2021-5}
  &\boldsymbol{Z}_{ii}=-\theta_i, & i=1,2,3,
  \\
  \label{22gennaio2021-6}
 &  \boldsymbol{Z}_{ij}(x)=z_j(x)-\frac{\widetilde{v}_j(x)z_j(x)}{\widetilde{v}_i(x)z_i(x)}\bigl(z_i(x)+\theta_i\bigr),& 1\leq i \neq j \leq 3.
 \end{align}
 \epr 
 
 \begin{proof}
 It follows from the factorization \eqref{20gennaio2021-1} and $G^{-1}VG=\hbox{\rm diag}(\mu_1,0,\mu_3)$  that
 $$
 \hbox{\rm diag}(\mu_1,0,\mu_3)= \left[(u_3-u_1)^{-\Theta} G(u) \right]^{-1}\Omega(x)~\left[ (u_3-u_1)^{-\Theta} G(u)\right].
 $$
  The above implies, for some functions $\delta_j(u)$, the factorization
 $$
  G(u)= (u_3-u_1)^{\Theta}\cdot \widetilde{G}(x) \cdot \hbox{\rm diag}
 \bigl( \delta_1(u),\delta_2(u) ,\delta_3(u)\bigr)
  ,
  $$ where $ \widetilde{G}(x)$ diagonalizes $\Omega(x)$. 
    From the above, we receive the following factorization
\begin{align*}
&  a_j(u)=\widetilde{a}_j(x) \frac{(u_3-u_1)^{-\theta_j} }{\delta_1(u)}, \quad\quad b_j(u)=\widetilde{b}_j(x) (u_3-u_1)^{-\theta_j} \delta_1(u)
   ,  
   \\
   &
 c_j(u)=\widetilde{c}_j(x) \frac{(u_3-u_1)^{-\theta_j} }{\delta_3(u)}
,\quad\quad   d_j(u)=\widetilde{d}_j(x) (u_3-u_1)^{\theta_j} \delta_3(u),\quad\quad j=1,2,3.
  \end{align*}
  Hence, \eqref{20gennaio2021-25} becomes
  $$\mathcal{A}_k =
   \begin{pmatrix}
  \widetilde{ a}_k(x)\widetilde{b}_k(x) -\dfrac{\theta_k}{2} 
  &
   \widetilde{a}_k(x)\widetilde{d}_k(x) \cdot \dfrac{\delta_3(u)}{\delta_1(u)}
   \\
   \noalign{\medskip}
  \dfrac{\delta_1(u)}{\delta_3(u)}\cdot   \dfrac{(\theta_k-\widetilde{a}_k(x)\widetilde{b}_k(x))\widetilde{b}_k(x)}{\widetilde{d}_k(x)}	& \dfrac{\theta_k}{2}-\widetilde{ a}_k(x)\widetilde{b}_k(x)
   \end{pmatrix}
   $$
Comparison with \eqref{22gennaio2021-1} yields 
  $$ 
  \frac{\delta_1(u)}{\delta_3(u)}=(u_3-u_1)^{\theta_\infty},\quad\quad 
 \widetilde{ \mathcal{A}}_k= \begin{pmatrix}
  \widetilde{ a}_k(x)\widetilde{b}_k(x) -{\theta_k}/{2} & \widetilde{a}_k(x)\widetilde{d}_k(x) 
   \\
   \noalign{\medskip}
   {\bigl(\theta_k-\widetilde{a}_k(x)\widetilde{b}_k(x)\bigr)\widetilde{b}_k(x)}/{\widetilde{d}_k(x)}	& {\theta_k}/{2}-\widetilde{ a}_k(x)\widetilde{b}_k(x)
   \end{pmatrix}. 
  $$
  Comparison with $\mathcal{A}_k$ in  \eqref{22gennaio2021-3} shows that 
  $$ 
  z_k(u)\equiv z_k(x)\quad \Longrightarrow \quad v_k(u)=\widetilde{v}_k(x) \dfrac{\delta_3(u)}{\delta_1(u)}
  \quad\hbox{ and } \quad \widetilde{\mathcal{A}}_k= \begin{pmatrix}
z_k(x)+{\theta_k}/{2} & -\widetilde{v}_k(x) z_k(x)
\\
\noalign{\medskip}
{(z_k(x)+\theta_k)}/{\widetilde{v}_k(x)} & -z_k(x)-{\theta_k}/{2}
\end{pmatrix}.
  $$
In conclusion $V$ in  \eqref{22gennaio2021-3} becomes
$$ 
V_{ij}(u)= \left(z_j(x)-\frac{\widetilde{v}_j(x)z_j(x)}{\widetilde{v}_i(x)z_i(x)}\bigl(z_i(x)+\theta_i\bigr)\right)\frac{\widetilde{d}_i(x)}{\widetilde{d}_j(x)}(u_3-u_1)^{\theta_i-\theta_j}, \quad 1\leq i\neq j\leq 3.
$$
The above expression again proves the factorization \eqref{20gennaio2021-1}, and shows that 
  $$
  \Omega_{ij}(x)= \left(z_j(x)-\frac{\widetilde{v}_j(x)z_j(x)}{\widetilde{v}_i(x)z_i(x)}\bigl(z_i(x)+\theta_i\bigr)\right)\frac{\widetilde{d}_i(x)}{\widetilde{d}_j(x)}, \quad 1\leq i\neq j\leq 3.
  $$
  Therefore, 
  \be
  \label{22gennaio2021-7}
  \Omega(x)=\hbox{\rm diag}\Bigl(\widetilde{d}_1(x),\widetilde{d}_2(x),\widetilde{d}_3(x)\Bigr) \cdot \boldsymbol{Z}(x) \cdot \hbox{\rm diag}\Bigl(\widetilde{d}_1(x),\widetilde{d}_2(x),\widetilde{d}_3(x)\Bigr)^{-1},
  \ee
  where $\boldsymbol{Z}(x)$  has entries as in \eqref{22gennaio2021-5}-\eqref{22gennaio2021-6}. Now, \eqref{22gennaio2021-7} coincides with \eqref{22gennaio2021-10}, with 
  $ 
    k_1(x):=\widetilde{d}_1(x)/ \widetilde{d}_3(x)$ and  $k_2(x):= \widetilde{d}_2(x)/\widetilde{d}_3(x)$. 
  \end{proof}

  \subsubsection{Step 4. From \eqref{23gennaio2021-1} to PVI}
  The last step is the one-to-one correspondence  between equivalence classes 
\be
\label{17giugno2021-15}
  \{\mathcal{E}^{-1}\widetilde{\mathcal{A}}_1 \mathcal{E},~\mathcal{E}^{-1}\widetilde{\mathcal{A}}_2 \mathcal{E},~\mathcal{E}^{-1}\widetilde{\mathcal{A}}_3 
\mathcal{E},\quad\mathcal{E}=\hbox{\rm diag}(\varepsilon_1,\varepsilon_3),
\quad
 \varepsilon_1,\varepsilon_3\in\mathbb{C}\backslash\{0\}\}
 \ee
 of solutions of the Schlesinger equations  \eqref{23gennaio2021-1}   with constraints  \eqref{13giugno2021-2}, 
 corresponding to the classes of Lemma \ref{17giugno2021-14},  and  solutions of PVI. This equivalence has been known 
 since the work of R. Fuchs \cite{Fuchs} and L. Schlesinger \cite{Schles}. We will use its  formulation as in appendix C of \cite{JM}. 
  
 Formula (C.47) of  \cite{JM} provides the  parameterization of the matrices satisfying the conditions   \eqref{13giugno2021-2}, which   coincides with our  \eqref{22gennaio2021-8}.  
 Our  $x$ coincides with $t$ used in \cite{JM}. The matrices $A_0(t),A_t(t),A_1(t)$ in  (C.47) are related to ours by the identifications 
 $$ 
 A_0(t)-\frac{\theta_0}{2}=\widetilde{\mathcal{A}}_1(x),\quad  A_t(t)-\frac{\theta_t}{2}=\widetilde{\mathcal{A}}_2(x),\quad  A_1(t)-\frac{\theta_1}{2}=\widetilde{\mathcal{A}}_3(x),\quad t=x,
 $$ 
 $$
 (\theta_0,\theta_t,\theta_1) \hbox{ in (C.47) } =  (\theta_1,\theta_2,\theta_3)  \hbox{ in our notations.} 
 $$
 $$
 (z_0,z_t,z_1)\hbox{ in (C.47) } = (z_1,z_2,z_3)  \hbox{ in our notations.}
 $$ 
 $$ 
 (u, w, v) \hbox{ in (C.47) } = (\widetilde{v}_1,\widetilde{v}_2,\widetilde{v}_3)  \hbox{ in our notations.}
 $$
 
 The condition  $\widetilde{\mathcal{A}}_1+\widetilde{\mathcal{A}}_2+\widetilde{\mathcal{A}}_3=-\mathcal{A}_\infty$ becomes 
 $$
 z_1+z_2+z_3+\frac{\theta_1+\theta_2+\theta_3+\theta_\infty}{2}=0,
 \quad
\frac{z_1+\theta_1}{\widetilde{v}_1}
+
 \frac{z_2+\theta_2}{\widetilde{v}_2}
+\frac{z_3+\theta_3}{\widetilde{v}_3}=0,
\quad 
\widetilde{v}_1z_1+\widetilde{v}_2z_2+\widetilde{v}_3z_3=0.
$$
It can be verified, as in \cite{JM},    that  the matrices $\widetilde{\mathcal{A}}_k$ are parameterized by  the 7+1   independent parameters 
 $$ 
 \theta_1, \quad \theta_2, \quad \theta_3,\quad \theta_\infty, \quad k,\quad y,\quad z,\quad x,
 $$ 
 where $k$, $y$ and $z$  are respectively defined by 
 $$ 
 k=(1+x)\widetilde{v}_1z_1+\widetilde{v}_2z_2+x\widetilde{v}_3z_3,\quad \quad 
  y=\frac{x\widetilde{v}_1z_1}{k},
  $$
  $$
   z=\frac{z_1+\theta_1}{y}+\frac{z_2+\theta_2}{y-x}+\frac{z_3+\theta_3}{y-1}\equiv \frac{(\widetilde{\mathcal{A}}_1)_{11}+\theta_1/2}{y}+\frac{(\widetilde{\mathcal{A}}_2)_{11}+\theta_2/2}{y-x}+\frac{(\widetilde{\mathcal{A}}_3)_{11}+\theta_3/2}{y-1}.  
   $$
   Here, do not confuse the parameter $z$ -- whose symbol we borrow from \cite{JM} --  with the independent variable $z$ in  system \eqref{19febbraio2021-11}. 
The explicit parameterization of $\widetilde{v}_k$ and $z_k$, $k=1,2,3$ in terms of $ 
 \theta_1, \theta_2,  \theta_3, \theta_\infty,  k,y, z, x$ is in formulae (C.51)-(C.52) of \cite{JM}, where $\tilde{z}=z-\theta_1/y-\theta_2/(y-x)  -\theta_3/(y-1)$ is used in place of $z$. 
  
The Schlesinger equations \eqref{23gennaio2021-1} become a first order non-linear system of three differential equations for $y, z, k$, reported in formula (C.55) of \cite{JM}. 
 As for $k(x)$, it  is computable as the exponential of a  quadrature in $dx$ involving $y(x)$ and $z(x)$, so it is determined  up to a multiplicative constant $k^0$, which is identified with  $\varepsilon_3/\varepsilon_1$ in the equivalence class \eqref{17giugno2021-15} in  Lemma \ref{17giugno2021-14-bis} and  \ref{17giugno2021-14}.
 Eliminating $z$ from the remaining first order system for $y$ and $z$, we see that   $y$ solves   PVI. 
If $y=y(x)$  is a solution, then  
\be
\label{23gennaio2021-3}
z(x)= \frac{1}{2}\left(\frac{x(x-1)}{y(x)(y(x)-1)(y(x)-x)} \frac{dy(x)}{dx}
+\frac{\theta_1}{y(x)}+\frac{\theta_2-1}{y(x)-x}+\frac{\theta_3}{y(x)-1}\right).
\ee

 \bre{\rm 
   If we express $(\widetilde{v}_1,\widetilde{v}_2,\widetilde{v}_3)$, $(z_1,z_2,z_3)$ in terms of $(k,y,z)$ we see that 
   {\small
   $$ 
   \left(\frac{\widetilde{\mathcal{A}}_1}{\lambda}+\frac{\widetilde{\mathcal{A}}_2}{\lambda-x}+\frac{\widetilde{\mathcal{A}}_3}{\lambda-1}\right)_{\hbox{\rm entry }12}=\frac{k\left(\lambda-\dfrac{x\widetilde{v}_1z_1}{k}\right)}{\lambda(\lambda-x)(\lambda-1)}=\frac{k\left(\lambda-y\right)}{\lambda(\lambda-x)(\lambda-1)}.
   $$}
   Thus, $y$ is the solution for the unknown  $\lambda$ of the equation 
{\small
$
 \left(\frac{\widetilde{\mathcal{A}}_1}{\lambda}+\frac{\widetilde{\mathcal{A}}_2}{\lambda-x}+\frac{\widetilde{\mathcal{A}}_3}{\lambda-1}\right)_{12}=0$.
}
  Namely, 
\be
\label{18giugno2021-2}
   y(x)=\frac{x~(\widetilde{\mathcal{A}}_1)_{12}}{x[(\widetilde{\mathcal{A}}_1)_{12}+(\widetilde{\mathcal{A}}_3)_{12})]-(\widetilde{\mathcal{A}}_3)_{12}}.
\ee
The above shows again that $y(x)$ does not depend on $k^0=\varepsilon_3/\varepsilon_1$. 
}
\ere
%%%%%%%%%%%%%%
\subsubsection{Step 5. Completion of the proof of Theorem \ref{20gennaio2021-2} and explicit formulae}
   From the above discussion,  we are able to explicitly write $ z_1,z_2,z_3$ and $\widetilde{v}_1,\widetilde{v}_2,\widetilde{v}_3$ in terms of PVI transcendents $y(x)$ and $dy(x)/dx$. 
 The entries of $\boldsymbol{Z}(x)$ can then be expressed through \eqref{22gennaio2021-6}   in terms of $y(x)$.
  Firstly, substituting  (C.51), we see that 
 \begin{align*}
 &
 \boldsymbol{Z}_{12}=z_2-(z_1+\theta_1)\frac{x-y}{(1-x)y}, & \boldsymbol{Z}_{13}=z_3-(z_1+\theta_1)\frac{x(y-1)}{(1-x)y},
 \\
 &
 \boldsymbol{Z}_{21}=z_1-(z_2+\theta_2)\frac{(1-x)y}{x-y}, & \boldsymbol{Z}_{23}=z_3-(z_2+\theta_2)\frac{x(y-1)}{x-y},
 \\
 &
 \boldsymbol{Z}_{31}=z_1-(z_3+\theta_3)\frac{(1-x)y}{x(y-1)}, & \boldsymbol{Z}_{32}=z_2-(z_3+\theta_3)\frac{x-y}{x(y-1)},
 \end{align*} 
 Substituting  (C.52) in the above, where $z$ is given by \eqref{23gennaio2021-3}, we obtain the explicit expressions of $ \boldsymbol{Z}(x)$ in terms of $y(x)$ and $dy(x)/dx$,  that  yield  $\Omega_{ij}(x)=\boldsymbol{Z}_{ij}(x)k_i(x)/k_j(x) $, $1\leq i\neq j \leq 3$, $k_3:=1$, as  in the statement of 
  Theorem \ref{20gennaio2021-2}.
 
  To complete the proof of Theorem \ref{20gennaio2021-2}, we find the differential equations for $k_1(x)$ and $k_2(x)$.  
 The factorization \eqref{22gennaio2021-10}  implies that 
 $$
 \widehat{\Omega}_2=K \boldsymbol{Z}_2 K^{-1},\quad\quad
 \boldsymbol{Z}_2
 =
 \begin{pmatrix}
 0 & -\boldsymbol{Z}_{12}/x & 0 
 \\
- \boldsymbol{Z}_{21}/x & 0 & \boldsymbol{Z}_{23}/(1-x)
\\
0 & \boldsymbol{Z}_{32}/(1-x) & 0 
 \end{pmatrix}
 $$ 
 Substituting the above and \eqref{22gennaio2021-10} into \eqref{20gennaio2021-26}  we find 
 $$ 
 \Bigl[ K^{-1} \frac{dK}{dx}~,~ \boldsymbol{Z} \Bigr]=[\boldsymbol{Z},\boldsymbol{Z}_2]-\frac{d\boldsymbol{Z}}{dx},
 $$
 namley
 $$
  \frac{d \ln k_1}{dx}= \left\{
  \begin{array}{c}
  \dfrac{1}{\boldsymbol{Z}_{13}}\left([\boldsymbol{Z},\boldsymbol{Z}_2]_{13}-\dfrac{d\boldsymbol{Z}_{13}}{dx} \right)
  \\
  \noalign{\medskip}
  \dfrac{1}{\boldsymbol{Z}_{31}}\left(\dfrac{d\boldsymbol{Z}_{31}}{dx} -[\boldsymbol{Z},\boldsymbol{Z}_2]_{31}\right)
  \end{array}
  \right. ,
  \quad 
   \frac{d \ln k_2}{dx}= \left\{
  \begin{array}{c}
  \dfrac{1}{\boldsymbol{Z}_{23}}\left([\boldsymbol{Z},\boldsymbol{Z}_2]_{23}-\dfrac{d\boldsymbol{Z}_{23}}{dx} \right)
  \\
  \noalign{\medskip}
  \dfrac{1}{\boldsymbol{Z}_{32}}\left(\dfrac{d\boldsymbol{Z}_{32}}{dx} -[\boldsymbol{Z},\boldsymbol{Z}_2]_{32}\right)
  \end{array}
  \right. 
  $$
  Substituting the expressions of the entries of $ \boldsymbol{Z}$ and $ \boldsymbol{Z}_2$ in terms of $y$ and $dy/dx$, a   computations  show that the r.h.s. of the above expressions respectively  are $l_1(x)$ and $l_2(x)$ in the statement of Theorem \ref{20gennaio2021-2}.   Notice that $k_1(x)$ and $k_2(x)$ contain the multiplicative integration constants $k_1^0$ and $k_2^0$ responsible for the correspondence 
  $$
  y\quad \longleftrightarrow \quad \hbox{ class } \{K^0\cdot  V \cdot (K^0)^{-1},\quad K^0=\hbox{\rm diag}  (k_1^0,k_2^0,1),\quad k_1^0,k_2^0\in\mathbb{C}\backslash\{0\} \}
$$
The proof of Theorem \ref{20gennaio2021-2} is complete. $\Box$

\bre
{\rm
From \eqref{18giugno2021-2} and \eqref{22gennaio2021-8}:
$$ 
y= \frac{x(\widetilde{v}_1z_1)/(\widetilde{v}_3z_3)}{x(1+(\widetilde{v}_1z_1)/(\widetilde{v}_3z_3))-1}.
$$
Computing $(\widetilde{v}_1z_1)/(\widetilde{v}_3z_3)$ by solving the six equations \eqref{22gennaio2021-6} for the ratios $(\widetilde{v}_j(x)z_j(x))/(\widetilde{v}_i(x)z_i(x))$,  we obtain 
$$ 
y(x)= \frac{xR(x)}{x(1+R(x))-1},\quad R(x):=-\frac{\theta_1+Z_{31}(x)}{\theta_3+Z_{13}(x)} =-\frac{\theta_1+k_1(x)\Omega_{31}(x)}{\theta_3+\Omega_{13}(x)/k_1(x)}
$$
}
\ere

%%%%%%%%%%%%%%%%%

\subsection{A final remark}
\label{25giugno2021-1}
Another strategy is followed in \cite{MazzoIrr}  to obtain the parameterization of $V$  in terms of PVI transcendents: the Schlessinger equations \eqref{21gennaio2021-1}  are reduced to a system for three functions $q,p,k$:
\begin{equation}
\label{reducedschlesinger}
{\small
\begin{cases}
\begin{aligned}
&\frac{\partial q}{\partial u_i}=\frac{P(q)}{P'(u_i)}\left[2p+\frac{1}{q-u_i}-\sum_{k=1}^3\frac{\theta_k}{q-u_k}\right],
\\[2ex]
&\frac{\partial p}{\partial u_i}=-\Bigg\{P'(q)p^2+\Bigg[2q+u_i-\sum_{k=1}^3\theta_k\Bigg(2q+u_k-\sum_{j=1}^3u_j\Bigg)\Bigg]p,
\\[2ex]
&\qquad\quad+\frac{1}{4}\Bigg(\sum_{k=1}^3\theta_k-\theta_\infty\Bigg)\Bigg(\sum_{k=1}^3\theta_k+\theta_\infty-2\Bigg)\Bigg\}\frac{1}{P'(u_i)},
\\[2ex]
&\frac{\partial\log(k)}{\partial u_i}=(\theta_\infty-1)\frac{q-u_i}{P'(u_i)},
\end{aligned}
\end{cases}
i=1,2,3,}
\end{equation}
which are in turn reduced to  PVI  for a function $y(x)$, where $x=(u_2-u_1)/(u_3-u_1)$ and $y=(q-u_1)/(u_3-u_1)$. The functions $q,p,k$ are defined by
\[
(\mathcal{A}(q,u))_{12}=0,\quad p=\sum_{k=1}^{3}\frac{(\mathcal{A}_k)_{11}+\theta_k/2}{q-u_k},\quad k=\frac{2P(\lambda)(\mathcal{A}(\lambda,u))_{12}}{\theta_\infty(q-\lambda)},
\]
where
\[
\mathcal{A}(z,u)=\frac{\mathcal{A}_1(u)}{z-u_1}+\frac{\mathcal{A}_2(u)}{z-u_2}+\frac{\mathcal{A}_3(u)}{z-u_3}\quad\mbox{and}\quad P(\lambda)=(\lambda-u_1)(\lambda-u_2)(\lambda-u_3),
\]
so that the matrices $\mathcal{A}_i$ are parameterized by
\begin{equation}
\label{a11}
\begin{aligned}
&(\mathcal{A}_i)_{11}=\frac{1}{\theta_\infty P'(u_i)}\left\{P(q)(q-u_i)p^2+P(q)(q-u_i)p\left(\frac{\theta_\infty}{q-u_i}-\sum_{k=1}^3\frac{\theta_k}{q-u_k}\right)\right.\\[2ex]
                                     &\qquad\qquad+(q-u_i)\left[\frac{\theta_\infty^2}{4}\left(q+2u_i-\sum_{k=1}^3u_k\right)+\sum_{k=1}^3\frac{\theta_k^2}{4}\left(q+2u_k-\sum_{j=1}^3u_j\right)\right]\\[2ex]
                                     &\qquad\qquad\left.+\frac{q-u_i}{4}\sum_{k\ne j\ne l}\left(\theta_j\theta_k(q-u_l)\right)-\frac{\theta_\infty}{2}P(q)\sum_{k=1}^3\frac{\theta_k}{q-u_k}\right\},\\[2ex]
&(\mathcal{A}_i)_{12}=\theta_\infty\frac{q-u_i}{2P'(u_i)}k,\quad (\mathcal{A}_i)_{21}=\frac{1}{(\mathcal{A}_i)_{12}}\left(\frac{\theta_i^2}{4}-(\mathcal{A}_i)_{11}^2\right),\quad (\mathcal{A}_i)_{22}=-(\mathcal{A}_i)_{11}.
\end{aligned}
\end{equation}
By comparison of~\eqref{a11} with~\eqref{20gennaio2021-25} we can find $a_k,b_k$ as functions of $u,q,p,k$, and hence, inserting into~\eqref{20gennaio2021-23} we get the following parameterization of the off-diagonal elements of the matrix $V$ in terms of $u,q,p$:
\[
V_{ij}=\frac{1}{u_j-u_l}\left[(q-u_j)(q-u_l)p+\frac{\theta_\infty-\theta_1-\theta_2-\theta_3}{2}q+\frac{\theta_i-\theta_j+\theta_l-\theta_\infty}{2}u_j+\theta_ju_l\right]\frac{d_i}{d_j},
\]
$i,j,l=1,2,3$, $i\ne j\ne l$. From the isomonodromicity condition \eqref{20gennaio2021-6}-\eqref{20gennaio2021-7}, the functions $d_i/d_j$ are determined by quadratures from
\begin{align*}
&\frac{\partial}{\partial u_1}\log\left(\frac{d_1}{d_3}\right)=\frac{\theta_1-\theta_3}{u_1-u_3}+\frac{p(q-u_2)-\theta_2}{u_1-u_2},
&&
\frac{\partial}{\partial u_1}\log\left(\frac{d_2}{d_3}\right)=\frac{u_2-u_3}{(u_1-u_2)(u_1-u_3)}[p(q-u_1)-\theta_1],
\\[2ex]
&\frac{\partial}{\partial u_2}\log\left(\frac{d_1}{d_3}\right)=\frac{u_1-u_3}{(u_2-u_1)(u_2-u_3)}[p(q-u_2)-\theta_2],
&&
\frac{\partial}{\partial u_2}\log\left(\frac{d_2}{d_3}\right)=\frac{\theta_2-\theta_3}{u_2-u_3}+\frac{p(q-u_1)-\theta_1}{u_2-u_1},
\\[2ex]
&\frac{\partial}{\partial u_3}\log\left(\frac{d_1}{d_3}\right)=-\frac{\theta_1-\theta_3}{u_1-u_3}+\frac{p(q-u_2)-\theta_2}{u_2-u_3},
&&
\frac{\partial}{\partial u_3}\log\left(\frac{d_2}{d_3}\right)=-\frac{\theta_2-\theta_3}{u_2-u_3}+\frac{p(q-u_1)-\theta_1}{u_1-u_3}.
\end{align*}
The explicit parameterization in terms of PVI transcendents follows from Remark \ref{par-reduction}.

\bre
{\rm 
\label{par-reduction}
The gauge transformation $Y(z,u)=e^{zu_1}\widehat{Y}(z,u)$ and change of independent variable $z=\widehat{z}/(u_3-u_1)$ 
change the $z$-component \eqref{systempainleve-bis} to 
$$ 
\frac{d\widehat{Y}}{dz}=\left(\widehat{U}+\frac{V}{\widehat{z}}\right) Y,\quad\quad \widehat{U}=\mbox{\rm diag}\left(0,\frac{u_2-u_1}{u_3-u_1},1\right)$$ 
Thus we can assume $u_1=0$, $u_2=x$, $u_3=1$ without loss of generality. With this setting, the $z$-component of  \eqref{19febbraio2021-11} becomes system 
\eqref{systempainleve}.
}\ere  
  
  \section{Monodromy Data}
  \label{25gennaio2021-1}
  
  As mentioned  in the introduction, the monodromy data of the $2\times 2$ Fuchsian system 
  \eqref{20gennaio2021-21} are used   to parameterize Painlev\'e VI transcendents and solve their non-linear connection problem. The ring of  invariant functions of $SL(2,\mathbb{C})^3/(\mathcal{M}_j\mapsto C^{-1}\mathcal{M}_j C, ~\det C\neq 0)$ is generated by  the traces 
         $$
     p_{jk}=\hbox{\rm tr}(\mathcal{M}_j\mathcal{M}_k), \quad j\neq k\in\{1,2,3\},
     $$ 
     $$
    p_\infty=\hbox{\rm Tr}( \mathcal{M}_{j_3} \mathcal{M}_{j_2} \mathcal{M}_{j_1})= 2\cos(\pi\theta_\infty),\quad j_1\prec j_2\prec j_3;\quad\quad p_j=\hbox{Tr} \mathcal{M}_j=   2\cos(\pi\theta_j).
    $$
     The ordering relation above is explained in \eqref{16aprile2021-3}. 
The goal of this section is  formula  \eqref{3febbraio2021-2}, expressing the $p_{jk}$ {\it in terms of the Stokes matrices} of the $3\times 3$ system \eqref{systempainleve-bis}.

\bde The {\bf Stokes rays} associated with $U(u)=\hbox{\rm diag}(u_1,u_2,u_3)$ are the infinitely many half-lines in the universal covering of the punctured $z$-plane  $\mathbb{C}\backslash\{0\}$, issuing from $z=0$ towards $\infty$, defined by 
 $\Re((u_j-u_k)z)=0$, $\Im((u_j-u_k)z)<0$, for $u_j\neq u_k$. 
 \ede

In the following, we will work under  the analyticity assumptions of  \underline{Case 1} or \underline{Case 2} on the domain $\mathbb{D}$ explained in the beginning of Section \ref{19giugno2021-6}, with  the following  refinement of the size of $\mathbb{D}$.

 \underline{Case 1}. $\mathbb{D}=\mathbb{D}(u^0)$.   Let $\eta^{(0)}\in\mathbb{R}$ satisfy 
     \be
 \label{27gennaio2021-2}
 \eta^{(0)}\neq \arg(u_i^0-u_j^0)  \hbox{ mod }\pi, ~ \forall~i\neq j,
 \ee
 so that $$\tau^{(0)}:=3\pi/2-\eta^{(0)}$$ is an {\it admissible direction} in the $z$-plane for the Stokes rays of $U(u^0)$, that is  no such rays have directions  $\tau^{(0)}+h\pi$, $h\in\mathbb{Z}$. 
The size of  $\mathbb{D}(u^0)$ is so  small that  the Stokes rays of $U(u)$ in the $z$-plane do not  cross the directions  $\tau^{(0)}+h\pi$, as $u$ varies in $\mathbb{D}(u^0)$.

  \underline{Case 2}. $\mathbb{D}=\mathbb{D}(u^c)$, so that $u^c=(u_1^c,u_2^c,u_3^c)$ has only two distinct components\footnote{ $u^c=(\lambda_1,\lambda_1,\lambda_2)$, or $(\lambda_1,\lambda_2,\lambda_2)$, or $(\lambda_1,\lambda_2,\lambda_1)$.}  $\lambda_1,\lambda_2$.  
    Let $ \eta\in\mathbb{R}$ satisfy
  $
 \eta\neq \arg(u_i^c-u_j^c) \hbox{ mod }\pi$, $\forall~i\neq j$ such that $u_i^c\neq u_j^c$, 
that is  
  \be
\label{27gennaio2021-3}
\eta\neq  \arg(\lambda_1-\lambda_2) \hbox{ mod $\pi$}.
\ee 
Then, $$\tau:=3\pi/2-\eta$$ is an {\it admissible direction} in the $z$-plane for the Stokes rays of $U(u^c)$. 
The size of $\mathbb{D}(u^c)$ must  be sufficiently small so that no Stokes rays  associated with pairs $(u_j,u_k)$ such that $u_j^c\neq u_k^c$ cross the admissible directions $\arg z= \tau +h\pi$, $h\in\mathbb{Z}$, as $u$ varies in $\mathbb{D}(u^c)$. 

\bre
  \label{16aprile2021-2}
  {\rm 
   The bound for  the size  of 
  $
  \mathbb{D}(u^c)=\{u\in\mathbb{C}^3~|~ \max_{1\leq l\leq 3} |u_j-u_j^c|\leq \epsilon_0\} 
  $
  is described in \cite{CDG,guz2021} as follows. Consider  in the $\lambda$-plane
      half-lines  $\mathcal{L}_1(\eta)$,  $\mathcal{L}_2(\eta)$ respective from $\lambda_1$ and $\lambda_2$ 
    to infinity in direction $\eta$. Let $\delta=  \min_{\rho>0}|\lambda_1-\lambda_2+\rho \exp\{\eta \sqrt{-1}\}|$ be their distance. Then $\epsilon_0<\delta/2$. 
 }   \ere

    The Pfaffian systems  \eqref{19febbraio2021-11} and \eqref{20gennaio2021-4} are related by the gauge transformation 
  $ 
  Y_{\eqref{19febbraio2021-11}}= z~ Y_{\eqref{20gennaio2021-4}}$. 
Therefore, the connection and Stokes matrices  are the same for  the two systems   \be
\label{16aprile2021-1}
   \frac{dY}{dz}=\left(U+\frac{W}{z}\right) Y, \quad \quad  \quad W=
   \left\{
   \begin{array}{c}
  \hbox{ either } V, 
   \\
 \hbox{ or }  \widetilde{V}:=V-I.
   \end{array}
   \right.
   \ee
With the above assumptions on $\mathbb{D}$,   according to \cite{CDG} system \eqref{16aprile2021-1}  admits a unique formal solution
  \be
  \label{27gennaio2021-5}
Y_F(z,u)=  (I+\sum_{k\geq 1} F_k(u)z^{-k}) z^{{\rm diag} (W)}e^{zU},
\ee
with matrix coefficients $F_k(u)$ holomorphic on $\mathbb{D}$.  To it, there correspond unique fundamental matrix solutions $Y_1(z,u)$, $Y_2(z,u)$, $Y_3(z,u)$, that are sometimes called {\it canonical solutions},  such that
\be
 \label{27gennaio2021-6}
Y_j(z,u)\sim Y_F(z,u), \quad \quad  \hbox{ $z\to \infty$ in $\mathcal{S}_j$},\quad j=1,2,3,
\ee
where, depending on Case 1 or 2, 
\begin{align*}
&\mathcal{S}_1:\quad (\tau^{(0)}-\pi)-\varepsilon<\arg z <  \tau^{(0)}+\varepsilon, & \hbox{ or } \quad&
 (\tau-\pi)-\varepsilon<\arg z <  \tau+\varepsilon,
   \\
& \mathcal{S}_2:\quad \tau^{(0)}-\varepsilon<\arg z <  (\tau^{(0)}+\pi)+\varepsilon, & \hbox{ or } \quad &
  \tau-\varepsilon<\arg z <  (\tau+\pi)+\varepsilon,
  \\
& \mathcal{S}_3:\quad (\tau^{(0)}+\pi)-\varepsilon<\arg z <  (\tau^{(0)}+2\pi)+\varepsilon ,
&  \hbox{ or } 
  \quad &
  (\tau+\pi)-\varepsilon<\arg z <  (\tau+2\pi)+\varepsilon.
 \end{align*}

\bre
\label{7luglio2021-3}
{\rm
In the case  $\mathbb{D}=\mathbb{D}(u^c)$, the notion of  {\bf partial resonance} for $V(u^c)$ is introduced  in corollary 4.1 of  \cite{CDG}. The name  is due to \cite{sabbah}.  In our case, with vanishing conditions \eqref{19giugno2021-3}, partial resonance simply means
   $$ 
   \theta_i-\theta_j\in\mathbb{Z}\backslash\{0\} \quad \hbox{ for $i\neq j$ such that $u_i^c=u_j^c$}.
   $$
 In case of partial resonance, system \eqref{16aprile2021-1} at fixed $u=u^c$ has a family of formal 
 solutions depending on a finite number of parameters, and only one of them is $Y_F(z,u^c)$.  
 With no partial resonance, the formal solution  is unique and coincides
  with $Y_F(z,u^c)$. See\footnote{ The  statement of corollary 1.1. of \cite{CDG} is imprecise: it is not that ``the diagonal entries 
of $\widehat{A}_1(0)$ do not differ by non-zero integers'', but the elements   of each sequence 
 $(\widehat{A}_1(0))_{j_1j_1},(\widehat{A}_1(0))_{j_2j_2},...,(\widehat{A}_1(0))_{j_\ell j_\ell}$ 
 corresponding to $u_{j_1}(0)=u_{j_2}(0)=...=u_{j_\ell}(0)$, which   precisely is the partial resonance in the context of \cite{CDG}.}    corollaries 1.1 and 4.1 of \cite{CDG}. 
 In the specific cases considered in this paper for $x=(u_2-u_1)/(u_3-u_1)\to 0$, partial resonance means $\theta_1-\theta_2\in\mathbb{Z}\backslash\{0\}$. }
   \ere
   
  \bde The {\bf Stokes matrices} are the connection matrices such that   
\be
\label{27gennaio2021-9}
  Y_2(z,u)=Y_1(z,u)\mathbb{S}_1,
  \quad\quad Y_3(z,u)=Y_2(z,u)\mathbb{S}_{2}.
  \ee 
  \ede
  They are constant on $\mathbb{D}$ by the integrability of  \eqref{19febbraio2021-11} or \eqref{20gennaio2021-4}.
    This result is standard in case of $\mathbb{D}(u^0)$, essentially following \cite{JMU}, while it follows from  theorem 1.1 of  \cite{CDG} in case of $\mathbb{D}(u^c)$. In the latter case,   
$$ 
(\mathbb{S}_{1})_{ij}=(\mathbb{S}_{1})_{ji}=(\mathbb{S}_{2})_{ij}=(\mathbb{S}_{2})_{ji}=0
\quad \hbox{ if $u_i^c=u_j^c$}.
$$
For example, in the case $u^c=(\lambda_1,\lambda_1,\lambda_2)$, we have $(\mathbb{S})_{12}=(\mathbb{S})_{21}=0$.

\bre
{\rm
    The matrices  $\mathbb{S}_1$ and $\mathbb{S}_2$ correspond to $\mathbb{S}_\nu$ and $\mathbb{S}_{\nu+\mu}$ in  \cite{CDG,guz2021}. The sectors in \eqref{27gennaio2021-6} correspond to the sectors $\mathcal{S}_\nu(\mathbb{D}(u^*))$,  $\mathcal{S}_{\nu+\mu}(\mathbb{D}(u^*))$ and  $\mathcal{S}_{\nu+2\mu}(\mathbb{D}(u^*))$ of \cite{guz2021}, where $u^*:=u^0 $ or $ u^c$,
   depending on the type of polydisc considered. 
    }
    \ere

\begin{figure}
\centerline{\includegraphics[width=0.5\textwidth]{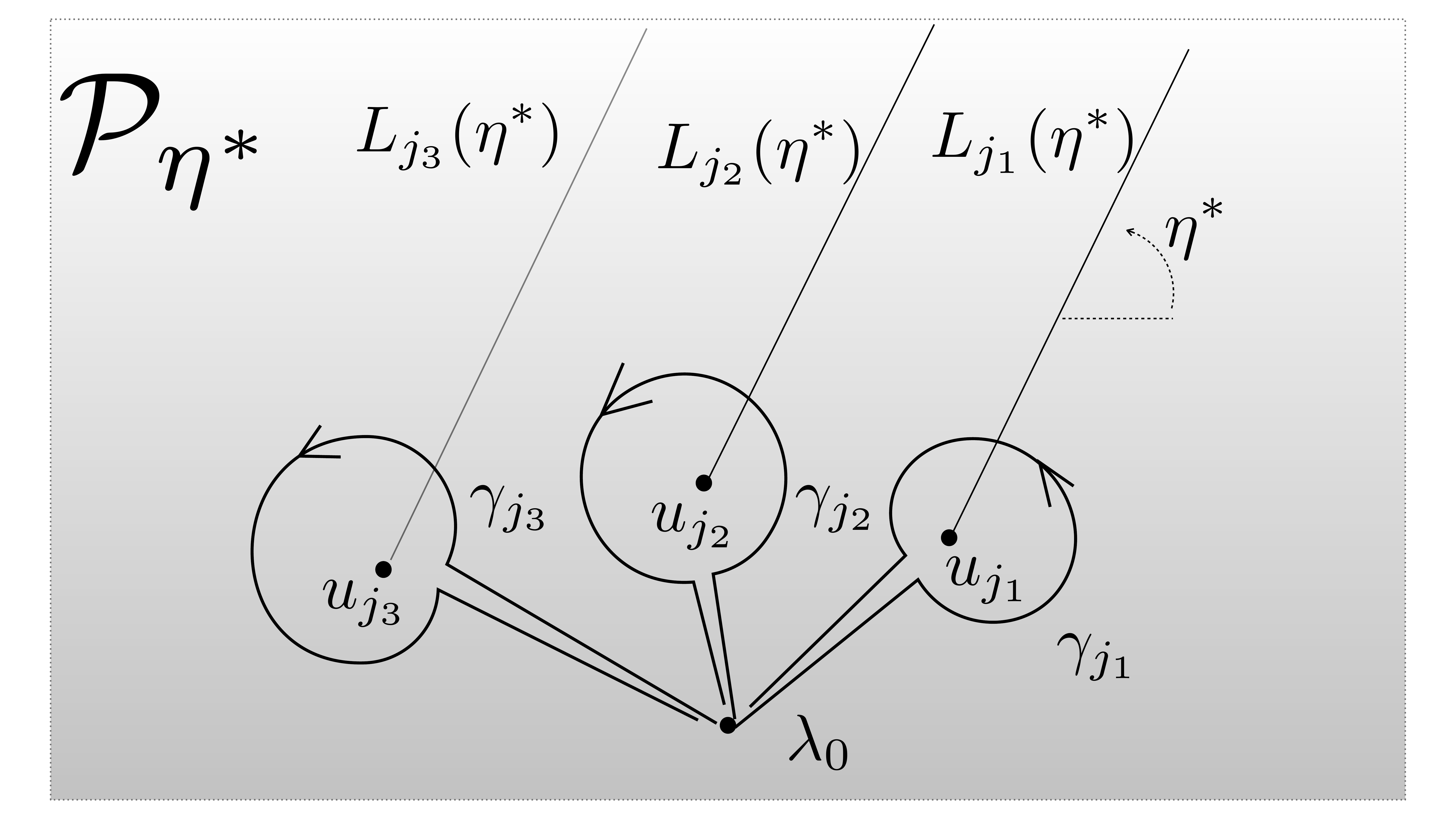}}
\caption{$\mathcal{P}_{\eta^*}$, branch-cuts and basic loops.}
\label{3agosto2021-6}
\end{figure}

   To proceed, we need an important fact:  the Stokes matrices can be expressed in terms of the  connection coefficients (defined below) of certain {\it selected column vector solutions} $\vec{\Psi}_1(\lambda,u),\vec{\Psi}_2(\lambda,u),\vec{\Psi}_3(\lambda,u)$ of system \eqref{20gennaio2021-3}. In order to introduce the selected solutions, let an admissible direction as in \eqref{27gennaio2021-2}  or \eqref{27gennaio2021-3} be denoted for short by
  $$\eta^*: =\eta^{(0)} \hbox{  or } \eta.
  $$ 
 In the $\lambda$-plane, let $L_1(\eta^*),L_2(\eta^*),L_3(\eta^*)$ be  oriented  half-lines from $u_1,u_2,u_3$ respective to infinity in direction $\eta^*$.  Following \cite{BJL4}, let  $\mathcal{P}_{\eta^*}(u)$ be the $\lambda$-plane with these {\it branch cuts} and with the determinations
  $$ 
  \eta^*-2\pi <\arg(\lambda-u_k)<\eta^*.
  $$
  of  $\ln(\lambda-u_k)=\ln|\lambda-u_k|+i\arg(\lambda-u_k)$.  See Figure \ref{3agosto2021-6}. The domain of definition  of  the  solutions of \eqref{20gennaio2021-3} is the set 
  $$ 
  \mathcal{P}_{\eta^*}(u)\hat{\times} \mathbb{D}(u^*):= \{(\lambda,u)~|~u\in\mathbb{D},~\lambda\in \mathcal{P}_{\eta^*}(u)\}
  \equiv 
  \bigcup_{u\in\mathbb{D}}\Bigl(  \mathcal{P}_{\eta^*}(u)\times\{u\}\Bigr).
  $$
  The symbol $\hat{\times}$ for the $u$-dependent Cartesian product is borrowed from \cite{JMU}. 
  Let $\mathbb{Z}_{-}:=\{-1,-2,-3,\dots\}$ be the negative integers. 
 The three selected  vector solutions are  uniquely determined  in theorem 5.1 of \cite{guz2021}, to which we refer (with $n=3$ and  the identification $\theta_k=-\lambda_k^\prime-1$,  being $\lambda_k^\prime$ used in \cite{guz2021}).   They are holomorphic on $\mathcal{P}_\eta(u)\hat{\times} \mathbb{D}(u^*)$, where $u^*=u^c$ or $u^0$, depending on the polydisc considered. 
A  solution $\vec{\Psi}_k$ has a branching point at $\lambda=u_k$  in case $\theta_k\not\in \mathbb{Z}_{-}$, and for $\lambda\in\mathcal{P}_{\eta^*}(u)$ its monodromy  corresponding to a loop $\gamma_j: (\lambda-u_j)\longmapsto (\lambda-u_j)e^{2\pi i }$ is given in  \cite{guz2021} by:
\be
\label{24gennaio2021-1}
\vec{\Psi}_k\longmapsto 
\begin{cases}
\begin{aligned}
&e^{2\pi i \theta_k} \vec{\Psi}_k , & j=k,
\\[1ex]
& \vec{\Psi}_k +\alpha_j  c_{jk} \vec{\Psi}_j , & j\neq k,
\end{aligned}
\end{cases}
\quad \quad \alpha_j:= 
\left\{
\begin{array}{cc}
e^{2\pi i \theta_j}-1,& \theta_j\not\in\mathbb{Z}
\\
\noalign{\medskip}
2\pi i ,& \theta_j\in\mathbb{Z}
\end{array}
\right.,
\quad 
\quad j,k\in\{1,2,3\}
\ee
with certain  ``connection'' coefficients $c_{jk}$. Since $\alpha_k=e^{2\pi i \theta_k}-1$ for $\theta_k\not\in\mathbb{Z}$ and $\alpha_k=2\pi i$ for $\theta_k\in\mathbb{Z}$, the above formulae imply that  
$$ 
c_{kk}=1 \hbox{ for $\theta_k\not\in\mathbb{Z}$;}\quad \quad  c_{kk}=0\hbox{ for $\theta_k\in\mathbb{Z}$.}
$$
If  $\theta_k\in \mathbb{Z}_{-}$, in some particular cases which depend on the specific $V$ it may happen that $\vec{\Psi}_k\equiv 0$. 
It is also deducible from  \cite{guz2016,guz2021} that 
\begin{itemize}
\item If $\theta_k\in\mathbb{Z}_{-}=\{-1,-2,\dots\}$ and $\vec{\Psi}_k\equiv 0$, then $c_{jk}=0$ for every $j$. 
\item If $\theta_k\in\mathbb{Z}_{+}=\{1,2,\dots\}$ and there is no singular solution at $\lambda= u_k$, then  $c_{kj}=0$ for every $j$. 
\end{itemize}
It is proved in \cite{guz2021} that the coefficients $c_{jk}$ are  independent of  $u\in \mathbb{D}$, so they are called  {\bf isomonodromic connection coefficients}. 
 As a consequence, the  transformation \eqref{24gennaio2021-1} holds for every $u$ in the polydisc.  In particular, in case of coalescences, we have
$$
c_{jk}=0\quad\hbox{ for $j\neq k$ such that $u_j^c=u_k^c$}.
$$

\bre
{\rm
From  $c_{jk}=0$,  the monodromy transformations \eqref{24gennaio2021-1} make sense also in case of coalescences in $\mathbb{D}=\mathbb{D}(u^c)$.  The base point $\lambda_0$ for the loops $\gamma_j$ is taken in $\bigcap_{u\in\mathbb{D}}\mathcal{P}_{\eta^*}(u)$ as in figure \ref{3agosto2021-6}. In case of coalescence $u_j^c=u_k^c$,    the two branch cuts $L_j(\eta)$ and $L_k(\eta)$ can overlap as $u$ varies in $\mathbb{D}(u^c)$, but this causes no difficulties because $c_{jk}=0$, so that $\vec{\Psi}_k$ has trivial monodromy at $\lambda=u_j$.}
\ere
It is important to notice that the matrix solution 
\be
\label{4febbraio2021-1} 
\Psi(\lambda,u)=\Bigl( \vec{\Psi}_1(\lambda,u)~\Bigl|~  \vec{\Psi}_2(\lambda,u)~\Bigr|~ \vec{\Psi}_3(\lambda,u) \Bigr)
\ee
of \eqref{20gennaio2021-3} has constant monodromy, but it is not necessarily a fundamental matrix.\footnote{ In case $u^*=u^c$, it is proved in \cite{guz2021} that if $u_j^c=u_k^c$ for some $j\neq k$, then $\vec{\Psi}_j$ and $\vec{\Psi}_k$ are either linearly independent or at least one of them is identically zero (possibly if $\theta_j$ or $\theta_k$ is in $\mathbb{Z}_{-}$).} 
It is fundamental if  for example 
$V$ has no integer eigenvalues. If it happens that $V$ has some  
integer eigenvalues and $\Psi(\lambda,u)$ is fundamental, then necessarily there is at least one  $\theta_k\in\mathbb{Z}$. This follows  from  \cite{BJL4} for the generic case and \cite{guz2016,guz2021} for the general case.
 \vskip 0.3 cm 
 
If $\mathbb{D}$ is as small has we have previously specified,   an ordering relation is well defined in $\{1,2,3\}$:  
 \begin{itemize}
 \item In case  $\mathbb{D}=\mathbb{D}(u^0)$,  the ordering 
\be
\label{16aprile2021-3}
 j\prec k\quad\Longleftrightarrow\quad \Re(e^{i\tau^{(0)}}(u_j-u_k))<0,\quad  j\neq k .
 \ee
  is well defined for all  $u\in \mathbb{D}(u^0)$. It coincides with the ordering defined by  $\Re(e^{i\tau^{(0)}}(u_j^0-u_k^0))<0$. Equivalently,   $j\prec k$ if and only if  $L_j(\eta^{(0)})$ is to the right of  $L_k(\eta^{(0)})$ in the $\lambda$-plane, where right or left refer to the orientation in direction $\eta^{(0)}$ towards infinity. 
 
  \item In case of $\mathbb{D}=\mathbb{D}(u^c)$,  there is no ordering relation for $j,k$ such that $u_j^c=u_k^c$, while
 \be
 \label{16aprile2021-4}
 j\prec k\quad\Longleftrightarrow\quad \Re(e^{i\tau}(u_j-u_k))<0,\quad j\neq k \hbox{ and } u_j^c\neq u_k^c.
 \ee
  If $\epsilon_0$ is as in Remark \ref{16aprile2021-2}, then the partial ordering is well defined for all $u\in\mathbb{D}(u^c)$.  Notice that  $j\prec k$ if and only if  $L_j(\eta)$ is to the right of  $L_k(\eta)$ in the $\lambda$-plane, or equivalently if  $\mathcal{L}_j(\eta)$ is to the right of  $\mathcal{L}_k(\eta)$ (see Remark \ref{16aprile2021-2}).\footnote{ In our case, if for example $u_1^c=u_2^c=\lambda_1  $ and $u_3^c=\lambda_2$, we have $1\prec 3$ iff $\mathcal{L}_1(\eta)$ is to the right of $\mathcal{L}_2(\eta)$ in the $\lambda$-plane, where  $\mathcal{L}_1(\eta)$, $\mathcal{L}_2(\eta)$ are defined in Remark \ref{16aprile2021-2}. There is no relation like $1\prec 2$ or $2\prec 1$.}
 \end{itemize}

Following  \cite{guz2016,guz2021}, the connection coefficients are related to the entries of the Stokes matrices of system \eqref{16aprile2021-1} by
  \be
  \label{27gennaio2021-1} 
   c_{jk}=
   \left\{
   \begin{array}{cc}
    \dfrac{e^{2\pi i \theta_k}}{\alpha_k}(\mathbb{S}_1)_{jk}, & j\prec k,
    \\
    \\
     -\dfrac{ (\mathbb{S}_{2}^{-1})_{jk}}{e^{2\pi i(\theta_j-\theta_k)}\alpha_k}, & j\succ k.
   \end{array}
   \right.
   \ee
We are ready to state the main result of the section. 
\bth
\label{16aprile2021-5}
Let $\mathbb{D}$ be  $\mathbb{D}(u^0)$, or $\mathbb{D}(u^c)$.
For every point $u^\bullet$ in  $\mathbb{D}(u^0)$, or in $\mathbb{D}(u^c)\backslash  \Delta_{\mathbb{C}^3}$, there is a neighborhood $\mathcal{U}$ of $u^\bullet$ in $\mathbb{D}$  and a  fundamental matrix solution $\Phi_{\rm hol}(\lambda,u)$ of system  \eqref{20gennaio2021-21}, holomorphic of $(\lambda,u)\in \mathcal{P}_{\eta^*}(u) \hat{\times}~ \mathcal{U}$ (here $\eta^*=\eta^{(0)}$ or $\eta$), whose monodromy invariants 
$$
p_{jk}:=\hbox{\rm tr} (\mathcal{M}_j(u)\mathcal{M}_k(u)),
$$
 are independent of $u$. Here,  $\mathcal{M}_j(u)$ is the monodromy matrix at $\lambda=u_j$ of $\Phi_{\rm hol}(\lambda,u)$.  
 They are   expressed in terms of the Stokes matrices of system \eqref{16aprile2021-1}:
  \be
   \label{3febbraio2021-2}
   p_{jk}
   =
   \left\{
  \begin{array}{cc}
 2\cos\pi(\theta_j-\theta_k)-e^{i\pi(\theta_j-\theta_k)} (\mathbb{S}_1)_{jk}(\mathbb{S}_{2}^{-1})_{kj}, & j\prec k,
  \\
  \\
 2\cos\pi(\theta_j-\theta_k)-e^{i\pi(\theta_k-\theta_j)}(\mathbb{S}_1)_{kj}(\mathbb{S}_{2}^{-1})_{jk},
  , & j\succ k.
  \end{array}
  \right.
  \ee
 with ordering  \eqref{16aprile2021-3} or \eqref{16aprile2021-4} according to the polydisc being either $\mathbb{D}(u^0)$ or  $\mathbb{D}(u^c)$. In the latter case, 
 $$ 
   p_{jk}= 2\cos\pi(\theta_j-\theta_k)\quad \hbox{ for $j\neq k$ such that $u_j^c=u_k^c$}.
   $$
\eth

To appreciate the general validity of \eqref{3febbraio2021-2}, notice that 
 every  fundamental solution  $\Phi^\prime(\lambda,u) $  of system  \eqref{20gennaio2021-21}, defined at $u^\bullet$, is  
 $
 \Phi^\prime(\lambda,u)  =\Phi_{\rm hol}(\lambda,u^\bullet) C^\prime(u)$, with $ \det C^\prime(u)\neq 0$, 
   so that its $\hbox{\rm tr}(\mathcal{M}^\prime_j\mathcal{M}^\prime_k)$  coincide with $p_{jk}$ in \eqref{3febbraio2021-2}. 
 
   \bcr
 \label{18giugno2021-5}
 The results of Theorem \ref{16aprile2021-5} hold for the monodromy invariants  $p_{jk}$ of system  \eqref{23gennaio2021-8}  below, in terms of the Stokes matrices of system \eqref{systempainleve}, where $\Omega$ is as in Theorem \ref{20gennaio2021-2}.
 \ecr 
 
 \begin{proof} 
 It is elementary to check that the $2\times 2$ Fuchsian system 
 \eqref{20gennaio2021-21} 
 is equivalent to 
\be
\label{23gennaio2021-8}
 \frac{d \widetilde{\Phi}}{d\widetilde{\lambda}} =\left(\frac{\widetilde{\mathcal{A}_1}(x)}{\widetilde{\lambda}}+\frac{\widetilde{\mathcal{A}_2}(x)}{\widetilde{\lambda}-x}+\frac{\widetilde{\mathcal{A}_3}(x)}{\widetilde{\lambda}-1} \right) \widetilde{\Phi},\quad\quad x=\frac{u_2-u_1}{u_3-u_1},
 \ee
 through the  gauge transformation $\widetilde{\Phi}(\widetilde{\lambda},x)=(u_3-u_1)^{\mathcal{A}_\infty}\Phi(\lambda,u)$  and change of variables
 $ \widetilde{\lambda}=(\lambda-u_1)/(u_3-u_1)$. 
 The Schlesinger system \eqref{23gennaio2021-1} (equivalent to PVI) expresses isomonodromy of \eqref{23gennaio2021-8}. 
 Now, we just take  $u_1=0$, $u_2=x$, $u_3=1$ in the above discussion.  
\end{proof}

\subsection{Proof of Theorem \ref{16aprile2021-5}}

 To the  selected vector solutions $\vec{\Psi}_1,\vec{\Psi}_2,\vec{\Psi}_3$, we  associate the   solutions $\vec{X}_k(\lambda,u):=G^{-1} \vec{\Psi}_k$, $k=1,2,3$,  of system \eqref{20gennaio2021-5}. They have the same monodromy \eqref{24gennaio2021-1}. Let 
\be
 \label{4febbraio2021-2} 
 X(\lambda,u):=G^{-1}\Psi(\lambda,u)=\Bigl( \vec{X}_1(\lambda,u)~\Bigl|~  \vec{X}_2(\lambda,u)~\Bigr|~ \vec{X}_3(\lambda,u) \Bigr)
  =\begin{pmatrix}
 X_1^1 & X_2^1 & X_3^1
 \\
 \noalign{\medskip}
 X_1^2 & X_2^2 & X_3^2
 \\
 \noalign{\medskip}
 X_1^3 & X_2^3 & X_3^3
 \end{pmatrix}
\ee
be the matrix solution of system \eqref{20gennaio2021-5} corresponding to \eqref{4febbraio2021-1}. 
  From the entries in the first and third row we  extract six  matrix solutions (not necessarily fundamental) of  the 2-dimensional system \eqref{20gennaio2021-22}
  $$ 
  X^{[jk]}:=\begin{pmatrix}
  X_j^1 & X_k^1 
  \\
  \noalign{\medskip}
  X_j^3 & X_k^3
  \end{pmatrix},\quad [jk]=[12],~[13],~[23],~[21],~[31],~[32].
  $$ 
   Clearly, they reduce to $X^{[12]}$, $X^{[13]}$, $X^{[23]}$, because
   $ 
   X^{[kj]}=X^{[jk]}P$, $ P:=\begin{pmatrix} 0 & 1 \\ 1 & 0 \end{pmatrix}
   $. 
The column vectors  of $X^{[jk]}$ will be written in bold symbol:
$$ 
\boldsymbol{X}_j:=
\begin{pmatrix}
 X_j^1  
  \\
  \noalign{\medskip}
  X_j^3 
\end{pmatrix},
\quad\quad
 \hbox{ so that }\quad 
 X^{[jk]}(\lambda,u)=\Bigl(\boldsymbol{X}_j(\lambda,u)~\Bigr|~\boldsymbol{X}_k(\lambda,u)\Bigr).
 $$
 The monodromy properties \eqref{24gennaio2021-1} induce the following transformations:
\begin{align}
\label{26gennaio2021-2}
& \gamma_1: \quad\boldsymbol{X}_1\longmapsto e^{2\pi i \theta_1} \boldsymbol{X}_1, &&\boldsymbol{X}_2\longmapsto  \boldsymbol{X}_2+\alpha_1 c_{12} \boldsymbol{X}_1, && \boldsymbol{X}_3\longmapsto  \boldsymbol{X}_3+\alpha_1 c_{13} \boldsymbol{X}_1;
\\
\noalign{\medskip}
\label{26gennaio2021-3}
& \gamma_2: \quad \boldsymbol{X}_1\longmapsto  \boldsymbol{X}_1+\alpha_2 c_{21} \boldsymbol{X}_2   , && \boldsymbol{X}_2\longmapsto e^{2\pi i \theta_2} \boldsymbol{X}_2, && \boldsymbol{X}_3\longmapsto  \boldsymbol{X}_3+\alpha_2 c_{23} \boldsymbol{X}_2;
\\
\noalign{\medskip}
\label{26gennaio2021-4}
& \gamma_3: \quad \boldsymbol{X}_1\longmapsto  \boldsymbol{X}_1+\alpha_3 c_{31} \boldsymbol{X}_3 , &&\boldsymbol{X}_2\longmapsto  \boldsymbol{X}_2+\alpha_3 c_{32} \boldsymbol{X}_3, && \boldsymbol{X}_3\longmapsto e^{2\pi i \theta_3} \boldsymbol{X}_3;
\end{align}

\noindent
We distinguish two cases:

\begin{itemize}

\item[1)]  there are two linearly independent $\boldsymbol{X}_j$, $\boldsymbol{X}_k$, $j\neq k$;

\item[2)] all the $\boldsymbol{X}_j$ are linearly dependent.

\end{itemize}

\vskip 0.2 cm 
We  consider the linearly independent case first (it occurs for example if   \eqref{4febbraio2021-1} is a fundamental matrix solution). Without loss of generality, we  assume that $\boldsymbol{X}_1$ and $\boldsymbol{X}_3$ are linearly independent (otherwise, the analogous discussion holds  for an independent pair  $\boldsymbol{X}_j, \boldsymbol{X}_k$). For the loops  $\gamma_1,\gamma_2,\gamma_3$ we compute the monodromy matrices $\mathcal{M}_1,\mathcal{M}_2,\mathcal{M}_3$ of the  following fundamental matrix solution of system \eqref{20gennaio2021-21}
 \be
    \label{19giugno2021-7}
\Phi_{\rm hol}(\lambda,u):=\Phi^{[13]}(\lambda,u)=\prod_{j=1}^3 (\lambda-u_j)^{-\theta_j/2} X^{[13]}(\lambda,u).
\ee
From \eqref{26gennaio2021-2} and \eqref{26gennaio2021-4}  we receive
  $$
\mathcal{M}_1=\begin{pmatrix}
e^{\pi i \theta_1} & e^{-\pi i \theta_1}\alpha_1 c_{13} 
\\
\noalign{\medskip}
0 &e^{-\pi i \theta_1}
\end{pmatrix},
\quad\quad  \mathcal{M}_3=\begin{pmatrix}
e^{-\pi i \theta_3}& 0
\\
\noalign{\medskip}
e^{-\pi i \theta_3}\alpha_3 c_{31}  &e^{\pi i \theta_3}
\end{pmatrix},
$$
so  that 
$$
p_{13}=e^{-i\pi(\theta_1+\theta_3)}\alpha_1\alpha_3c_{13}c_{31} +2\cos(\pi(\theta_1-\theta_3)).
$$
Now,    $ 
\boldsymbol{X}_2=a\boldsymbol{X}_1+b\boldsymbol{X}_3
$ for some $a,b\in\mathbb{C}$, so that for the loop $\gamma_2$  the transformation  \eqref{26gennaio2021-3} yields
$$
\mathcal{M}_2=e^{-i\pi\theta_2}\begin{pmatrix}
1+\alpha_2 c_{21} a & \alpha_2 c_{23} a
\\
\noalign{\medskip}
\alpha_2 c_{21}b &1+\alpha_2 c_{23}b
\end{pmatrix}.
$$
In order to find $a,b$, we use that trace and determinant, being invariant by  conjugation, must be 
\be
\label{19giugno2021-8}
   \hbox{\rm Tr}  \mathcal{M}_2=2\cos\pi\theta_j,\quad \quad \det \mathcal{M}_2=1.
   \ee
   Using that  $\alpha_2=e^{2i\pi  \theta_2}-1$ for $\theta_2\not\in\mathbb{Z}$ and $\alpha_2=2\pi i $ for $\theta_2\in\mathbb{Z}$, both relations \eqref{19giugno2021-8} turn out to be equivalent to 
   \be
   \label{26gennaio2021-5}
   c_{21}a+c_{23}b=\left\{\begin{array}{cc}
   1, & \hbox{ if } \theta_2\not\in\mathbb{Z};
   \\
   \noalign{\medskip}
   0, &\hbox{ if } \theta_2\in\mathbb{Z}.
   \end{array}
   \right.
   \ee
      In order to find other consitions on $a$ and $b$, we consider the transformation of $\boldsymbol{X}_2$ in \eqref{26gennaio2021-2} along the  loop $\gamma_1$:
   \begin{align*}
   \boldsymbol{X}_2&~\longmapsto && \boldsymbol{X}_2+\alpha_1 c_{12} \boldsymbol{X}_1\equiv  a\boldsymbol{X}_1+b\boldsymbol{X}_3+\alpha_1 c_{12} \boldsymbol{X}_1,
 \\
    \boldsymbol{X}_2=a\boldsymbol{X}_1+b\boldsymbol{X}_3 &~\longmapsto && ae^{2i\pi \theta_1}\boldsymbol{X}_1+b(\boldsymbol{X}_3+\alpha_1 c_{13} \boldsymbol{X}_1).
    \end{align*}
   The above holds if and only if 
   $ a(e^{2i\pi\theta_1}-1 )+\alpha_1(bc_{13}-c_{12})=0
   $.
   Thus
\be
   \label{26gennaio2021-6}    c_{13}b=\left\{\begin{array}{cc}
   c_{12}-a, & \hbox{ if } \theta_2\not\in\mathbb{Z};
   \\
   \noalign{\medskip}
   c_{12}, &\hbox{ if } \theta_2\in\mathbb{Z}.
   \end{array}
   \right.
   \ee
    For the loop $\gamma_3$ in \eqref{26gennaio2021-4}:
     \begin{align*}
   \boldsymbol{X}_2&~\longmapsto && \boldsymbol{X}_2+\alpha_3 c_{32} \boldsymbol{X}_3\equiv  a\boldsymbol{X}_1+b\boldsymbol{X}_3+\alpha_3 c_{32} \boldsymbol{X}_3,
 \\
    \boldsymbol{X}_2=a\boldsymbol{X}_1+b\boldsymbol{X}_3 &~\longmapsto && a(\boldsymbol{X}_1+\alpha_3 c_{31} \boldsymbol{X}_3)+be^{2\pi i \theta_3} \boldsymbol{X}_3.
    \end{align*}   
   The above is true if and only if
   $
   b(e^{2i\pi \theta_3}-1)+\alpha_3(ac_{31}-c_{32})=0$,  
   namely 
   \be
   \label{26gennaio2021-7}    c_{31}a=\left\{\begin{array}{cc}
   c_{32}-b, & \hbox{ if } \theta_2\not\in\mathbb{Z};
   \\
   \noalign{\medskip}
   c_{32}, &\hbox{ if } \theta_2\in\mathbb{Z}.
   \end{array}
   \right.
   \ee
   For the loop $\gamma_2$ in \eqref{26gennaio2021-3}, the analogous procedure yields again \eqref{26gennaio2021-5}. 
   
   We are ready to  compute $\hbox{\rm Tr}(\mathcal{M}_1  \mathcal{M}_2)$.  In case $\boldsymbol{X}_1$ and $\boldsymbol{X}_2$ are also independent (i.e. $b\neq 0$), since $\hbox{\rm Tr }(  \mathcal{M}_1  \mathcal{M}_2)$ is  invariant by conjugation, we can compute it using  $\boldsymbol{X}_1$ and $\boldsymbol{X}_2$  as a basis, and this is done as above for the case  $\boldsymbol{X}_1,\boldsymbol{X}_3$, yielding 
$$
p_{12}=e^{-i\pi(\theta_1+\theta_2)}\alpha_1\alpha_2c_{12}c_{21} +2\cos(\pi(\theta_1-\theta_2)).
$$
 In case $\boldsymbol{X}_1$ and $\boldsymbol{X}_2$ are not  independent, then 
 $b=0$ and 
 $$
\mathcal{M}_2
=
\begin{pmatrix}
e^{-i\pi\theta_2}(1+\alpha_2 c_{21} a )& e^{-i\pi\theta_2}\alpha_2 c_{23} a 
\\
\noalign{\medskip}
0 &e^{-i\pi\theta_2}
\end{pmatrix}.
$$
From  \eqref{26gennaio2021-5}, \eqref{26gennaio2021-6}, \eqref{26gennaio2021-7}  we receive 
$$
c_{21}a=\left\{\begin{array}{cc}
   1, & \hbox{ if } \theta_2\not\in\mathbb{Z},
   \\
   \noalign{\medskip}
   0, &\hbox{ if } \theta_2\in\mathbb{Z},
   \end{array}
   \right.
   \quad  0=\left\{\begin{array}{cc}
   c_{12}-a, & \hbox{ if } \theta_2\not\in\mathbb{Z},
   \\
   \noalign{\medskip}
   c_{12}, &\hbox{ if } \theta_2\in\mathbb{Z},
   \end{array}
   \right.
\quad \quad \quad c_{31}a=c_{32}.
   $$
   
   For $\theta_2\not\in\mathbb{Z}$:
   $$a=c_{12}=\frac{1}{c_{21}}, \quad\quad \hbox{ and }\quad c_{12}c_{21}=1,\quad c_{32}=c_{31}c_{12}.
   $$
           Using $\alpha_2=\exp\{2\pi i \theta_2\} -1$, we receive 
   $$
\mathcal{M}_2
=
\begin{pmatrix}
e^{i\pi\theta_2}& 2i\sin(\pi\theta_2) c_{12}c_{23} 
\\
\noalign{\medskip}
0 &e^{-i\pi\theta_2}
\end{pmatrix},
$$
so that 
$$
 p_{12}=2\cos(\pi(\theta_1+\theta_2))
\underset{c_{12}c_{21}=1}\equiv2c_{12}c_{21} (\cos \pi(\theta_1+\theta_2)-\cos\pi(\theta_1-\theta_2))+2\cos\pi(\theta_1-\theta_2)
$$

    For $\theta_2\in\mathbb{Z}$:
   $$ 
   ac_{21}=0,\quad c_{12}=0.
   $$
   Recalling that  $\alpha_2=2\pi i$, we receive  
   $$
\mathcal{M}_2
=
(-1)^{\theta_2}\begin{pmatrix}
1& 2\pi i c_{23}a  
\\
\noalign{\medskip}
0 &1
\end{pmatrix},\quad c_{31}a=c_{32}.
$$
   Therefore
$$
 p_{12}=(-1)^{\theta_2} 2\cos(\pi\theta_1)
\underset{c_{12}=0}\equiv2c_{12}c_{21} (\cos \pi(\theta_1+\theta_2)-\cos\pi(\theta_1-\theta_2))+2\cos\pi(\theta_1-\theta_2).
$$

  The computation of $p_{23}=\hbox{\rm Tr}(\mathcal{M}_2\mathcal{M}_3)$ can be done in an analogous way. 
   In conclusion,  all the possibilities considered are summarized  in the   formula   
   \be
 \label{3febbraio2021-1} 
   p_{jk}=e^{-i\pi(\theta_j+\theta_k)}\alpha_j\alpha_kc_{jk}c_{kj} +2\cos(\pi(\theta_j-\theta_k)).
    \ee
     Finally, we substitute \eqref{27gennaio2021-1} and receive \eqref{3febbraio2021-2} in full generality.

   \vskip 0.2 cm 
   $\bullet$ We consider the  linearly-dependent case.   We do a gauge transformation 
   $
   {}_\gamma Y:=z^{-\gamma}Y$, $\gamma\in\mathbb{C}$,
   which transforms  \eqref{20gennaio2021-4} into 
$$
d( {}_\gamma Y)=\omega(z,u;\gamma)~ {}_\gamma Y,\quad\quad \omega(z,u;\gamma)=\left(U+\frac{\widetilde{V}-\gamma I}{z}\right)dz +\sum_{k=1}^3 (zE_k+V_k)du_k,$$
where $ \widetilde{V}:=V-I$. 
This changes  $\theta_k\mapsto \theta_k+\gamma$, while $\theta_\infty$ is unchanged, and 
\eqref{5febbraio2021-1} changes to 
  $\mu_1[\gamma]= (\theta_\infty - \theta_1 - \theta_2 - \theta_3-3\gamma)/2$, $\mu_3[\gamma]=(-\theta_\infty - \theta_1 -\theta_2 - \theta_3-3\gamma)/2$.
  
  There exists  $\gamma_0>0$  sufficiently small such that      $V-\gamma$ and $\widetilde{V}-\gamma$ have non-integer eigenvalues and non-integer diagonal entries for 
  $
   0<|\gamma|<\gamma_0
  $.  
Hence,   the analogous of the matrix  \eqref{4febbraio2021-1},  here called     ${}_\gamma\Psi(\lambda,u)=\bigl( {}_\gamma\vec{\Psi}_1(\lambda,u)~\bigl|~  {}_\gamma\vec{\Psi}_2(\lambda,u)~\bigr|~ {}_\gamma\vec{\Psi}_3(\lambda,u) \bigr)$,   is fundamental. The connection coefficients in \eqref{24gennaio2021-1}, depending on $\gamma$, will be called   $c_{jk}[\gamma]$,   with  
 $$ 
 \alpha_k[\gamma]=e^{2\pi i (\theta_k+\gamma)}-1.
 $$
  Consequently, there are two independent column vectors in the triple  ${}_\gamma\boldsymbol{X}_1(\lambda,u),{}_\gamma\boldsymbol{X}_2(\lambda,u),{}_\gamma\boldsymbol{X}_3(\lambda,u)$.   
The discussion of the independent case  can be repeated, with the $c_{jk}[\gamma]$ and $ \alpha_k[\gamma] $ in the transformations \eqref{26gennaio2021-2}-\eqref{26gennaio2021-4}. 
 We can  assume that ${}_\gamma\boldsymbol{X}_1$ and ${}_\gamma\boldsymbol{X}_3$ are independent, otherwise the  discussion is analogous  for another independent pair  ${}_\gamma\boldsymbol{X}_j$, ${}_\gamma\boldsymbol{X}_k$. After the gauge
  $$ 
{}_\gamma{\Phi}=\prod_{j=1}^3 (\lambda-u_j)^{-(\theta_j+\gamma)/2}~{}_\gamma\boldsymbol{X},
$$
we obtain the analogous of system \eqref{20gennaio2021-21}:
\be
\label{5febbraio2021-2}
 \frac{d({}_\gamma\Phi)}{d\lambda} =\sum_{k=1}^3 \frac{\mathcal{A}_k[\gamma]}{\lambda-u_k} ~{}_\gamma\Phi, \quad
  \quad
   \mathcal{A}_k:=A_k[\gamma]-\frac{\theta_k+\gamma}{2}.
   \ee 
   Let  
   $$
 {}_\gamma  \Phi^{[1,3]}(\lambda, u)=\prod_{j=1}^3 (\lambda-u_j)^{-\theta_j/2} \Bigl( {}_\gamma\boldsymbol{X}_1(\lambda,u)~\Bigr|~ {}_\gamma\boldsymbol{X}_3(\lambda,u)\Bigr) 
   $$ 
   be the analogous of \eqref{19giugno2021-7}, and let  $\mathcal{M}_j[\gamma]$ be its monodromy matrix at $u_j$.  
The same procedure which has proved \eqref{3febbraio2021-1} yields
\begin{equation}
\label{19giugno2021-9}
\begin{aligned}
   p_{jk}[\gamma]
   &:=\hbox{\rm Tr}(\mathcal{M}_j[\gamma]\mathcal{M}_k[\gamma])
 \\
 \noalign{\medskip}
 &  =e^{-i\pi(\theta_j+\theta_k+2\gamma)}\alpha_j[\gamma]\alpha_k[\gamma]c_{jk}[\gamma]c_{kj}[\gamma] +2\cos(\pi(\theta_j-\theta_k)),\quad j\neq k.
  \end{aligned}
  \end{equation}
In general the $\gamma$-dependent objects $c_{jk}[\gamma]$,  ${}_\gamma\Psi$,  ${}_\gamma\boldsymbol{X}_k$ and ${}_\gamma  \Phi^{[1,3]}$ diverge for $\gamma\to 0$. Therefore, the  monodromy matrices $\mathcal{M}_1[\gamma]$,  $\mathcal{M}_2[\gamma]$, $\mathcal{M}_3[\gamma]$ generate the monodromy group for $0<|\gamma|<\gamma_0$, but may be not defined at $\gamma=0$. 
To overcome the problem, we use a relation  proved in full generality  in \cite{guz2016}, and in \cite{BJL4} in a generic case. In case $\mathbb{D}=\mathbb{D}(u^0)$, the relation says that at any $u\in\mathbb{D}(u^0)$ 
\be
\label{19giugno2021-10}
\alpha_k c_{jk}=\begin{cases}
\begin{aligned} 
e^{-2\pi i \gamma}\alpha_k[\gamma]~c_{jk}[\gamma], &&\hbox{ if }k\succ j,
\\[1ex]
\alpha_k[\gamma]~c_{jk}[\gamma] ,&&\hbox{ if } k\prec j,
\end{aligned}
\end{cases}
\quad \quad \hbox{for real  } 0<\gamma<\gamma_0.
\ee
The ordering $ k\prec j$ is \eqref{16aprile2021-3}, well defined and the same at any  $u\in\mathbb{D}(u^0)$.  
In case $\mathbb{D}=\mathbb{D}(u^c)$, the same relation holds at any $u\in \mathbb{D}(u^c)\backslash   \Delta_{\mathbb{C}^3}$   for $j\neq k$ such that $u_j^c\neq u_k^c$, the ordering relation \eqref{16aprile2021-4} being well defined  for every $u\in \mathbb{D}(u^c)\backslash   \Delta_{\mathbb{C}^3}$. For  $j\neq k$ such that $u_j^c= u_k^c$ the ordering relation is not defined, but $
c_{jk}=c_{jk}[\gamma]=0
$, so that we can  state that   \eqref{19giugno2021-10} still holds. Using \eqref{19giugno2021-10}, \eqref{19giugno2021-9} becomes
\be
\label{21aprile2021-1}
 p_{jk}[\gamma]=e^{i \pi  \gamma}~e^{-i\pi(\theta_j+\theta_k)}\alpha_j\alpha_kc_{jk}c_{kj} +2\cos(\pi(\theta_j-\theta_k)),\quad\quad 0<\gamma<\gamma_0 \hbox{ real},
\ee
 for both $u\in \mathbb{D}(u^0)$ and   $u\in\mathbb{D}(u^c)\backslash   \Delta_{\mathbb{C}^3}$ (in the latter case,  \eqref{21aprile2021-1} is true  also for $j\neq k$ such that $u_j^c= u_k^c$, because it just reduces to  the identity $2\cos(\pi(\theta_j-\theta_k))=2\cos(\pi(\theta_j-\theta_k))$). 
 Since both the  $c_{jk}$ and $c_{jk}[\gamma]$ are constant, \eqref{21aprile2021-1} extends analytically  at $ \Delta_{\mathbb{C}^3}$.

\vskip 0.2 cm 
The r.h.s. of  \eqref{21aprile2021-1}  depends holomorphically on $\gamma\in\mathbb{C}$, while the l.h.s. $p_{jk}[\gamma]$ has been computed from  ${}_\gamma \Phi^{[13]}$,  which  is not in general  defined for $\gamma=0$.
We show that $p_{jk}[\gamma]$ can also be obtained from a fundamental matrix solution of \eqref{5febbraio2021-2} which is  holomorphic of $\gamma$ in a neighbourhood of  $\gamma=0$.
 To do that, we need to recall  from \cite{CDG} that if  $\mathbb{D}=\mathbb{D}(u^c)$,  the choice of an  admissible direction  $\tau$ determines a cell decomposition of $\mathbb{D}(u^c)$ into topological cells, called {\bf $\tau$-cells}. They  are  the connected components of $\mathbb{D}(u^c)\backslash( \Delta_{\mathbb{C}^3}\cup X(\tau))$, where  $X(\tau)$ is the locus of points $u=(u_1,u_2,u_3)\in\mathbb{D}(u^c)$ such that $\Re(e^{i\tau}(u_j-u_k))=0$. 

For $\mathbb{D}= \mathbb{D}(u^0)$, let $u^\bullet\in \mathbb{D}(u^0)$. For $\mathbb{D}= \mathbb{D}(u^c)$, let $u^\bullet $ belong to a $\tau$-cell of $\mathbb{D}(u^c)$.     Then, there is a sufficiently small neighbourhood $\mathcal{U}$ of  $u^\bullet$ such that, as $u$ varies in $\mathcal{U}$, then   $u_k$ represented in the $\lambda$-plane remains    inside a closed disc $D_k$ centered at $u^\bullet_k$,  with $D_j\cap D_k=\emptyset$ for $1\leq j\neq k\leq 3$. Consider  the simply connected domain 
$$
\mathcal{B}_{u^\bullet}:=\mathcal{P}_{\eta^*}(u^\bullet)\backslash (D_1\cup D_2\cup D_3).$$ 
   Since system \eqref{5febbraio2021-2}  holomorphically  depends on the  parameters $(u,\gamma)\in \mathcal{U}\times  \{\gamma\in\mathbb{C} \hbox{ s.t. } |\gamma|<\gamma_0\}$, according to a general result (see for example \cite{Wasow}) it has a fundamental matrix solution  $$\Phi_{\rm hol}^{(u^\bullet)}(\lambda,u,\gamma)$$ holomorphic of  
$
(\lambda,u,\gamma)\in 
\mathcal{B}_{u^\bullet}
\times
 \mathcal{U}
 \times 
 \{\gamma\in\mathbb{C} \hbox{ s.t. } |\gamma|<\gamma_0\}
$. 
If $\mathcal{U}$ is sufficiently small, it    holomorphically extends  to
$
\bigl(\mathcal{P}_{\eta^*}(u)\widehat{\times}~ \mathcal{U}\bigr) \times \{\gamma\in\mathbb{C} \hbox{ s.t. } |\gamma|<\gamma_0\}$.  
 For some invertible connection matrix $C(u,\gamma)$ we have 
  $$ 
  \Phi_{\rm hol}^{(u^\bullet)}(\lambda,u,\gamma)={}_\gamma \Phi^{[13]}(\lambda,u)\cdot C(u,\gamma), \quad 0<|\gamma|<\gamma_0.
  $$ 
Now, $C(u,\gamma)$ is holomorphic of $u\in  \mathcal{U} $, for any  $0<|\gamma|<\gamma_0$, but   may diverge as $\gamma\to 0$. 
 On the other hand, the monodromy matrix $\mathcal{M}_k^{\rm hol}(u,\gamma)$ of $\Phi_{\rm hol}^{(u^\bullet)}(\lambda,u,\gamma)$   at $\lambda=u_k$,\footnote{ For example represented by the loop going around $\partial D_k$.}   is holomorphic of $(u,\gamma)\in \mathcal{U}\times  \{\gamma\in\mathbb{C} \hbox{ s.t. } |\gamma|<\gamma_0\}$, including $\gamma=0$. Moreover, 
  $$ 
  \mathcal{M}_k[\gamma]=C(u,\gamma) \cdot \mathcal{M}_k^{\rm hol}(u,\gamma) \cdot C(u,\gamma)^{-1}.  
  $$
     Since
      $\hbox{\rm Tr}(\mathcal{M}_j^{\rm hol}(u,\gamma)\mathcal{M}_k^{\rm hol}(u,\gamma))=\hbox{\rm Tr}(\mathcal{M}_j[\gamma]\mathcal{M}_k[\gamma])\equiv p_{jk}[\gamma]$, 
      we see from \eqref{21aprile2021-1}  that  
 $$
 \hbox{\rm Tr}(\mathcal{M}_j^{\rm hol}(u,\gamma)\mathcal{M}_k^{\rm hol}(u,\gamma))  =
 e^{i \pi  \gamma}~e^{-i\pi(\theta_j+\theta_k)}\alpha_j\alpha_kc_{jk}c_{kj} +2\cos(\pi(\theta_j-\theta_k)),\quad\quad 0<\gamma<\gamma_0. 
$$
 Now,  both $ \hbox{\rm Tr}(\mathcal{M}_j^{\rm hol}(u,\gamma)\mathcal{M}_k^{\rm hol}(u,\gamma))$ and $p_{jk}[\gamma]$ are continuous of $\gamma$ in a neighbourhood of $\gamma=0$. 
  Therefore  we are allowed to take the limit $\gamma\to 0_{+}$ and obtain 
   \be
   \label{21aprile2021-2}
   p_{jk}:=\hbox{\rm Tr}(\mathcal{M}_j^{\rm hol}(u,0)\mathcal{M}_k^{\rm hol}(u,0)) = ~e^{-i\pi(\theta_j+\theta_k)}\alpha_j\alpha_kc_{jk}c_{kj} +2\cos(\pi(\theta_j-\theta_k)).
   \ee
 In case $\mathbb{D}=\mathbb{D}(u^c)$, we can repeat the above discussion also   for  $u^\bullet$ on the boundary of one $\tau$-cell (provided that $u^\bullet\not \in  \Delta_{\mathbb{C}^3}$), because  $\tau=3\pi/2-\eta$ can be slightly changed   without affecting the properties of fundamental solutions. 
 
Now, also at $\gamma=0$ the matrix
  $$\Phi_{\rm hol}^{(u^\bullet)}(\lambda,u,0)$$ is a fundamental  solution of  system \eqref{20gennaio2021-21}. 
  Therefore, we have proved that for every $u^\bullet\in \mathbb{D}(u^0)$ or  $u^\bullet \in\mathbb{D}(u^c)\backslash  \Delta_{\mathbb{C}^3}$,  we can find a fundamental solution $\Phi_{\rm hol}^{(u^\bullet)}(\lambda,u,0)$  
 of system  \eqref{20gennaio2021-21},  holomorphically 
 depending
  on    $(\lambda,u)\in\mathcal{P}_{\eta^*}(u)\widehat{\times} ~\mathcal{U}$, with  $\mathcal{U}$  small enough, 
    such that its monodromy invariants are the $p_{jk}$  in \eqref{21aprile2021-2}. 
     Thus, we conclude that the formulae 
   \eqref{3febbraio2021-2} always 
    hold. 
    $\Box$
    
    \vskip 0.2 cm 
    Notice that in the linearly dependent case,     $\Phi_{\rm hol}(\lambda,u)$  in the statement of the theorem is precisely $\Phi_{\rm hol}^{(u^\bullet)}(\lambda,u,0)$ above, while in the linearly independent case it is a fundamental matrix like \eqref{19giugno2021-7}.     

 %%%%%%%%%%%%%%%%%%%%%

\section{Classification of  analytic solutions of PVI   satisfying theorem 1.1 of \cite{CDG}. Reduction to special functions}
\label{27luglio2021-3}

We study system \eqref{systempainleve} whose  $\Omega(x)$ is associated, through Theorem \ref{20gennaio2021-2}, with PVI transcendents that admit a Taylor expansion near a critical point. By the symmetries of PVI,  it suffices to only  consider  $x=0$. 
Our  first goal is to classify  those transcendents with behaviour $y(x)=\sum_{n=0}^\infty 
b_n x^n$ at $x=0$  such that  theorem 1.1  of \cite{CDG} applies, namely such that $\Omega(x)$ is holomorphic at $x=0$ and 
  \begin{equation}
\label{coalcond}
\lim_{x\to 0}\Omega_{12}(x)=\lim_{x\to 0}\Omega_{21}(x)=0 \quad \hbox{ holomorphically}.
\end{equation}

As already explained in the Introduction, since Jimbo's work \cite{Jimbo}, the strategy to solve the  non-linear connection problem for a transcendent $y(x)$   has been  to parameterize  the integration constants determining the critical behaviour at $x=0, 1 $ or $\infty$   in terms of $\theta_1$, $\theta_2$, $\theta_3$, $\theta_\infty$ and the traces $p_{jk}$ of the $2\times 2$ Fuchsian system \eqref{23gennaio2021-8} associated with  $y(x)$.  Theorem \ref{16aprile2021-5} and Corollary \ref{18giugno2021-5} allow us to 
obtain the $p_{jk}$'s in terms of the Stokes matrices of \eqref{systempainleve} (with $\Omega$ of Theorem \ref{20gennaio2021-2}).  According to theorem 
1.1  of \cite{CDG}, when \eqref{coalcond} holds  all the monodromy data of system \eqref{systempainleve} on a polydisc   
can be computed from the system at the fixed coalescence point $x=0$:
\begin{equation}
\label{mainsystem}
\frac{dY}{dz}=\left(U_0+\frac{\Omega_0}{z}\right)Y,\quad \quad U_0:=\mbox{diag}(0,0,1),\quad \Omega_0:=\Omega(0).
\end{equation}
The latter is simpler than  \eqref{systempainleve}, because $\Omega_{12}(0)=\Omega_{21}(0)=0$.

 The second goal of this section is to show that for the transcendents with Taylor series at $x=0$ such that \eqref{coalcond} holds,  system \eqref{mainsystem} is reduced to the confluent hypergeometric equation (equivalently, a Whittaker equation, and for special cases a Bessel equation) or  to the generalized hypergeometric equation of type $(q,p)=(2,2)$. This  implies that the Stokes matrices can be computed using  classical special functions. Some  selected examples of this computation will be the object of Section \ref{20giugno2021-1}. 

\vskip 0.2 cm 

The first column of the table below reproduces the table of \cite{guz2012} for transcendents with Taylor series at $x=0$, and corrects  minor misprints there.    The parameters  $\theta_1,\theta_2,\theta_3,\theta_\infty$   and  the free parameter (the integration constant) $y_0$ or $y_0^{(|N|)}$ in $y(x)$, if any, must satisfy the conditions in the second column. 
In the third column, we give our classification according to the fulfilment of the vanishing conditions \eqref{coalcond} and   indicate the classical special functions in terms of which \eqref{mainsystem} can be solved.  In the subsequent sub-sections we present  the details of this classification. For the classical special functions, we refer to Appendix \ref{specfunc}. For $N\in\mathbb{Z}\backslash\{0\}$  we  also define 
\be
\label{24giugno2021-1}
\mathcal{N}_N:= \begin{cases} \{0,2,4,\ldots,|N|-1\}\cup\{-2,-4,\ldots,-(|N|-1)\}, &\mbox{if } N\mbox{ is odd}\\ \{1,3,\ldots,|N|-1\}\cup\{-1,-3,\ldots,-(|N|-1)\}, &\mbox{if } N\mbox{ is even}.\end{cases}
\ee

\begin{longtable}{|p{4.9cm}|p{5cm}|l|}
\hline
\minitab[l]{\bf{Taylor series}} & \minitab[l]{\bf{Conditions on}\\\bf{parameters}} & \minitab[l]{\bf{Classical special}\\ \bf{functions}}\\
\hline
\minitab[l]{$\displaystyle{y(x)=y_0+\sum_{n\ge 1}b_n(y_0)~x^n}$} & \minitab[l]{{}
\\
$\displaystyle{y_0=\frac{\theta_\infty-1+\theta_3}{\theta_\infty-1}}$
\\
 $\theta_\infty+\theta_3\notin\mathbb{Z}$ 
 \\
  $\theta_\infty\ne 1$
  \\
   or
    \\
      $\displaystyle{y_0=\frac{\theta_\infty-1-\theta_3}{\theta_\infty}}$
       \\
         $\theta_\infty-\theta_3\notin\mathbb{Z}$\\ $\theta_\infty\ne 1$
         \\
         {} 
         } 
 & \multirow{2}*{\minitab[l]{{}\\The vanishing\\ condition \eqref{coalcond}\\ does not hold}} \\
 \cline{1-2} \minitab[l]{$y(x)=y_0+\displaystyle\sum_{n\ge 1}b_n(y_0)~x^n$} & \minitab[l]{{}\\ $\theta_3=0,\,\theta_\infty=1$ \\ $y_0$ free parameter\\{}} &\\
\hline
\multirow{2}*{\minitab[l]{ 
$\displaystyle{y(x)=\sum_{n=0}^{|N|-1}b_nx^n+\frac{y_0^{(|N|)}}{(|N|)!}x^{|N|}}$
\\
\noalign{\medskip} 
 $\quad\quad   \displaystyle{+\sum_{n\ge |N|+1}b_n(y_0^{(|N|)})~x^n }$ \\  (T1)}}
  &
   \minitab[l]{{}\\$\displaystyle{b_0=\frac{N}{\theta_\infty-1}},~\theta_\infty\ne 1.$ 
  \\
   $y_0^{(|N|)}
    \mbox{ free parameter}.$
    \\
    \noalign{\medskip}
     $\theta_\infty-1+\theta_3=N\in\mathbb{Z}\setminus\{0\}$ 
     \\ 
     or
       \\ $\theta_\infty-1-    \theta_3=N\in\mathbb{Z}\setminus\{0\}$, 
       \\
       \noalign{\medskip}
      and 
       \\ 
        $\theta_1-\theta_2 \in \mathcal{N}_N$\\{}}
 & 
    \minitab[l]{{\bf Confluent}\\ {\bf hypergeometric} \\ {}}\\
\cline{2-3} & \minitab[l]{{}
\\
$\displaystyle{b_0=\frac{N}{\theta_\infty-1}},~\theta_\infty\ne 1.$
\\
 $y_0^{(|N|)}\mbox{ free parameter.}$ 
 \\
 \noalign{\medskip}
  $\theta_\infty-1+\theta_3=N\in\mathbb{Z}\setminus\{0\}$ 
  \\
   or
    \\
     $\theta_\infty-1-\theta_3=N\in\mathbb{Z}\setminus\{0\},$
      \\
       \noalign{\medskip}
   and either 
        \\
         $\theta_\infty-1\in\begin{cases} \{-1,\ldots,-|N|+1\}
          \\
           \{1,\ldots,|N|-1\} \end{cases}$
           \\
           or $  \theta_1+\theta_2 \in \mathcal{N}_N$
         \\
           {}} 
           &\minitab[l]{ The vanishing\\ condition \eqref{coalcond}\\ does not hold} \\
\hline
\multirow{2}*{\minitab[l]{$\displaystyle{y(x)=y_0^\prime x+\sum_{n\ge 2}b_n(y_0^\prime)~x^n}$ 
\\
{}
\\
 (T2)}} & \minitab[l]{$\displaystyle{y_0'=\frac{\theta_1}{\theta_1-\theta_2},\quad \theta_1-\theta_2\not\in\mathbb{Z}}$} &\minitab[l]
 {{}
 \\
 $(2,2)$ {\bf - Generalized} \\ {\bf hypergeometric}
 \\
 {}}
  \\
\cline{2-3} & \minitab[l]{$\displaystyle{y_0'=\frac{\theta_1}{\theta_1+\theta_2}},\quad \theta_1+\theta_2\not\in\mathbb{Z}$} &\minitab[l]{\\ For $\theta_2=0$ \\ case above \\ {} \\ Otherwise: \\ the vanishing\\ condition \eqref{coalcond}\\ does not hold} \\
\hline
\minitab[l]{{}\\$\displaystyle{y(x)=y_0'x+\sum_{n\ge 2}b_n(y_0')~x^n}$\\ {}\\ (T3)\\{}} & \minitab[l]{{}\\ $\theta_1=\theta_2=0$ \\ $y_0'$ free parameter\\{}} & \minitab[l]{{}\\ {\bf Confluent} \\ {\bf hypergeometric} \\{}} \\
\hline
\multirow{4}*{\minitab[l]{
$\displaystyle{y(x)=\sum_{n=1}^{|N|}b_nx^n}$
\\
\noalign{\medskip}
$\quad\quad\displaystyle{+\frac{y_0^{(|N|+1)}}{(|N|+1)!}x^{|N|+1}}$
 \\
 \noalign{\medskip}
  $\quad\quad\displaystyle{+\sum_{n\ge |N|+2}b_n(y_0^{(|N|+1)})~x^n}$
   \\
   {}
    \\ (T4)}
    } 
&
 \minitab[l]{{}\\$\theta_1-\theta_2=N\in\mathbb{Z}\setminus\{0\}$ \\ $\theta_1\in\begin{cases} \{-1,-2,\ldots,-|N|+1\} \\ \{1,2,\ldots,|N|-1\}  \end{cases}$ \\ $y_0^{(|N|+1)}$ free parameter \\{}} & \minitab[l]{{\bf $(2,2)$ - Generalized}\\ {\bf hypergeometric}} 
 \\
\cline{2-3} & \minitab[l]{{}\\$\theta_1-\theta_2=N\in\mathbb{Z}\setminus\{0\}$ \\ or \\ $\theta_1+\theta_2=N\in\mathbb{Z}\setminus\{0\}$ \\ and \\ $\theta_1=0,N,$ \\  $y_0^{(|N|+1)}$ free parameter \\{} } & \minitab[l]{{\bf Confluent} \\ {\bf hypergeometric}}\\
\cline{2-3} & \minitab[l]{{}\\$\theta_1-\theta_2=N\in\mathbb{Z}\setminus\{0\}$ \\ $\{(\theta_3+\theta_\infty-1),$\\ $(-\theta_3+\theta_\infty-1)\}\cap\mathcal{N}_N\ne\emptyset$ \\  $y_0^{(|N|+1)}$ free parameter\\{} } & \minitab[l]{{}\\{\bf $(2,2)$ - Generalized}\\ {\bf hypergeometric}\\{}}\\
\cline{2-3} & \minitab[l]{{}\\$\theta_1+\theta_2=N\in\mathbb{Z}\setminus\{0\}$ \\  $\theta_1\in\begin{cases} \{-1,-2,\ldots,-|N|+1\} \\ \{1,2,\ldots,|N|-1\} \end{cases}$ \\ or \\ $\{(\theta_3+\theta_\infty-1),$\\ $(-\theta_3+\theta_\infty-1)\}\cap\mathcal{N}_N\ne\emptyset$ \\  $y^{(|N|+1)}_0$ free parameter \\{}  } & \minitab[l]{The vanishing\\ condition \eqref{coalcond}\\ does not hold} \\
\hline
\end{longtable}
\bre
{\rm
In the table of  \cite{guz2012} there is a missprint. Corresponding to the branches (45), the correct condition is $\sqrt{-2\beta}\in \{-1,-2,N+1\}$ for $N<0$, and $\sqrt{-2\beta}\in \{1,2,N-1\}$ for $N>0$. In (61), the correct condition is 
$\sqrt{2\alpha}\in \{-1,-2,N+1\}$ for $N<0$, and $\sqrt{2\alpha}\in \{1,2,N-1\}$ for $N>0$. 
}
\ere

%%%%%%%%%%%%

\subsection{Transcendents with Taylor expansion  $y(x)=y_0+O(x)$}
\label{primorif}
We study the behaviour of $\Omega(x)$ as given in Theorem \ref{20gennaio2021-2} when $y(x)=y_0+O(x)$, where $O(x)=\sum_{n\ge 1}b_nx^n$ is a Taylor series in a neighborhood of $x=0$.  The functions $k_1$, $k_2$ can be written as
\begin{align*}
&k_1(x)=k_1^0(1+O(x))\frac{\sqrt{y}(x-1)^{(1+\theta_2)/2}}{\sqrt{y-1}x^{(1+\theta_2)/2}}x^{(\theta_1+\theta_3y_0/(y_0-1))/2},\\[1ex]
&k_2(x)=k_2^0(1+O(x))\frac{\sqrt{y-x}(x-1)^{\theta_2-\theta_3}}{\sqrt{y-1}x^{(1+\theta_1)/2}}x^{(\theta_2+\theta_3y_0/(y_0-1))/2},&k_1^0,k_2^0\in\mathbb{C}\backslash\{0\}.
\end{align*}
Hence, the off diagonal elemets of the matrix $\Omega(x)$ have the structure
\begin{align*}
\Omega_{12}(x)&=\omega_{12}(x)x^{\theta_1-\theta_2},
&\Omega_{13}(x)&=\omega_{13}(x)x^{(\theta_1+\theta_3y_0/(y_0-1)-\theta_2-1)/2},
\\
\Omega_{21}(x)&=\omega_{21}(x)x^{\theta_2-\theta_1},
&\Omega_{23}(x)&=\omega_{23}(x)x^{(\theta_2+\theta_3y_0/(y_0-1)-\theta_1-1)/2},
\\
\Omega_{31}(x)&=\omega_{31}(x)x^{-(\theta_1+\theta_3y_0/(y_0-1)-\theta_2-1)/2-1},
&\Omega_{32}(x)&=\omega_{32}(x)x^{-(\theta_2+\theta_3y_0/(y_0-1)-\theta_1-1)/2-1},
\end{align*}
where $\omega_{ij}(x)$ are holomorphic functions at $x=0$, explicitly computed from the formulae of Theorem \ref{20gennaio2021-2}. 

There are {\it three classes} of solutions $y(x)=y_0+O(x)$, obtained in \cite{guz2006} and  classified in the tables of \cite{guz2012}. The generic one is
\begin{equation}
\label{genericregular}
y(x)=y_0(\theta_\infty,\theta_3)+\sum_{n=1}^\infty b_n(\vec{\theta}) x^n,\quad \quad \theta_\infty\ne 1, \quad\vec{\theta}=(\theta_1,\theta_2,\theta_3,\theta_\infty),
\end{equation}
with two possibilities
\[ 
y_0=\frac{\theta_\infty-1+\theta_3}{\theta_\infty-1},
\quad 
\theta_\infty+\theta_3\notin\mathbb{Z};
\quad \quad \hbox{
or
}
\quad \quad 
y_0=\frac{\theta_\infty-1-\theta_3}{\theta_\infty-1},
\quad 
\theta_\infty-\theta_3\notin\mathbb{Z}.
\]
The coefficients $\{b_n\}_{n\ge 1}$ are uniquely determined when  $y_0$ is decided between the two possibilities. The other two classes consist of the following  one-parameter families of solutions:  the family
\begin{equation}
\label{nongeneric1}
y(x)=y_0+\frac{(1-y_0)(1+\theta_1^2-\theta_2^2)}{2}x+\sum_{n=2}^\infty b_n(y_0,\theta_1,\theta_2)x^n,\quad\theta_3=0,\quad\theta_\infty=1,
\end{equation}
where $y_0$ is a free parameter, and the family
\begin{equation}
\label{nongeneric2}
y(x)=\sum_{n=0}^{|N|-1}b_n(\vec{\theta})x^n+\frac{y_0^{(|N|)}}{(|N|)!}x^{|N|}+\sum_{n=|N|+1}^{\infty}b_n(y_0^{(|N|)},\vec{\theta})x^n,
\end{equation}
where $y_0^{(|N|)}$ is a free parameter, $N\in\mathbb{Z}\setminus\{0\}$, the leading order coefficient is $b_0=N/(\theta_\infty-1)$,  $\theta_\infty\ne 1$, and the conditions 
$$\theta_\infty-1+\theta_3=N\quad \hbox{ or }\quad \theta_\infty-1-\theta_3=N
$$
hold,  with either
\begin{equation}
\label{condizione1}
\theta_\infty-1\in \begin{cases} \{-1,-2,\ldots,N+1\} &\mbox{if }N<0 \\ \{1,2,\ldots,N-1\} &\mbox{if } N>0 \end{cases}
\end{equation}
or  (here,  $\mathcal{N}_N$ is \eqref{24giugno2021-1})
\begin{equation}
\label{condizione2}
\{(\theta_1+\theta_2),\,(\theta_1-\theta_2)\}\cap\mathcal{N}_{N}\ne\emptyset.
\end{equation}

\subsubsection{Generic case~\eqref{genericregular}}
 We show that the vanishing condition \eqref{coalcond} does not hold.
Let us consider  the solution \eqref{genericregular} with
\[
y_0=\frac{\theta_\infty-1-\theta_3}{\theta_\infty-1},\quad \theta_\infty-\theta_3\notin\mathbb{Z},\quad \theta_\infty\ne 1.
\]

If $\operatorname{Re}(\theta_1-\theta_2)=0$, then the condition \eqref{coalcond} 
is equivalent to $
\theta_1=\theta_2$ and $\theta_3=\theta_\infty$, 
so that $\theta_\infty-\theta_3=0$, which is a contradiction.

If $\operatorname{Re}(\theta_1-\theta_2)>0$, then $\lim_{x\to 0}\Omega_{12}(x)=0$ and for the requirement $\lim_{x\to 0}\Omega_{21}(x)=0$ to be fulfilled it is necessary that $\omega_{21}(0)=0$, that is (from the explicit formulae)
\begin{equation}
\label{vanishingcondition}
\theta_1=\theta_2+\theta_3-\theta_\infty,
\end{equation}
implying $\operatorname{Re}(\theta_3-\theta_\infty)>0$. The functions $\Omega_{31}(x)$ and $\Omega_{21}(x)$ are then
\begin{align*}
&\Omega_{31}(x)=\left(\frac{\omega''_{31}(0)}{2}x^2+O(x^3)\right)x^{-(\theta_3-\theta_\infty+1)}=x\left(\frac{\omega''_{31}(0)}{2}+O(x)\right)x^{-(\theta_3-\theta_\infty)},\\
&\Omega_{21}(x)=\left(\frac{\omega''_{21}(0)}{2}x^2+O(x^3)\right)x^{\theta_2-\theta_1}=x^2\left(\frac{\omega''_{21}(0)}{2}+O(x)\right)x^{-(\theta_3-\theta_\infty)}.
\end{align*}
The entry  (3,1)
\begin{equation}
\frac{d\Omega_{31}}{dx}=\frac{\Omega_{32}\Omega_{21}}{x(x-1)},
\end{equation}
  of the isomonodromic deformation equation \eqref{20gennaio2021-26}   gives 
\[
\frac{\omega''_{31}(0)}{2}(\theta_\infty-\theta_3+1)+O(x)=0,
\quad 
\hbox{ which implies }
\theta_\infty=\theta_3-1,
\]
giving a contradiction.

If $\operatorname{Re}(\theta_1-\theta_2)<0$,  than $\lim_{x\to 0}\Omega_{21}(x)=0$ and a necessary condition for $\Omega_{12}(x)$ to vanish as $x\to 0$ is $\omega_{12}(0)=0$, that is (from the explicit formulae)
\begin{equation}
\label{vanishingconditionneg}
\theta_1=\theta_2-\theta_3+\theta_\infty,
\end{equation}
which implies also $\operatorname{Re}(-\theta_3+\theta_\infty)<0$. With condition~\eqref{vanishingconditionneg} it follows that $\Omega_{32}(x)=\omega_{32}(x) x^{-(\theta_3-\theta_\infty+1)}$, with $\omega_{32}(0)=\omega'_{32}(0)=0$. Hence from the entry $(3,2)$ of the  isomonodromic deformation equation  \eqref{20gennaio2021-26}
\[
\frac{d\Omega_{32}}{dx}=(\theta_3-\theta_2)\frac{\Omega_{32}}{x-1}-\frac{\Omega_{12}\Omega_{31}}{x},
\]
we receive
\[
\frac{\omega''_{32}(0)}{2}(1-\theta_3+\theta_\infty)+O(x)=0,
\quad \hbox{ which implies }\theta_\infty=\theta_3-1,
\]
obtaining again a contradiction.
\noindent
Analogous results hold  for the solution with
\[
y_0=\frac{\theta_\infty-1+\theta_3}{\theta_\infty-1},
\quad 
\theta_\infty+\theta_3\notin\mathbb{Z},\quad \theta_\infty\ne 1.
\]
\noindent
We conclude that the vanishing condition \eqref{coalcond} does not hold.

\subsubsection{Solutions~\eqref{nongeneric1}}
If $\operatorname{Re}(\theta_1-\theta_2)=0$, then  the vanishing condition \eqref{coalcond}
holds if and only  $\theta_1=\theta_2$ and  $\omega_{12}(0)=\omega_{21}(0)=0$. From the explicit formulae  $\omega_{12}(0)=\omega_{21}(0)=0$ is equivalent to 
\[
\begin{cases}
1-\theta_1+\theta_2=0\\
-1-\theta_1+\theta_2=0
\end{cases}
\]
which is a contraddiction. If $\operatorname{Re}(\theta_1-\theta_2)>0$, the requirement $\lim_{x\to 0}\Omega_{21}(x)=0$ is satisfied only if $\omega_{21}(0)=0$, or equivalently (from the explicit formulae)  $
\theta_1-\theta_2=-1$, a contradiction. If $\operatorname{Re}(\theta_1-\theta_2)<0$, the condition $\lim_{x\to 0}\Omega_{12}(x)=0$ is satisfied only if   $\omega_{12}(0)=0$, or equivalently $
\theta_1-\theta_2=1$, 
a contradiction. 
We conclude that the vanishing condition \eqref{coalcond} does not hold.

\subsubsection{Solutions~\eqref{nongeneric2} - Case $\theta_\infty=\theta_3+N+1$} 

Let us start with the case $\operatorname{Re}(\theta_1-\theta_2)=0$. It is necessary for 
 \eqref{coalcond}  that $\theta_1=\theta_2$ and $\omega_{12}(0)=\omega_{21}(0)=0$, which  from the explicit formulae is equivalent to $\theta_1=\theta_2+N+1$ and $N=-1$, that is
\begin{equation}
\label{cccondition}
\theta_1=\theta_2\hbox{ and }
\theta_\infty=\theta_3.
\end{equation}
With this condition $\Omega(x)$ has holomorphic limit at $x=0$:
{\small

\[
\Omega_0=
\begin{pmatrix}
-\theta_2 & 0 & -\tilde{k}_1^0\sqrt{\theta_3}\\[2ex]
0 & -\theta_2 & \tilde{k}_2^0\sqrt{\theta_3}\\[2ex]
\displaystyle{\frac{1}{\tilde{k}_1^0\sqrt{\theta_3}}\left(\frac{\theta_3(1-\theta_2)}{2}+(1-\theta_3)y'_0\right)} & \displaystyle{\frac{1}{\tilde{k}_2^0\sqrt{\theta_3}}\left(\frac{\theta_3(1+\theta_2)}{2}+(1-\theta_3)y'_0\right)}  & -\theta_3
\end{pmatrix}
,
\]
}
where
\begin{equation}
\label{par1}
\tilde{k}_1^0=ik_1^0e^{i\pi\theta_2/2}\quad \tilde{k}_2^0=-k_2^0e^{i\pi(\theta_2-\theta_3)}.
\end{equation}
A diagonalizing matrix of $\Omega_0$ is
{\small
\begin{equation}
\label{g0ex1}
G_0=
\begin{pmatrix}
\displaystyle{-\frac{k_1}{k_2}\frac{\theta_3(1+\theta_2)+2(1-\theta_3)y'_0}{\theta_3(1-\theta_2)+2(1-\theta_3)y'_0}} & -\displaystyle{\frac{k_1\sqrt{\theta_3}}{\theta_2}} & \displaystyle{\frac{k_1}{\sqrt{\theta_3}}}\\[2ex]
1 & \displaystyle{\frac{k_2\sqrt{\theta_3}}{\theta_2}} & \displaystyle{-\frac{k_2}{\sqrt{\theta_3}}}\\[2ex]
0 & 1 & 1
\end{pmatrix}
,
\end{equation}
}
thus system~\eqref{mainsystem} 
for  $Y=G_0\widetilde{Y}$
 is
\begin{equation}
\label{sysex1}
\frac{d\widetilde{Y}}{dz}=\left[\frac{1}{\theta_2+\theta_3}
\begin{pmatrix}
0 & 0 & 0\\
0 & \theta_2 & \theta_2 \\
0 & \theta_3 & \theta_3
\end{pmatrix}
-\frac{1}{z}
\begin{pmatrix}
\theta_2 & 0 & 0\\
0 & 0 & 0\\
0 & 0 & \theta_2+\theta_3
\end{pmatrix}
\right]\widetilde{Y}.
\end{equation} 
If $(\tilde{y}_1,\tilde{y}_2,\tilde{y}_3)^T$ is a column of $\widetilde{Y}$, then
{\small 
\[
\begin{cases}
\displaystyle{\frac{d\tilde{y}_1}{dz}=-\frac{\theta_2}{z}\tilde{y}_1}\\[2ex]
\displaystyle{\frac{d\tilde{y}_2}{dz}=\frac{\theta_2}{\theta_2+\theta_3}\left(\tilde{y}_2+\tilde{y}_3\right)}\\[2ex]
\displaystyle{\frac{d\tilde{y}_3}{dz}=\frac{\theta_3}{\theta_2+\theta_3}\left(\tilde{y}_2+\tilde{y}_3\right)-\frac{\theta_2+\theta_3}{z}\tilde{y}_3}.
\end{cases}
\]
}
From the first  and second equations  we obtain
\begin{equation}
\label{y1ex1}
\tilde{y}_1(z)=Cz^{-\theta_2},\quad C\in\mathbb{C},\quad\quad \tilde{y}_3=\frac{\theta_2+\theta_3}{\theta_2}\frac{d\tilde{y}_2}{dz}-\tilde{y}_2,
\end{equation}
and plugging  this expression into the third equation we get  a {\bf confluent hypergeometric equation} with parameters $a=\theta_2$ and $b=\theta_2+\theta_3$:
\begin{equation}
\label{y2ex1}
z\frac{d^2\tilde{y}_2}{dz^2}+(\theta_2+\theta_3-z)\frac{d\tilde{y}_2}{dz}-\theta_2\tilde{y}_2=0.
\end{equation}

If $\operatorname{Re}(\theta_1-\theta_2)>0$, then $\lim_{x\to 0}\Omega_{12}(x)=0$ holomorphically for $\theta_1-\theta_2$ positive integer,  and for the requirement $\lim_{x\to 0}\Omega_{21}(x)=0$ holomorphically to be fulfilled it is necessary that $\omega_{21}(0)=0$, which is equivalent to
\[
\theta_1=\theta_2-N-1,\quad \hbox{ with } N\le -2.
\]
 Furthermore, from the explicit formulae, we have 
\[
\lim_{x\to 0}\Omega_{23}(x)=\omega_{23}(0)=\tilde{k}_2^0\sqrt{|N|\theta_3},
\quad \hbox{ and }\quad 
\omega_{32}(0)=0,\quad \omega_{31}(0)=0.
\]
Now $\Omega_{32}=\omega_{32}/x$, and so 
\[
\lim_{x\to 0}\Omega_{32}(x)=\left.\frac{d}{dx} \omega_{32}(x)\right|_{x=0}=\frac{\theta_2}{\tilde{k}_2^0}\sqrt{\frac{\theta_3}{|N|}}
\]
 The entries $(2,1)$ and $(3,1)$ of the isomonodromic deformation equations \eqref{20gennaio2021-26} give 
 \begin{equation}
\label{dv21dv31}
\begin{cases}
\displaystyle{\frac{d\Omega_{21}}{dx}=(\theta_2-\theta_1)\frac{\Omega_{21}}{x}+\frac{\Omega_{31}\Omega_{23}}{x-1}}\\[2ex]
\displaystyle{\frac{d\Omega_{31}}{dx}=\frac{\Omega_{21}\Omega_{32}}{x(x-1)}}
\end{cases}
\quad
\Longrightarrow
\quad 
\begin{cases}
\displaystyle{\frac{d\omega_{21}}{dx}=\frac{\omega_{31}\omega_{23}}{x(x-1)}}\\[2ex]
\displaystyle{\frac{d\omega_{31}}{dx}}+N\displaystyle{\frac{\omega_{31}}{x}=\frac{\omega_{21}\omega_{32}}{x(x-1)}.}
\end{cases}
\end{equation}
From the second equation  we  compute $\omega_{31}'(0)=\lim_{x\to 0}\omega_{31}(x)/x$ by taking its limit for $x\to 0$:
\begin{align*}
\omega'_{31}(0)-|N|\lim_{x\to 0}\frac{\omega_{31}(x)}{x}&=\lim_{x\to 0}\frac{\omega_{21}(x)\omega_{32}(x)}{x(x-1)}
\\
&
=\lim_{x\to 0}\frac{\left(\omega_{21}'(0)x+O(x^2)\right)\left(\omega_{32}'(0)x+O(x^2)\right)}{x}=0,
\end{align*}
so that we have  (recall that $N\le-2$)
$$\omega'_{31}(0)(1-|N|)=0
\quad \Longrightarrow\quad \omega'_{31}(0)=0.
$$ 
From the above and  first equation we obtain $\omega'_{21}(0)=0$, because
$$
\omega'_{21}(0)=\lim_{x\to 0}\frac{\omega_{31}(x)\omega_{23}(x)}{x(x-1)}=\lim_{x\to 0}\frac{\left(\omega_{31}'(0)x+O(x^2)\right)\left(\omega_{23}(0)+O(x)\right)}{x(x-1)}= -\omega_{31}^\prime(0)~ \tilde{k}_2^0\sqrt{|N|\theta_3}=0.
$$
 By differentiating $n$ times equations~\eqref{dv21dv31}, it follows that 
\[
\frac{d^{n+1}}{dx^{n+1}}(\omega_{21})=\sum_{j=0}^n\binom{n}{j}\frac{d^j}{dx^j}\left(\frac{\omega_{31}}{x}\right)\frac{d^{n-j}}{dx^{n-j}}\left(\frac{\omega_{23}(x)}{x-1}\right)
\]
and
\[
\frac{d^{n+1}}{dx^{n+1}}(\omega_{31})\left(1-\frac{|N|}{n+1}\right)=\sum_{j=0}^{n}\binom{n}{j}\frac{d^j}{dx^j}(\omega_{21}(x))\frac{d^{n-j}}{dx^{n-j}}\left(\frac{\omega_{32}(x)}{x(x-1)}\right).
\]
Thus, taking the limit $x\to 0$, we get
\[
\omega_{31}^{(n+1)}(0)\left(1-\frac{|N|}{n+1}\right)=\lim_{x\to 0}\sum_{j=0}^n\binom{n}{j}\frac{d^j}{dx^j}\left(\frac{\omega_{21}(x)}{x}\right)\frac{d^{n-j}}{dx^{n-j}}\left(\frac{\omega_{32}(x)}{x-1}\right),
\]
so that  we can repeat the computations till $1-|N|/(n+1)=0$, that is up to step $n=|N|-2$ included, and show that all the derivatives of $\omega_{21}(x)$ and $\omega_{31}(x)$ at $x=0$ up to the order $|N|-1$ vanish. In conclusion, 
\[
\omega_{21}(x)=x^{|N|}\left(\frac{\omega_{21}^{(|N|)}(0)}{(|N|)!}+O(x)\right),\quad\omega_{31}(x)=x^{|N|}\left(\frac{\omega_{31}^{(|N|)}(0)}{(|N|)!}+O(x)\right)
\]
and
\[
\lim_{x\to 0}\Omega_{21}(x)=0,\quad\lim_{x\to 0}\Omega_{31}(x)=\sqrt{\frac{\theta_3}{|N|}}\frac{K(\theta_2,\theta_3,N,A)}{\tilde{k}_1^0} \quad\hbox{ holomorphically},
\]
where 
\[
K(\theta_2,\theta_3,N,A)=\left.-\frac{1}{2}\frac{d^{|N|}}{dx^{|N|}}\left[\frac{x(x-1)dy/dx+\theta_3y(y-x)}{y-1}\right]\right|_{x=0}+|N|A,\quad\quad A:=\frac{y_0^{(|N|)}}{2}.
\]
Here $A$ is  the free parameter. 
The matrix $\Omega(x)$ at the coalescence point $x=0$  is then
{\small \[
\Omega_0=
\begin{pmatrix}
1-|N|-\theta_2 & 0 & 0\\[2ex]
0 & -\theta_2 & \tilde{k}_2^0\sqrt{|N|\theta_3}\\[2ex]
\displaystyle{\sqrt{\frac{\theta_3}{|N|}}\frac{K(\theta_2,\theta_3,N,A)}{\tilde{k}_1^0}} & \displaystyle{\frac{\theta_2}{\tilde{k}_2^0}\sqrt{\frac{\theta_3}{|N|}}} & -\theta_3
\end{pmatrix}
,
\]
}
where $\tilde{k}_1^0,\tilde{k}_2^0$ are as in~\eqref{par1}. A diagonalizing matrix for $\Omega_0$ is
{\small 
\[
G_0=
\begin{pmatrix}
\displaystyle{\frac{k_1}{K(\theta_2,\theta_3,N,A)}\sqrt{\frac{|N|}{\theta_3}}\left(\frac{\theta_2\theta_3}{|N|-1}+\theta_3-\theta_2-|N|+1\right)} & 0 & 0 \\[2ex]
\displaystyle{-k_2\frac{\sqrt{|N|\theta_3}}{|N|-1}} & 1& \displaystyle{-k_2\sqrt{\frac{|N|}{\theta_3}}}\\[2ex]
1 & \displaystyle{\frac{\theta_2}{k_2\sqrt{|N|\theta_3}}}  & 1
\end{pmatrix}
,
\]
}
thus, applying the gauge transformation $Y=G_0\widetilde{Y}$ to system~\eqref{mainsystem} we get
\[
\frac{d\widetilde{Y}}{dz}=\left[\frac{1}{\theta_2+\theta_3}
\begin{pmatrix}
0 & 0 & 0\\
\tilde{k}_2^0\sqrt{|N|\theta_3} & \theta_2 & \tilde{k}_2^0\sqrt{|N|\theta_3} \\
\theta_3 & \displaystyle{\frac{\theta_2}{\tilde{k}_2^0}\sqrt{\frac{\theta_3}{|N|}}} & \theta_3
\end{pmatrix}
-\frac{1}{z}
\begin{pmatrix}
\theta_2+|N|-1 & 0 & 0\\
0 & 0 & 0\\
0 & 0 & \theta_2+\theta_3
\end{pmatrix}
\right]\widetilde{Y}.
\]
Let $(\tilde{y}_1,\tilde{y}_2,\tilde{y}_3)^T$ be a column of $\widetilde{Y}$, then
{\small \[
\begin{cases}
\displaystyle{\frac{d\tilde{y}_1}{dz}=-\frac{\theta_2+|N|-1}{z}\tilde{y}_1}\\[2ex]
\displaystyle{\frac{d\tilde{y}_2}{dz}=\frac{1}{\theta_2+\theta_3}\left(\tilde{k}_2^0\sqrt{|N|\theta_3}\tilde{y}_1+\theta_2\tilde{y}_2+\tilde{k}_2^0\sqrt{|N|\theta_3}\tilde{y}_3\right)}\\[2ex]
\displaystyle{\frac{d\tilde{y}_3}{dz}=\frac{1}{\theta_2+\theta_3}\left(\theta_3\tilde{y}_1+\frac{\theta_2}{\tilde{k}_2}\sqrt{\frac{\theta_3}{|N|}}\tilde{y}_2+\theta_3\tilde{y}_3\right)-\frac{\theta_2+\theta_3}{z}\tilde{y}_3}.
\end{cases}
\]
}
The first and second  equations  give
\[
\tilde{y}_1(z)=C\,z^{1-|N|-\theta_2},\quad C\in\mathbb{C},
\quad
\quad 
\tilde{y}_3=\frac{\theta_2+\theta_3}{\tilde{k}_2^0\sqrt{|N|\theta_3}}\frac{d\tilde{y}_2}{dz}-\frac{\theta_2}{\tilde{k}_2^0\sqrt{|N|\theta_3}}\tilde{y}_2-Cz^{1-|N|-\theta_2}
\]
and inserting these expressions into the third equation we get
\[
z\frac{d^2\tilde{y}_2}{dz^2}+(\theta_2+\theta_3-z)\frac{d\tilde{y}_2}{dz}-\theta_2\tilde{y}_2-Lz^{1-(|N|+\theta_2)}=0,
\quad
\quad 
L:=C(1-|N|+\theta_3)\frac{\tilde{k}_2^0\sqrt{|N|\theta_3}}{\theta_2+\theta_3}.
\]
This is an inhomogeneous {\bf confluent hypergeometric equation} with parameters $a=\theta_2$ and $b=\theta_2+\theta_3$.

If $\operatorname{Re}(\theta_1-\theta_2)<0$, then $\lim_{x\to 0}\Omega_{21}(x)=0$ and for the requirement $\lim_{x\to 0}\Omega_{12}(x)=0$ to be fulfilled it is necessary that $\omega_{12}(0)=0$, which is equivalent to
$$
\theta_1=\theta_2+N+1.
$$
With this condition, $\operatorname{Re}(\theta_1-\theta_2)<0$ implies $N\le -2$ as before, thus
\begin{align*}
&\lim_{x\to 0}\Omega_{13}(x)=\omega_{13}(0)=-\tilde{k}_1^0\sqrt{|N|\theta_3},
& \lim_{x\to 0}\Omega_{23}(x)=\lim_{x\to 0}\left(\omega_{23}(x)x^{-(N+1)}\right)=0,\\
&\lim_{x\to 0}\Omega_{31}(x)=\left.\frac{d\omega_{31}}{dx}\right|_{x=0}=-\frac{\theta_2-|N|+1}{\tilde{k}_1^0}\sqrt{\frac{\theta_3}{|N|}}
\end{align*}
and $\omega_{32}(0)=0$. As we did before, from the isomonodromic deformations equations
\[
\displaystyle{\frac{d\Omega_{12}}{dx}=(\theta_1-\theta_2)\frac{\Omega_{12}}{x}-\frac{\Omega_{13}\Omega_{32}}{x-1}},
\quad\quad
\displaystyle{\frac{d\Omega_{32}}{dx}=(\theta_3-\theta_2)\frac{\Omega_{32}}{x-1}-\frac{\Omega_{12}\Omega_{31}}{x}},
\]
it is possible to show that 
\[
\left.\frac{d^n\omega_{12}}{dx^n}\right|_{x=0}=\left.\frac{d^n\omega_{32}}{dx^n}\right|_{x=0}=0
\]
for each $1\le n\le |N|-1$, thus
\begin{align*}
& \lim_{x\to 0}\Omega_{12}(x)=\lim_{x\to 0}\left(\omega_{12}(x)x^{1-|N|}\right)=0
\\
\noalign{\medskip}
&
\lim_{x\to 0}\Omega_{32}(x)=\lim_{x\to 0}\left(\omega_{32}(x)x^{-|N|}\right)=-\frac{H(\theta_2,\theta_3,N,A)}{\tilde{k}_2^0}\sqrt{\frac{\theta_3}{|N|}},
\end{align*}
where 
\[
H(\theta_2,\theta_3,N,A)=\left.\frac{1}{2}\frac{d^{|N|}}{dx^{|N|}}\left[\frac{x(x-1)dy/dx+\theta_3y(y-x)}{y-1}\right]\right|_{x=0}-|N|A,\quad \quad A:=\frac{y_0^{(|N|)}}{2}.
\]
Here,  $A$ is the free parameter. The matrix $\Omega(x)$ at the coalescence point $x=0$ is then
{\small\[
\Omega_0=
\begin{pmatrix}
-\theta_2+|N|-1 & 0 & -\tilde{k}_1^0\sqrt{|N|\theta_3}\\[2ex]
0 & -\theta_2 & 0\\[2ex]
-\displaystyle{\sqrt{\frac{\theta_3}{|N|}}\frac{\theta_2-|N|+1}{\tilde{k}_1^0}} & -\displaystyle{\sqrt{\frac{\theta_3}{|N|}}\frac{H(\theta_2,\theta_3,N,A)}{\tilde{k}_2^0}} & -\theta_3
\end{pmatrix}
,
\]
}
where $\tilde{k}_1^0,\tilde{k}_2^0$ are as in~\eqref{par1}. A diagonalizing matrix of $\Omega_0$ is
{\small
\[
G_0=
\begin{pmatrix}
\displaystyle{\frac{\tilde{k}_1^0\sqrt{|N|\theta_3}}{|N|-1}} & 1 & \displaystyle{\tilde{k}_1^0\sqrt{\frac{|N|}{\theta_3}}}\\[2ex]
\displaystyle{\frac{\tilde{k}_2^0\theta_2}{H(\theta_2,\theta_3,N,A)}\sqrt{\frac{|N|}{\theta_3}}\left(1-\frac{\theta_3}{|N|-1}\right)} & 0 & 0\\[2ex]
1 & \displaystyle{-\frac{\theta_2-|N|+1}{\tilde{k}_1^0\sqrt{|N|\theta_3}}} & 1
\end{pmatrix}
,
\]
}
hence applying the gauge transformation  to system~\eqref{mainsystem} we get
{\small
$$
\frac{d\widetilde{Y}}{dz}=\left[
\frac{1}{\theta_2+\theta_3-|N|+1}
\begin{pmatrix}
0 & 0 & 0\\
-\tilde{k}_1^0\sqrt{|N|\theta_3} & \theta_2-|N|+1 & -\tilde{k}_1^0\sqrt{|N|\theta_3} \\
\theta_3 & -\displaystyle{\frac{\theta_2-|N|+1}{\tilde{k}_1^0}\sqrt{\frac{\theta_3}{|N|}}} & \theta_3
\end{pmatrix}
 -\frac{1}{z}
\begin{pmatrix}
\displaystyle{\theta_2} & 0 & 0\\
0 & 0 & 0\\
0 & 0 & \displaystyle{\theta_2+\theta_3-|N|+1}
\end{pmatrix}
\right]
\widetilde{Y}.
$$
}
Let $(\tilde{y}_1,\tilde{y}_2,\tilde{y}_3)^T$ be a column of $\widetilde{Y}$, then
{\small \[
\begin{cases}
\displaystyle{\frac{d\tilde{y}_1}{dz}=-\frac{\theta_2}{z}\tilde{y}_1}\\[2ex]
\displaystyle{\frac{d\tilde{y}_2}{dz}=\frac{1}{\theta_2+\theta_3-|N|+1}\left(-\tilde{k}_1^0\sqrt{|N|\theta_3}\tilde{y}_1+(\theta_2-|N|+1)\tilde{y}_2-\tilde{k}_1^0\sqrt{|N|\theta_3}\tilde{y}_3\right)}\\[2ex]
\displaystyle{\frac{d\tilde{y}_3}{dz}=\frac{1}{\theta_2+\theta_3-|N|+1}\left(\theta_3\tilde{y}_1-\frac{\theta_2-|N|+1}{\tilde{k}_1^0}\sqrt{\frac{\theta_3}{|N|}}\tilde{y_2}+\theta_3\tilde{y}_3\right)-\frac{\theta_2+\theta_3-|N|+1}{z}\tilde{y}_3}.
\end{cases}
\]
}
The first and second equations give
\[
\tilde{y}_1(z)=C\,z^{-\theta_2},\quad C\in\mathbb{C},
\quad\quad
\tilde{y}_3=\frac{\theta_2-|N|+1}{\tilde{k}_1^0\sqrt{|N|\theta_3}}\tilde{y}_2-Cz^{-\theta_2}-\frac{\theta_2+\theta_3-|N|+1}{\tilde{k}_1^0\sqrt{|N|\theta_3}}\frac{d\tilde{y}_2}{dz}
\]
and the third equation becomes
\[
z\frac{d^2\tilde{y}_2}{dz^2}+\left(\theta_2+\theta_3-|N|+1-z\right)\frac{d\tilde{y}_2}{dz}-\left(\theta_2-|N|+1\right)\tilde{y}_2=-Lz^{-\theta_2},
\quad
L:=\frac{C\tilde{k}_1^0\sqrt{|N|\theta_3} ~(\theta_3-|N|+1)}{\theta_2+\theta_3-|N|+1}.
\]
This is an inhomogeneous {\bf confluent hypergeometric equation} with parameters $a=\theta_2-|N|+1$ and $b=\theta_2+\theta_3-|N|+1$.

\subsubsection{Solutions~\eqref{nongeneric2} - Case $\theta_\infty=-\theta_3+N+1$} 
For this case, the tecniques to find the solutions for which condition~\eqref{coalcond} is satisfied  are the same as those used in the previous sections, so we will just report the main results.

For $\operatorname{Re}(\theta_1-\theta_2)=0$, the matrix $\Omega(x)$ holomorphically  has  limit $\Omega_0$ as $x\to 0$ with vanishing $\Omega_{12}(0)=\Omega_{21}(0)=0$ if and only if $\theta_1=\theta_2+N+1$ and $N=-1$, that is
\[
\theta_1=\theta_2,\quad 
\theta_\infty=-\theta_3.
\]
We have
\begin{equation}
\label{par2}
{\small
\Omega_0=
\begin{pmatrix}
-\theta_2 & 0 & \displaystyle{\frac{\tilde{k}_1^0}{\sqrt{\theta_3}}\left(\frac{\theta_3(1+\theta_2)}{2}+(1-\theta_3)y'_0\right)}\\[2ex]
0 & -\theta_2 & \displaystyle{\frac{\tilde{k}_2^0}{\sqrt{\theta_3}}\left(\frac{\theta_3(1-\theta_2)}{2}+(1-\theta_3)y'_0\right)}\\[2ex]
\displaystyle{\frac{\sqrt{\theta_3}}{\tilde{k}_1}^0} & \displaystyle{-\frac{\sqrt{\theta_3}}{\tilde{k}_2^0}} & -\theta_3
\end{pmatrix},
}
\quad
\begin{array}{c}
\tilde{k}_1^0=k_1^0e^{i\pi\theta_2/2}
,
\\
\\
\tilde{k}_2^0=-ik_2^0e^{i\pi(\theta_2-\theta_3)}. 
\end{array}
\end{equation}
A diagonalizing matrix of $\Omega_0$ is
{\small
\[
G_0=
\begin{pmatrix}
\displaystyle{-\frac{\tilde{k}_1^0}{\theta_3\sqrt{\theta_3}}\frac{\theta_3(1+\theta_2)+2(1-\theta_3)y'_0}{2}} & \displaystyle{\frac{\tilde{k}_1^0}{\theta_2\sqrt{\theta_3}}\frac{\theta_3(1+\theta_2)+2(1-\theta_3)y'_0}{2}} & 1 \\[2ex]
\displaystyle{-\frac{\tilde{k}_2^0}{\theta_3\sqrt{\theta_3}}\frac{\theta_3(1-\theta_2)+2(1-\theta_3)y'_0}{2}} & \displaystyle{\frac{\tilde{k}_2^0}{\theta_2\sqrt{\theta_3}}\frac{\theta_3(1-\theta_2)+2(1-\theta_3)y'_0}{2}} & \displaystyle{\frac{\tilde{k}_2^0}{\tilde{k}_1^0}} \\[2ex]
1 & 1 & 0
\end{pmatrix}
.
\]
}
With the gauge $Y=G_0\widetilde{Y}$ system~\eqref{mainsystem} reduces to
\begin{equation}
\label{sysex1bis}
\frac{d\widetilde{Y}}{dz}=\left[\frac{1}{\theta_2+\theta_3}
\begin{pmatrix}
\theta_3 & \theta_3 & 0\\
\theta_2 & \theta_2 & 0\\
0 & 0 & 0\\
\end{pmatrix}
-\frac{1}{z}
\begin{pmatrix}
\theta_2+\theta_3 & 0 & 0\\
0 & 0 & 0\\
0 & 0 & \theta_2
\end{pmatrix}
\right]\widetilde{Y},
\end{equation}
which is again integrable in terms of {\bf confluent hypergeometric}  functions (see \eqref{sysex1}) with parameters $a=\theta_2$ and $b=\theta_2+\theta_3$.

If $\operatorname{Re}(\theta_1-\theta_2)>0$, condition~\eqref{coalcond} holds if and only if 
\[
\theta_1=\theta_2-N-1,
\quad  N\le-2,
\]
and the matrix $\Omega(x)$ has holomorphic limit 
{\small 
\[
\Omega_0=
\begin{pmatrix}
1-|N|-\theta_2 & 0 & -\tilde{k}_1^0\displaystyle{\sqrt{\frac{\theta_3}{|N|}}(\theta_2+|N|-1)}\\[2ex]
0 & -\theta_2 & -\tilde{k}_2^0\displaystyle{\sqrt{\frac{|N|}{\theta_3}}\widetilde{H}(\theta_2,\theta_3,N,A)} \\[2ex]
\displaystyle{-\frac{\sqrt{|N|\theta_3}}{\tilde{k}_1^0}} & 0 & -\theta_3
\end{pmatrix}
,
\]
}
where
\begin{align*}
&\widetilde{H}(\theta_2,\theta_3,N,A)=\frac{1}{2}\left.\frac{d^{|N|}}{dx^{|N|}}\left[\frac{x(x-1)dy/dx-(1-\theta_2)y(y-1)}{y-x}\right]\right|_{x=0}-(\theta_2+\theta_3+|N|-1)A,
\\
&A=\frac{y^{(|N|)}_0}{2},
\end{align*}
 and $\tilde{k}_1^0,\tilde{k}_2^0$ are as in~\eqref{par2}. A diagonalizing matrix of $\Omega_0$ is
{\small 
\[
G_0=
\begin{pmatrix}
\displaystyle{\frac{\tilde{k}_1^0}{\sqrt{|N|\theta_3}}(\theta_2+|N|-1)} & \displaystyle{-\tilde{k}_1^0\sqrt{\frac{\theta_3}{|N|}}} & 0 \\[2ex]
\displaystyle{\tilde{k}_2^0\sqrt{\frac{|N|}{\theta_3}}\frac{\widetilde{H}(\theta_2,\theta_3,N,A)}{\theta_3+|N|-1}} & -\displaystyle{\frac{\tilde{k}_2^0}{\theta_2}\sqrt{\frac{|N|}{\theta_3}}\widetilde{H}(\theta_2,\theta_3,N,A)} & 1\\[2ex]
1 & 1 & 0
\end{pmatrix}
,
\]
}
thus with the gauge $Y=G_0\widetilde{Y}$ system~\eqref{mainsystem} reduces to
\[
\frac{d\widetilde{Y}}{dz}=\left[\frac{1}{\theta_2+\theta_3+|N|-1}
\begin{pmatrix}
\theta_3 & \theta_3 & 0\\
\theta_2+|N|-1 & \theta_2+|N|-1 & 0\\
\widehat{H}(\theta_2,\theta_3,N,A) & \widehat{H}(\theta_2,\theta_3,N,A) & 0\\
\end{pmatrix}
-\frac{1}{z}
\begin{pmatrix}
\theta_2+\theta_3+|N|-1 & 0 & 0\\
0 & 0 & 0\\
0 & 0 & \theta_2
\end{pmatrix}
\right]\widetilde{Y}
,
\]
where
\[
\widehat{H}(\theta_2,\theta_3,N,A)=-\tilde{k}_2^0\widetilde{H}(\theta_2,\theta_3,N,A)\sqrt{\frac{|N|}{\theta_3}}\frac{(1-|N|)(\theta_2+\theta_3+|N|-1)}{\theta_2(\theta_3+|N|-1)}.
\]
The differential system for a column of $\widetilde{Y}$ is
{\small
\[
\begin{cases}
\displaystyle{\frac{d\tilde{y}_1}{dz}=\frac{\theta_3}{\theta_2+\theta_3+|N|-1}(\tilde{y}_1+\tilde{y}_2)-\frac{\theta_2+\theta_3+|N|-1}{z}\tilde{y}_1}\\[2ex]
\displaystyle{\frac{d\tilde{y}_2}{dz}=\frac{\theta_2+|N|-1}{\theta_2+\theta_3+|N|-1}(\tilde{y}_1+\tilde{y}_2)}\\[2ex]
\displaystyle{\frac{d\tilde{y}_3}{dz}=\frac{\widehat{H}(\theta_2,\theta_3,N,A)}{\theta_2+\theta_3+|N|-1}(\tilde{y}_1+\tilde{y}_2)-\frac{\theta_2}{z}\tilde{y}_3}.
\end{cases}
\]
}
The first two equations reduce to $
\tilde{y}_1=\frac{\theta_2+\theta_3+|N|-1}{\theta_2+|N|-1}\frac{d\tilde{y}_2}{dz}-\tilde{y}_2
$
and
\[
z\frac{d^2\tilde{y}_2}{dz^2}+(\theta_2+\theta_3+|N|-1-z)\frac{d\tilde{y}_2}{dz}-(\theta_2+|N|-1)\tilde{y}_2=0,
\]
which is a {\bf confluent hypergeometric equation} with parameters $a=\theta_2+|N|-1$ and $b=\theta_2+\theta_3+|N|-1$. The last equation is integrable by quadratures. 

If $\operatorname{Re}(\theta_1-\theta_2)<0$, condition \eqref{coalcond} holds if and only if
\[
\theta_1=\theta_2+N+1,\quad N\le-2
\]
and we have
{\small 
\[
\Omega_0=
\begin{pmatrix}
-\theta_2+|N|-1 & 0 & -\tilde{k}_1^0\displaystyle{\sqrt{\frac{|N|}{\theta_3}}\widetilde{K}(\theta_2,\theta_3,N,A)} \\[2ex]
0 & -\theta_2 & \displaystyle{\tilde{k}_2^0\theta_2\sqrt{\frac{\theta_3}{|N|}}} \\[2ex]
0 & \displaystyle{\frac{\sqrt{|N|\theta_3}}{\tilde{k}_2^0}} & -\theta_3
\end{pmatrix}
,
\]
}
where 
\[
\widetilde{K}(\theta_2,\theta_3,N,A)=\left.\frac{1}{2}\frac{d^{|N|}}{dx^{|N|}}\left[\frac{x((x-1)dy/dx+\theta_2-|N|+1)}{y}\right]\right|_{x=0}-(\theta_3+|N|)A,
\quad 
A=\frac{y_0^{(|N|)}}{2},
\] and $\tilde{k}_1^0,\tilde{k}_2^0$ as in~\eqref{par2}. A diagonalizing matrix is
{\small
\[
G_0=
\begin{pmatrix}
\displaystyle{\tilde{k}_1^0\sqrt{\frac{|N|}{\theta_3}}\frac{\widetilde{K}(\theta_2,\theta_3,N,A)}{\theta_3+|N|-1}} & \displaystyle{-\tilde{k}_1^0\sqrt{\frac{|N|}{\theta_3}}\frac{\widetilde{K}(\theta_2,\theta_3,N,A)}{\theta_2-|N|+1}} & 1\\[2ex]
\displaystyle{-\frac{\tilde{k}_2^0\theta_2}{\sqrt{|N|\theta_3}}} & \displaystyle{\tilde{k}_2^0\sqrt{\frac{\theta_3}{|N|}}} & 0\\[2ex]
1 & 1 & 0
\end{pmatrix}.
\]
}
Hence, with the gauge transformation $Y=G_0\widetilde{Y}$ system~\eqref{mainsystem} becomes
\begin{align*}
\frac{d\widetilde{Y}}{dz}=\left[\frac{1}{\theta_2+\theta_3}
\begin{pmatrix}
 \theta_3 &  \theta_3 & 0\\
\theta_2 &  \theta_2 & 0\\
\widehat{K}(\theta_2,\theta_3,N,A) & \widehat{K}(\theta_2,\theta_3,N,A) & 0
\end{pmatrix}
-\frac{1}{z}
\begin{pmatrix}
\theta_2+\theta_3 & 0 & 0\\
0 & 0 & 0\\
0 & 0 & \theta_2-|N|+1
\end{pmatrix}
\right]\widetilde{Y},
\end{align*}
where
\[
\widehat{K}(\theta_2,\theta_3,N,A)=-\tilde{k}_1^0\widetilde{K}(\theta_2,\theta_3,N,A)\sqrt{\frac{\theta_3}{|N|}}\frac{(|N|+1)(\theta_2+\theta_3)}{(\theta_2-|N|+1)(\theta_3+|N|-1)}.
\]
Let $(\tilde{y}_1,\tilde{y}_2,\tilde{y}_3)^T$ be a column of $\widetilde{Y}$, then
{\small 
\[
\begin{cases}
\displaystyle{\frac{d\tilde{y}_1}{dz}=\frac{\theta_3}{\theta_2+\theta_3}(\tilde{y}_1+\tilde{y}_2)-\frac{\theta_2+\theta_3}{z}\tilde{y}_1}\\[2ex]
\displaystyle{\frac{d\tilde{y}_2}{dz}=\frac{\theta_2}{\theta_2+\theta_3}(\tilde{y}_1+\tilde{y}_2)}\\[2ex]
\displaystyle{\frac{d\tilde{y}_3}{dz}=\frac{\widehat{K}(\theta_2,\theta_3,N,A)}{\theta_2+\theta_3}(\tilde{y}_1+\tilde{y}_2)-\frac{\theta_2-|N|+1}{z}\tilde{y}_3}.
\end{cases}
\]
}
From the second equation we get $
\tilde{y}_1=\frac{\theta_2+\theta_3}{\theta_2}\frac{d\tilde{y}_2}{dz}-\tilde{y}_2
$, 
and inserting this expression into the first equation we obtain a {\bf confluent hypergeometric equation} with parameters $a=\theta_2$ and $b=\theta_2+\theta_3$.
\[
\tilde{z}\frac{d^2\tilde{y}_2}{d\tilde{z}^2}+(\theta_2+\theta_3-\tilde{z})\frac{d\tilde{y}_2}{d\tilde{z}}-\theta_2\tilde{y}_2=0.
\]
\begin{remark}
{\rm
In this and in the previous section we have not specified from the beginning the conditions on the parameters $\theta$, such as \eqref{condizione1} or \eqref{condizione2}, but we have just required the  vanishing conditions \eqref{coalcond} 
to be satisfied.  It has turned out that $\theta_1-\theta_2\in\mathcal{N}_N$  is equivalent to the vanishing conditions, and  that for conditions \eqref{condizione1}   or $\theta_1+\theta_2\in\mathcal{N}_N$ the limit \eqref{coalcond} does not hold. In the subsequent sections we will follow another strategy, consisting in the direct substitution  into $\Omega(x)$ of the tabulated Taylor expansions,   but the same  method as above can be equally applied.
}
\end{remark}

\subsection{Transcendents with Taylor expansion $y(x)=y'_0x+O(x^2)$}
\label{secondorif}
The special cases $y'_0=0,1$ will be dealt with separately in section \ref{nongeneric}.
If $y_0'\ne 0,1$, the functions $k_1,k_2$ can be written in a neighborhood of $x=0$ as
\begin{align*}
&k_1(x)=k_1^0(1+O(x))\sqrt{-y'_0+O(x)}(x-1)^{(1+\theta_2)/2}x^{(\theta_1(y'_0-1)/y'_0-\theta_2)/2},\\
&k_2(x)=k_2^0(1+O(x))\sqrt{1-y'_0+O(x)}(x-1)^{\theta_2-\theta_3}x^{(\theta_2y'_0/(y'_0-1)-\theta_1)/2},
&
k_1^0,k_2^0\in \mathbb{C}\backslash\{0\}. 
\end{align*}
Hence, the structure of the off-diagonal elements of the matrix $\Omega(x)$ is the following:
\begin{align*}
\Omega_{12}(x)&=\omega_{12}(x)x^{(2y'_0-1)(\theta_1/y'_0-\theta_2/(y'_0-1))},&
\Omega_{13}(x)&=\omega_{13}(x)x^{(\theta_1(y'_0-1)/y'_0-\theta_2)/2},
\\
\Omega_{21}(x)&=\omega_{21}(x)x^{-(2y'_0-1)(\theta_1/y'_0-\theta_2/(y'_0-1))},&
\Omega_{23}(x)&=\omega_{23}(x)x^{(\theta_2y'_0/(y'_0-1)-\theta_1)/2},
\\
\Omega_{31}(x)&=\omega_{31}(x)x^{-(\theta_1(y'_0-1)/y'_0-\theta_2)/2},&
\Omega_{32}(x)&=\omega_{32}(x)x^{-(\theta_2y'_0/(y'_0-1)-\theta_1)/2},
\end{align*}
where $\omega_{ij}(x)$ are holomorphic functions at $x=0$. 
There are three classes of solutions. The generic one is
\begin{equation}
\label{genericorder1}
y(x)=y'_0(\theta_1,\theta_2)x+\sum_{n=2}^\infty b_n(\vec{\theta}) x^n, \quad\vec{\theta}=(\theta_1,\theta_2,\theta_3,\theta_\infty),
\end{equation} 
where
\[
y'_0=\frac{\theta_1}{\theta_1-\theta_2},\quad \theta_1-\theta_2\ne 0
;
 \quad\quad 
 \hbox{ or }
\quad\quad 
y'_0=\frac{\theta_1}{\theta_1+\theta_2},\quad \theta_1+\theta_2\ne 0.
\]
The coefficients $\{b_n\}_{n\ge 2}$ are uniquely determined. Within this class, if $\theta_1=0$ then there is only the singular solution $y(x)\equiv 0$. The other two classes are one-parameter families of solutions. One family is \begin{equation}
\label{nongenericorder11}
y(x)=y'_0x+y'_0(y'_0-1)\frac{\theta_3^2-(\theta_\infty-1)^2-1}{2}x^2+\sum_{n=3}^\infty b_n(y'_0,\theta_3,\theta_\infty)x^n,\quad \theta_1=\theta_2=0,
\end{equation}
where $y'_0$ is a free parameter.  The expression~\eqref{nongenericorder11} for $y'_0=0$ or $y'_0=1$ is the singular solution $y(x)\equiv 0$ or $y(x)\equiv 1$, respectively.  The other family is 
\begin{equation}
\label{nongenericorder12}
y(x)=\sum_{n=1}^{|N|}b_n(\vec{\theta})x^n+\frac{y_0^{(|N|+1)}}{(|N|+1)!}x^{|N|+1}+\sum_{n=|N|+2}^\infty b_n(y_0^{(|N|+1)},\vec{\theta})x^n,\quad b_1=\frac{\theta_1}{N},
\end{equation}
where $y_0^{(|N|+1)}$ is a free parameter, $N\in\mathbb{Z}\setminus\{0\}$, the relation  
$$\theta_1-\theta_2=N \quad\hbox{  or }\quad \theta_1+\theta_2=N,
$$ holds, and either 
\begin{equation}
\label{cond}
\theta_1\in \begin{cases} \{0,-1,-2,\ldots,N\}, &\mbox{if } N<0 \\ \{0,1,2,\ldots,N\}, &\mbox{if } N>0 \end{cases}
\end{equation}
or ($\mathcal{N}_N$ is the set \eqref{24giugno2021-1})
\begin{equation}
\label{condd}
\{(\theta_3+\theta_\infty-1),(-\theta_3+\theta_\infty-1)\}\cap\mathcal{N}_N\ne\emptyset.
\end{equation}

%%%%%%%%%%%%%%%%%%%%%%

\subsubsection{Generic solution~\eqref{genericorder1} - Case $y'_0=\theta_1/(\theta_1-\theta_2)$}

The matrix $\Omega(x)$ is holomorphic  at $x=0$ with $\Omega_{1}(0)=\Omega_{21}(0)=0$, and we have
{\small
\[
\Omega_0=
\begin{pmatrix}
-\theta_1 & 0 & \displaystyle{ \frac{\tilde{k}_1^0}{2}(\theta_1-\theta_2-\theta_3-\theta_\infty)}
\\[2ex]
0 & -\theta_2 & \displaystyle{-\frac{\tilde{k}_2^0}{2}(\theta_1-\theta_2+\theta_3+\theta_\infty)}
\\[2ex]
\displaystyle{-\frac{1}{2\tilde{k}_1^0}(\theta_1-\theta_2-\theta_3+\theta_\infty)\frac{\theta_1}{\theta_1-\theta_2}} 
&
\displaystyle{- \frac{1}{2\tilde{k}_2^0}(\theta_1-\theta_2+\theta_3-\theta_\infty)\frac{\theta_2}{\theta_1-\theta_2}} & -\theta_3
\end{pmatrix}
,
\]
}
where
\[
\tilde{k}_1^0=k_1^0\sqrt{\frac{\theta_1}{\theta_1-\theta_2}}e^{i\pi\theta_2/2}\quad\mbox{and}\quad\tilde{k}_2^0=ik_2^0\sqrt{\frac{\theta_2}{\theta_1-\theta_2}}e^{i\pi(\theta_2-\theta_3)}.
\]
Applying the gauge transformation $Y=G_0\widetilde{Y}$ to system~\eqref{mainsystem}, where
{\small
\[
G_0=
\begin{pmatrix}
\displaystyle{\tilde{k}_1^0\frac{\theta_1-\theta_2-\theta_3-\theta_\infty}{\theta_1-\theta_2-\theta_3+\theta_\infty}} & \displaystyle{\tilde{k}_1^0\frac{\theta_1-\theta_2-\theta_3-\theta_\infty}{2\theta_1}}  & \tilde{k}_1^0\\[2ex]
\displaystyle{\tilde{k}_2^0\frac{\theta_1-\theta_2+\theta_3+\theta_\infty}{\theta_1-\theta_2+\theta_3-\theta_\infty}} & \displaystyle{-\tilde{k}_2^0\frac{\theta_1-\theta_2+\theta_3+\theta_\infty}{2\theta_2}} & \tilde{k}_2^0\\[2ex]
1 & 1 & 1
\end{pmatrix}
,
\]
}
we can reduce the system for a column $(\tilde{y}_1,\tilde{y}_2,\tilde{y}_3)$ of $\tilde{Y}$ to the {\bf generalized hypergeometric equation}
\[
z^2\frac{d^3w}{dz^3}+z(b_2+a_2z)\frac{d^2w}{dz^2}+(b_1+a_1z)\frac{dw}{dz}+a_0w=0,
\]
where $w=\tilde{y}_1\,z^{(\theta_1+\theta_2+\theta_3-\theta_\infty)/2}$ and the parameters are
\begin{align*}
&a_0=\frac{1}{4}(\theta_1+\theta_2+\theta_3-\theta_\infty)(\theta_1+\theta_2-\theta_3+\theta_\infty)-\theta_1\theta_2,\quad a_1=\theta_3-\theta_\infty-1,\quad a_2=-1,\\[2ex]
&b_1=-\frac{\theta_\infty}{2}(\theta_1+\theta_2+\theta_3-\theta_\infty),\quad b_2=\frac{2+3\theta_\infty-\theta_1-\theta_2-\theta_3}{2}.
\end{align*}

\subsubsection{Generic solution~\eqref{genericorder1} - Case $y'_0=\theta_1/(\theta_1+\theta_2)$}

We show that for the solution \eqref{genericorder1} with $y'_0=\theta_1/(\theta_1+\theta_2)$ the condition \eqref{coalcond} does not hold.   If  $\operatorname{Re}(\theta_1-\theta_2)=0$, then the vanishing  \eqref{coalcond} 
is equivalent to $\omega_{12}(0)=\omega_{21}(0)=0$, that is $\theta_1=0$ and $\theta_2=0$,
 a contradiction. If $\operatorname{Re}(\theta_1-\theta_2)>0$, it is necessary that $\omega_{21}(0)=0$, that is $
\theta_1=0$. 
 But the only solution within this class under this condition is $y(x)\equiv 0$.
 If $\operatorname{Re}(\theta_1-\theta_2)<0$, the requirement $\lim_{x\to 0}\Omega_{12}(x)=0$ implies necessarily $\omega_{12}(0)=0$, that is 
$\theta_2=0$. 
Under this condition, we can refer to the results of the previous case. In conclusion $\Omega(x)$ is not holomorphic at $x=0$. This can also be seen by substituting \eqref{genericorder1} into the explicit formulae of $\Omega$.

\subsubsection{Solutions~\eqref{nongenericorder11}, with $\theta_1=\theta_2=0$}
\label{t3}
In this case   $\Omega(x)$ is holomorphic  at $x=0$ with $\Omega_{12}(0)=\Omega_{21}(0)=0$, and we have
\be
\label{25giugno2021-3}
{\small
\Omega_0=
\begin{pmatrix}
0 & 0 & \displaystyle{\frac{\tilde{k}_1^0}{2}(\theta_3+\theta_\infty)}\\[2ex]
0 & 0 & \displaystyle{\frac{\tilde{k}_2^0}{2}(\theta_3+\theta_\infty)}\\
\displaystyle{\frac{1}{2\tilde{k}_1^0}(\theta_\infty-\theta_3)y'_0} & \displaystyle{\frac{1}{2\tilde{k}_2^0}(\theta_\infty-\theta_3)(1-y'_0)} & -\theta_3
\end{pmatrix}
},
\quad
\quad
\begin{array}{c}
\tilde{k}_1^0=k_1^0\sqrt{y'_0},
\\
\\
\tilde{k}_2^0=-k_2^0\sqrt{1-y'_0}e^{-i\pi\theta_3}.
\end{array}
\ee
Let $Y=G_0\widetilde{Y}$, where
{\small 
\begin{equation}
\label{g0t3}
G_0=
\begin{pmatrix}
-\tilde{k}_1^0\displaystyle{\frac{\theta_3+\theta_\infty}{\theta_3-\theta_\infty}} & -\displaystyle{\frac{\tilde{k}_1^0}{\tilde{k}_2^0}\frac{1-y'_0}{y'_0}} & -\tilde{k}_1^0\\[2ex]
-\displaystyle{\tilde{k}_2^0\frac{\theta_3+\theta_\infty}{\theta_3-\theta_\infty}} & 1 & -\tilde{k}_2^0\\[2ex]
1 & 0 & 1
\end{pmatrix}
,
\end{equation}
}
then
 system~\eqref{mainsystem}
  is
\begin{equation}
\label{syst3}
\frac{d\widetilde{Y}}{dz}=\left[\frac{1}{2\theta_\infty}
\begin{pmatrix}
-\theta_3+\theta_\infty & 0 & -\theta_3+\theta_\infty\\
0 & 0 & 0\\
\theta_3+\theta_\infty & 0 & \theta_3+\theta_\infty
\end{pmatrix}
+\frac{1}{2z}
\begin{pmatrix}
-\theta_3+\theta_\infty & 0 & 0\\
0 & 0 & 0\\
0 & 0 & -\theta_3-\theta_\infty
\end{pmatrix}
\right]\widetilde{Y}.
\end{equation}
If $(\tilde{y}_1,\tilde{y}_2,\tilde{y}_3)^T$ is a column of $\widetilde{Y}$, then
{\small 
\[
\begin{cases}
\displaystyle{\frac{d\tilde{y}_1}{dz}=\frac{\theta_\infty-\theta_3}{2\theta_\infty}(\tilde{y}_1+\tilde{y}_3)+\frac{\theta_\infty-\theta_3}{2z}\tilde{y}_1},\\[2ex]
\displaystyle{\frac{d\tilde{y}_2}{dz}=0},\\[2ex]
\displaystyle{\frac{d\tilde{y}_3}{dz}=\frac{\theta_\infty+\theta_3}{2\theta_\infty}(\tilde{y}_1+\tilde{y}_3)-\frac{\theta_\infty+\theta_3}{2z}\tilde{y}_3}.
\end{cases}
\]
}
Substituting the first equation into the third to eliminate $\tilde{y}_3$  and setting  $\tilde{y}_1=w\,z^{(\theta_\infty-\theta_3)/2}$ we receive 
\begin{equation}
\label{y1t3}
z\frac{d^2w}{dz^2}+(\theta_\infty-z)\frac{dw}{dz}-\frac{\theta_\infty-\theta_3}{2}w=0,
\end{equation}
which is a {\bf confluent hypergeometric equation} with parameters $a=(\theta_\infty-\theta_3)/2$ and $b=\theta_\infty$.

\subsubsection{Solution~\eqref{nongenericorder12} -  Case $\theta_1-\theta_2=N$ with condition~\eqref{cond}}\label{nongeneric}
Let $\theta_1=k$ integer, in the set \eqref{cond}, so that  $\theta_2=k-N$.

If $k\ne 0,N$, then the matrix $\Omega(x)$ is holomorphic at $x=0$, $\Omega_{12}(0)=\Omega_{21}(0)=0$, and 
{\small \[
\Omega_0=
\begin{pmatrix}
-k & 0 & \displaystyle{\frac{\tilde{k}_1^0}{2}(\theta_3+\theta_\infty-N)}\\[2ex]
0 & N-k & \displaystyle{\frac{\tilde{k}_2^0}{2}(\theta_3+\theta_\infty+N)}\\[2ex]
\displaystyle{\frac{1}{2\tilde{k}_1^0}\frac{k}{N}(\theta_\infty-\theta_3+N)} & \displaystyle{\frac{1}{2\tilde{k}_2^0}\frac{N-k}{N}(\theta_\infty-\theta_3-N)} & -\theta_3
\end{pmatrix}
,
\quad
\quad
\begin{array}{c}
\tilde{k}_1^0=k_1^0\sqrt{\frac{k}{N}}e^{i\pi\theta_2/2},
\\
\\
\tilde{k}_2^0=-k_2^0\sqrt{\frac{N-k}{N}}e^{i\pi(\theta_2-\theta_3)}.
\end{array}
\]
}
System~\eqref{mainsystem} reduces to the {\bf generalized hypergeometric equation}
\[
z^2\frac{d^3w}{dz}+z(a_2+b_2z)\frac{d^2w}{dz^2}+(a_1+b_1z)\frac{dw}{dz}+a_0w=0
\]
with parameters
\begin{align*}
&a_0=\frac{1}{2}(N^2+2k(-2+k-N+\theta_3-\theta_\infty)+N(2-\theta_3+\theta_\infty)),\quad a_1=\frac{1}{2}(-4-N+2k+3\theta_3-3\theta_\infty),\quad a_2=-1,\\[2ex]
&b_1=\frac{1}{2}(2+\theta_\infty)(2+N-2k-\theta_3+\theta_\infty),\quad b_2=\frac{1}{2}(8+N-2k-\theta_3+3\theta_\infty).
\end{align*}

\vskip 0.2 cm 
If $k=0$, that is $\theta_1=0$ and $\theta_2=-N$, then
\[
y(x)=Ax^{|N|+1}+O(x^{|N|+2}),\quad\quad A:=\frac{y_0^{(|N|+1)}}{(|N|+1)!} \hbox{ free parameter},
\]
and
\[
k_1(x)=k_1^0(1+O(x))x^{(N+|N|)/2}(x-1)^{(1-N)/2},
\quad
k_2(x)=k_2^0(1+O(x))(x-1)^{-N-\theta_3},\quad k_1^0,k_2^0\in\mathbb{C}.
\]
 We need to distinguish the cases with negative and positive $N$.  First, let $N\le -1$, then $\Omega(x)$ is holomorphic at $x=0$ and 
\be
\label{7luglio2021-1}
{\small \Omega_0=
\begin{pmatrix}
0 & 0 & -\frac{\tilde{k}_1^0}{2}(\theta_3+\theta_\infty+|N|)\\[2ex]
0 & -|N| & -\displaystyle{\frac{\tilde{k}_2^0}{2}(\theta_3+\theta_\infty-|N|)}\\[2ex]
0 & \displaystyle{\frac{1}{2\tilde{k}_2^0}(\theta_3-\theta_\infty-|N|)} & -\theta_3
\end{pmatrix}
}
,
\quad\quad
\tilde{k}_1^0=ik_1^0e^{i\pi|N|/2}\quad\tilde{k}_2^0=k_2^0e^{i\pi(|N|-\theta_3)}.
\ee
Applying the gauge transformation $Y=G_0\widetilde{Y}$, where
\[
G_0=
\begin{pmatrix}
\displaystyle{\tilde{k}_1^0\frac{\theta_3+\theta_\infty+|N|}{\theta_3-\theta_\infty+|N|}} & 1 & \tilde{k}_1^0 \\[2ex]
\displaystyle{\tilde{k}_2^0\frac{\theta_3+\theta_\infty-|N|}{\theta_3-\theta_\infty-|N|}} & 0 & \tilde{k}_2^0\\[2ex]
1 & 0 & 1 
\end{pmatrix}
\]
is a diagonalizing matrix of $\Omega_0$, system~\eqref{mainsystem} becomes
\[
\frac{d\widetilde{Y}}{dz}=\left[
\begin{pmatrix}
\displaystyle{-\frac{\theta_3-\theta_\infty-|N|}{2\theta_\infty}} & 0 & \displaystyle{-\frac{\theta_3-\theta_\infty-|N|}{2\theta_\infty}}\\[2ex]
\displaystyle{-\frac{2|N|\tilde{k}_1^0}{\theta_3-\theta_\infty+|N|}} & 0 & \displaystyle{-\frac{2|N|\tilde{k}_1^0}{\theta_3-\theta_\infty+|N|}}\\[2ex]
\displaystyle{\frac{\theta_3+\theta_\infty-|N|}{2\theta_\infty}} & 0 & \displaystyle{\frac{\theta_3+\theta_\infty-|N|}{2\theta_\infty}}
\end{pmatrix}
-\frac{1}{2z}
\begin{pmatrix}
\theta_3-\theta_\infty+|N| & 0 & 0\\[2ex]
0 & 0 & 0\\[2ex]
0 & 0 & \theta_3+\theta_\infty+|N|
\end{pmatrix}
\right]\widetilde{Y}.
\]
If $(\tilde{y}_1,\tilde{y}_2,\tilde{y}_3)^T$ is a column of $\widetilde{Y}$, then
{\small 
\[
\begin{cases}
\displaystyle{\frac{d\tilde{y}_1}{dz}=-\frac{\theta_3-\theta_\infty-|N|}{2\theta_\infty}(\tilde{y}_1+\tilde{y}_3)-\frac{1}{2z}(\theta_3-\theta_\infty+|N|)\tilde{y}_1}\\[2ex]
\displaystyle{\frac{d\tilde{y}_2}{dz}=-\frac{2|N|\tilde{k}_1}{\theta_3-\theta_\infty+|N|}(\tilde{y}_1+\tilde{y}_3)}\\[2ex]
\displaystyle{\frac{d\tilde{y}_3}{dz}=\frac{\theta_3+\theta_\infty-|N|}{2\theta_\infty}(\tilde{y}_1+\tilde{y}_3)-\frac{1}{2z}(\theta_3+\theta_\infty+|N|)}\tilde{y}_3.
\end{cases}
\]
}
From the first equation we get
\[
\tilde{y}_3=-\frac{2\theta_\infty}{\theta_3-\theta_\infty-|N|}\left(\frac{d\tilde{y}_1}{dz}+\frac{\theta_3-\theta_\infty+|N|}{2z}\tilde{y}_1\right)-\tilde{y}_1,
\]
Substitution  the first equation  into the third to eliminate $\tilde{y}_3$  and the change $\tilde{y}_1=w\,z^{(\theta_\infty-\theta_3-|N|)/2}$   give 
\[
z\frac{d^2w}{dz^2}+(\theta_\infty-z)\frac{dw}{dz}-\frac{\theta_\infty-\theta_3+|N|}{2}w=0,
\] 
which is a {\bf confluent hypergeometric equation} with parameters $a=(\theta_\infty-\theta_3+|N|)/2$ and $b=\theta_\infty$.

If $N\ge 1$, then
\be
\label{24giugno2021-2}
{\small \Omega_0=
\begin{pmatrix}
0 & 0 & 0\\[2ex]
0 & |N| & \displaystyle{-\frac{\tilde{k}_2^0(\theta_3+\theta_\infty+|N|)}{2}}\\[2ex]
\displaystyle{\frac{A(\theta_3-\theta_\infty-|N|)}{2\tilde{k}_1^0}} & \displaystyle{\frac{\theta_3-\theta_\infty+|N|}{2\tilde{k}_2^0}} & -\theta_3
\end{pmatrix}
}
,
\quad\quad
\begin{array}{c}
\tilde{k}_1^0=ik_1^0e^{-i\pi|N|/2},
\\
\\
\tilde{k}_2^0=k_2^0e^{i\pi(|N|+\theta_3)}.
\end{array}
\ee
Let $Y=G_0\widetilde{Y}$, with
{\small \[
G_0=
\begin{pmatrix}
0 & \displaystyle{-\frac{\tilde{k}_1^0}{2}\frac{\theta_3^2-\theta_\infty^2+|N|^2}{A|N|(\theta_3-\theta_\infty-|N|)}} & 0\\[2ex]
\displaystyle{\tilde{k}_2^0\frac{\theta_3+\theta_\infty+|N|}{\theta_3-\theta_\infty+|N|}} & \displaystyle{\frac{\tilde{k}_2^0}{2}\frac{\theta_3+\theta_\infty+|N|}{|N|}} & \tilde{k}_2^0\\[2ex]
1 & 1 & 1
\end{pmatrix}
.
\]
}
System~\eqref{mainsystem} becomes
{\small
$$
\frac{d\widetilde{Y}}{dz}=\left[\frac{1}{2\theta_\infty}
\begin{pmatrix}
-(\theta_3-\theta_\infty+|N|) & -(\theta_3-\theta_\infty+|N|) & -(\theta_3-\theta_\infty+|N|)\\[2ex]
0 & 0 & 0\\[2ex]
\theta_3+\theta_\infty+|N| & \theta_3+\theta_\infty+|N| & \theta_3+\theta_\infty+|N|
\end{pmatrix}
-\frac{1}{2z}
\begin{pmatrix}
\theta_3-\theta_\infty-|N| & 0 & 0\\[2ex]
0 & 0 & 0\\[2ex]
0 & 0 & \theta_3+\theta_\infty-|N|
\end{pmatrix}
\right]\widetilde{Y},
$$
}
which is the following system for a column $(\tilde{y}_1,\tilde{y}_2,\tilde{y}_3)^T$:
{\small \[
\begin{cases}
\displaystyle{\frac{d\tilde{y}_1}{dz}=-\frac{\theta_3-\theta_\infty+|N|}{2\theta_\infty}(\tilde{y}_1+\tilde{y_2}+\tilde{y}_3)-\frac{\theta_3-\theta_\infty-|N|}{2z}\tilde{y}_1}\\[2ex]
\displaystyle{\frac{d\tilde{y}_2}{dz}=0}\\[2ex]
\displaystyle{\frac{d\tilde{y}_3}{dz}=\frac{\theta_3+\theta_\infty+|N|}{2\theta_\infty}(\tilde{y}_1+\tilde{y}_2+\tilde{y}_3)-\frac{\theta_3+\theta_\infty-|N|}{2z}\tilde{y}_3}.
\end{cases}
\]
}
The second equation implies that $\tilde{y}_2$ is a constant $C\in\mathbb{C}$. Computing $\tilde{y}_3$ from the first and substituting into the third we receive
 the inhomogeneous {\bf confluent hypergeometric equation}
\[
z\frac{d^2w}{dz^2}+(\theta_\infty-z)\frac{dw}{dz}-\frac{\theta_\infty-\theta_3-|N|}{2}w+C\frac{(\theta_3+\theta_\infty-|N|)(\theta_3-\theta_\infty+|N|)}{2}=0
\]
with parameters $a=(\theta_\infty-\theta_3-|N|)/2$ and $b=\theta_\infty$, where $\tilde{y}_1=w\,z^{(\theta_\infty-\theta_3+|N|)/2}$.

If $k=N$, then
\[
y(x)=x+Ax^{|N|+1}+\sum_{n=|N|+2}^\infty b_n(A,\theta_3,\theta_\infty)x^n,
\]
where $A$ is the free parameter of the family of solutions. In this case $
k_1(x)=k_1^0(1+O(x))\sqrt{x-1}$ and $
k_2(x)=k_2^0(1+O(x))(x-1)^{-\theta_3}x^{(|N|-N)/2}$. 
We have to distinguish between the case with positive $N$ and the case with negative $N$.

If $N\ge 1$, then 
\begin{equation}
\label{par3}
{\small \Omega_0=
\begin{pmatrix}
-|N| & 0 & \displaystyle{-\frac{\tilde{k}_1^0(\theta_3+\theta_\infty-|N|)}{2}}\\[2ex]
0 & 0 & \displaystyle{-\frac{\tilde{k}_2^0(\theta_3+\theta_\infty+|N|)}{2}}\\[2ex]
\displaystyle{\frac{\theta_3-\theta_\infty-|N|}{2\tilde{k}_1^0}} & 0 & -\theta_3
\end{pmatrix}
}
,
\quad
\quad
\tilde{k}_1^0=ik_1^0,\quad\tilde{k}_2^0=k_2^0e^{-i\pi\theta_3}.
\end{equation}
We do gauge transformation  $\widehat{Y}=PY$ by  the permutation matrix $P$ below 
\[
P=
\begin{pmatrix}
0 & 1 & 0\\
1 & 0 & 0\\
0 & 0 & 1
\end{pmatrix}
\quad\Longrightarrow \quad P^{-1}U_0P=U_0,\quad 
P^{-1}\Omega_0P={\small 
\begin{pmatrix}
0 & 0 & \displaystyle{-\frac{\tilde{k}_2^0(\theta_3+\theta_\infty+|N|)}{2}}
\\[2ex]
0 & -|N| & \displaystyle{-\frac{\tilde{k}_1^0(\theta_3+\theta_\infty-|N|)}{2}}
\\[2ex]
0 & \displaystyle{\frac{\theta_3-\theta_\infty-|N|}{2\tilde{k}_1^0}} & -\theta_3
\end{pmatrix}}.
\]
This   transforms system~\eqref{mainsystem} with $\Omega_0$ in \eqref{par3} into a system with the same structure as that for the previous case  $k=0,N\le-1$, with the replacements $\tilde{k}_1^0\to \tilde{k}_2^0$ and $\tilde{k}_2^0\to \tilde{k}_1^0$ (see \eqref{7luglio2021-1}). We conclude that also on this case the computation of the fundamental matrix solution is reduced to the {\bf confluent hypergeometric equation} with parameters $a=(\theta_\infty-\theta_3+|N|)/2$ and $b=\theta_\infty$.

If $N\le-1$, then
{\small \[
\Omega_0=
\begin{pmatrix}
|N| & 0 & \displaystyle{-\frac{\tilde{k}_1^0}{2}(\theta_3+\theta_\infty+|N|)}
\\[2ex]
0 & 0 & 0
\\[2ex]
\displaystyle{\frac{1}{2\tilde{k}_1^0}(\theta_3-\theta_\infty+|N|)} & \displaystyle{-\frac{A}{2\tilde{k}_2^0}(\theta_3-\theta_\infty-|N|)}  & -\theta_3
\end{pmatrix}
,
\]
}
where $\tilde{k}_1^0,\tilde{k}_2^0$ are as in~\eqref{par3}. Again, if $P$ is the same permutation matrix as before, the gauge $\widehat{Y}=PY$ leads to a system with the same form \eqref{24giugno2021-2} of the previous case with $k=0,N\ge 1$, exchanging names of  $\tilde{k}_1^0$ and $\tilde{k}_2^0$. The computation of the fundamental matrix soution is then reduced to the same inhomogeneous  {\bf  confluent hypergeometric equation}.

\subsubsection{Solutions~\eqref{nongenericorder12} - Case $\theta_1-\theta_2=N$ with condition~\eqref{condd} }
The parameters $\theta_3,\theta_\infty$ satisfy condition \eqref{condd}: $\{(\theta_3+\theta_\infty-1),(-\theta_3+\theta_\infty-1)\}\cap\mathcal{N}_N\ne\emptyset$.  We have $\Omega_{ij}(x)=\omega_{ij}(x)$, $i,j=1,2,3$,$i\ne j$, thus the matrix $\Omega(x)$ is holomorphic at $x=0$, $\omega_{12}(0)=\omega_{21}(0)=0$ 
and
{\small
\[
 \Omega_0=
\begin{pmatrix}
-(N+\theta_2) & 0 & \displaystyle{\frac{\tilde{k}_1^0}{2}(\theta_3+\theta_\infty-N)}\\
0 & -\theta_2 & \displaystyle{-\frac{\tilde{k}_2^0}{2}(\theta_3+\theta_\infty+N)}\\
\displaystyle{\frac{1}{2\tilde{k}_1^0}\frac{N+\theta_2}{N}(\theta_\infty-\theta_3+N)} & \displaystyle{\frac{1}{2\tilde{k}_2^0}\frac{\theta_2}{N}(\theta_\infty-\theta_3-N)} & -\theta_3
\end{pmatrix}
,\quad
\quad 
\begin{array}{c}
\tilde{k}_1^0=k_1^0\sqrt{\frac{N+\theta_2}{N}}e^{i\pi\theta_2/2},
\\
\\
\tilde{k}_2^0=ik_2^0\sqrt{\frac{\theta_2}{N}}e^{i\pi(\theta_2-\theta_3)}.
\end{array}
\]
}
A diagonalizing matrix of $\Omega_0$ is
{\small \[
G_0=
\begin{pmatrix}
\displaystyle{-\tilde{k}_1^0\frac{\theta_3+\theta_\infty-N}{\theta_3-\theta_\infty-N}} & \displaystyle{\frac{\tilde{k}_1^0}{2}\frac{\theta_3+\theta_\infty-N}{\theta_2+N}} & -\tilde{k}_1^0 \\
\displaystyle{\tilde{k}_2^0\frac{\theta_3+\theta_\infty+N}{\theta_3-\theta_\infty+N}} & \displaystyle{-\frac{\tilde{k}_2^0}{2}\frac{\theta_3+\theta_\infty+N}{\theta_2}} & \tilde{k}_2^0 \\
1 & 1 & 1
\end{pmatrix}
.
\]
}
With the gauge $Y=G_0\widetilde{Y}$, our system becomes
\[
\frac{d\widetilde{Y}}{dz}=\left[\frac{1}{p}
\begin{pmatrix}
-q & -q & -q\\
4\theta_2(\theta_2+N) & 4\theta_2(\theta_2+N) & 4\theta_2(\theta_2+N)\\
p-4\theta_2(\theta_2+N)+q & p-4\theta_2(\theta_2+N)+q & p-4\theta_2(\theta_2+N)+q
\end{pmatrix}
-\frac{1}{z}
\begin{pmatrix}
l & 0  & 0\\
0 & 0 & 0 \\
0 & 0 & l+\theta_\infty
\end{pmatrix}
\right]\widetilde{Y},
\]
where
\begin{align*}
p&=2\theta_3(2\theta_2+\theta_3+N)+4\theta_2(\theta_2+N)-\theta_3^2-\theta_\infty^2+N^2,& l=\theta_2+\frac{\theta_3-\theta_\infty+N}{2},
\\
q&=\frac{2\theta_2+\theta_3+\theta_\infty+N}{2\theta_\infty}(\theta_3-\theta_\infty-N)(\theta_3-\theta_\infty+N),
\end{align*}
If $(\tilde{y}_1,\tilde{y}_2,\tilde{y}_3)^T$ is a column of $\widetilde{Y}$, then
{\small \[
\begin{cases}
\displaystyle{\frac{d\tilde{y}_1}{dz}=-\frac{q}{p}(\tilde{y}_1+\tilde{y}_2+\tilde{y}_3)-\frac{l}{z}\tilde{y}_1}\\[2ex]
\displaystyle{\frac{d\tilde{y}_2}{dz}=\frac{4\theta_2(\theta_2+N)}{p}(\tilde{y}_1+\tilde{y}_2+\tilde{y}_3)}\\[2ex]
\displaystyle{\frac{d\tilde{y}_3}{dz}=\frac{p-4\theta_2(\theta_2+N)+q}{p}(\tilde{y}_1+\tilde{y}_2+\tilde{y}_3)-\frac{l+\theta_\infty}{z}\tilde{y}_3},
\end{cases}
\]
}
hence
{\small \begin{align*}
\tilde{y}_1=\frac{p}{4\theta_2(\theta_2+N)}\frac{d\tilde{y}_2}{dz}-\tilde{y}_2-\tilde{y}_3,&&
\tilde{y}_2=\frac{1}{l(p-4\theta_2(\theta_2+N)+q)}\left\{pz\frac{d^2\tilde{y}_3}{dz^2}+\left[(2l+\theta_\infty)p-pz\right]\frac{d\tilde{y}_3}{dz}+\right.
\\
&&+\left.\left[-4\theta_2\theta_\infty(\theta_2+N)+\theta_\infty-lp+\frac{p(l+\theta_\infty)(l-1)}{z}\right]\tilde{y}_3\right\},
\end{align*}
}
and
{\small \begin{align*}
z^2\frac{d^3\tilde{y}_3}{dz^3}+z(2l+\theta_\infty+1-z)\frac{d^2\tilde{y}_3}{dz^2}+\left[(l+\theta_\infty)(l-1)+\frac{-4\theta_2(\theta_2+N)(\theta_\infty+l)+\theta_\infty q-p(l+1)}{p}z\right]\frac{d\tilde{y}_3}{dz}\\
-\left[\frac{4l\theta_2(\theta_2+N)(l+\theta_\infty)}{p}+\frac{(l+\theta_\infty)(l-1)}{z}\right]\tilde{y}_3=0.
\end{align*}
}
Defining the function $w$ by $\tilde{y}_3=w\,z^{1-l}$, the last equation becomes
\begin{equation}
\begin{aligned}
z^2\frac{d^3w}{dz^3}+z(4-l+\theta_\infty-z)\frac{d^2w}{dz^2}+\left[(1-l)(2+\theta_\infty)+\frac{p(l-3)-4\theta_2(\theta_2+N)(\theta_\infty+l)+\theta_\infty q}{p}z\right]\frac{dw}{dz}\\
+\frac{1}{p}(-4\theta_2(\theta_2+N)(\theta_\infty+l)+(\theta_\infty q-pl)(1-l))w=0,
\end{aligned}
\end{equation}
which is a {\bf generalized hypergeometric equation}, whose parameters in terms of $\theta_2,\theta_3,\theta_\infty,N$ are
\begin{align*}
&a_0=\frac{1}{2}(N^2+2\theta_2(-2+\theta_2+\theta_3-\theta_\infty)+N(-2+2\theta_2+\theta_3-\theta_\infty)),\quad a_1=\frac{1}{2}(-4+N+2\theta_2+3\theta_3-3\theta_\infty),\\[2ex]
&a_2=-1,\quad b_1=\frac{1}{2}(2+\theta_\infty)(2-N-2\theta_2-\theta_3+\theta_\infty),\quad b_2=\frac{1}{2}(8-N-2\theta_2-\theta_3+3\theta_\infty).
\end{align*}

\subsubsection{Solution~\eqref{nongenericorder12} - Case $\theta_1+\theta_2=N$}
Let us start assuming condition~\eqref{condd} for $\theta_3,\theta_\infty$ holds. For this family of solutions, \eqref{coalcond} does not hold. Indeed, the entries $(1,2)$ and $(2,1)$ of $\Omega(x)$ are
\[
\Omega_{12}(x)=\omega_{12}(x)\,x^{2(N-2\theta_2)}
\quad\hbox{ and }
\quad
\Omega_{21}(x)=\omega_{21}(x)\, x^{-2(N-2\theta_2)}.
\]
If $\operatorname{Re}(\theta_2)=N/2$, $N\in\mathbb{Z}\backslash\{0\}$,  then the conditions $\omega_{12}(0)=\omega_{21}(0)=0$ imply $N=0$ and $\theta_2=0$, a contraddiction. If $\operatorname{Re}(\theta_2)<N/2$, then condition $\omega_{21}(0)=0$ is equivalent to $\theta_2=N$, an absurd. If $\operatorname{Re}(\theta_2)>N/2$, then $\omega_{12}(0)$ imply $\theta_2=0$, which is again a contraddiction.

If condition~\eqref{cond} hold, by the same arguments,  $\Omega(x)$ is not holomorphic at $x=0$ for $k\ne 0,N$, while for $k=0$ or $k=N$ the problem can be traced back to the previous section.

%%%%%%%%%%%%%%%%%%%%%%%

   %%%%%%%%%%%%%%%%%%%%%

\section{Monodromy data -  Examples and applications}
\label{20giugno2021-1}

  For a selection of transcendents admitting Taylor expansion at $x=0$ such that  $\Omega(x)$ is holomorphic at $x=0$ with vanishing conditions \eqref{coalcond}, we compute the Stokes matrices of system \eqref{systempainleve}, and the corresponding  monodromy invariants $p_{jk}$ of the $2\times 2$ Fuchsian system \eqref{23gennaio2021-8}, using  formulae \eqref{3febbraio2021-2}.  

The Stokes matrices are defined  as in  \eqref{27gennaio2021-9}. Thanks to theorem 1.1 in \cite{CDG}, it suffices to consider  the simplified system \eqref{mainsystem}  at $x=0$ and its fundamental solutions $Y_1(z),Y_2(z),Y_3(z)$ with canonical asymptotics 
\be
\label{7luglio2021-2}
Y_j(z)\sim Y_F(z,0), \quad  z\to\infty \hbox{ in } \mathcal{S}_j,\quad j=1,2,3,
\ee
where $Y_F(z,0)$ is the value at $x=0$ of the unique formal solution 
$$Y_F(z,x)=  (I+\sum_{k\geq 1} F_k(x)z^{-k}) z^{-\Theta}e^{zU(x)},
$$ 
of system \eqref{systempainleve}, whose coefficients $F_k(x)$ are also holomorphic at $x=0$. Then, the Stokes matrices of \eqref{systempainleve} are just obtained from the relations
$$
  Y_2(z)=Y_1(z)\mathbb{S}_1,
  \quad\quad Y_3(z)=Y_2(z)\mathbb{S}_{2}.
 $$
 We will choose the basic Stokes sectors in the universal covering of $\mathbb{C}\backslash\{0\}$ to be 
\[
\mathcal{S}_1:\,-\frac{3\pi}{2}<\arg(z)<\frac{\pi}{2},
\quad\quad
\mathcal{S}_2:\,-\frac{\pi}{2}<\arg(z)<\frac{3\pi}{2},
\quad
\quad
\mathcal{S}_1:\,\frac{\pi}{2}<\arg(z)<\frac{5\pi}{2},
\] 
so that $\mathcal{S}_1$ contains the admissible ray in direction $\tau=0$, corresponding to $\eta=3\pi/2$ in the $\lambda$-plane. 
 Compared to \eqref{7luglio2021-2}, we will  do the gauge transformation  $Y_j=G_0 \widetilde{Y}_j $, which does not affect the Stokes matrices,  so that the canonical asymptotics will be 
$$
\widetilde{Y}_j(z)\sim \widetilde{Y}_F(z):=G_0^{-1}Y_F(z,0),\quad z\to\infty \hbox{ in } \mathcal{S}_j,\quad j=1,2,3.
$$

\bre
\label{28luglio2021-1}
{\rm  In the examples we study below, no partial resonance occurs, so that it is not necessary to  compute $Y_F(z,x)$ and evaluate it at $x=0$, being sufficient  to compute directly the formal solution of \eqref{mainsystem},  with behaviour 
$Y_F(z)=  (I+\sum_{k\geq 1} \mathring{F}_kz^{-k}) z^{-\Theta}e^{zU(0)}
$,  
which is unique (see Remark \ref{7luglio2021-3}) and automatically coincides with $Y_F(z,0)$, that is  $\mathring{F}_k=F_k(0)$.
    }
    \ere
%%%%%%%%%%%%%%%%%%%%%%
\subsection{Case (T1) of the table with $N=-1$, $\theta_1=\theta_2$, $\theta_3=\theta_\infty$}

For $ 
\theta_2=\theta_1$, $\theta_3=\theta_\infty$,   
we consider the transcendent  in case (T1) of the table
\be
\label{2maggio2021-10}
y(z)= \frac{1}{1-\theta_\infty} + y_0^\prime x + O(x^2), \quad x\to 0,
\ee 
 holomorhic at $x=0$. This is the solution appearing in formula (16) of \cite{guz2006}, whose monodromy data  $\{p_{jk},~1\leq j\neq k \leq 3\}$ for the $2\times 2$ Fuchsian system  are given in theorem 3 of \cite{guz2006}, for $\theta_1,\theta_\infty$ not integers.  It also coincides with   (63) of  \cite{guz2012}, with monodromy data  given at page 3258, especially in the footnote 4. 
Here, we compute the Stokes matrices  and apply the formulae of  Theorem \ref{16aprile2021-5}, so  re-obtain and verifying  the $p_{jk}$ of \cite{guz2006,guz2012}.

Since $\theta_1-\theta_2\not \in\mathbb{Z}\backslash\{0\}$, there is no partial resonance, so that the formal solution of~\eqref{sysex1}  is unique and   computed following \cite{CDG}, Section 4:

\begin{align*}
\widetilde{Y}_F(z)=&\left(1+O\left(\frac{1}{z}\right)\right)\cdot 
\\
&
\cdot \begin{pmatrix}
\displaystyle{\frac{\tilde{k}_2^0}{\tilde{k}_1^0}\frac{2y'_0(\theta_\infty-1)+\theta_\infty(\theta_1-1)}{2\theta_1\theta_\infty}z^{-\theta_1}} &  \displaystyle{\frac{2y'_0(\theta_\infty-1)+\theta_\infty(\theta_1-1)}{2\theta_1\theta_\infty}z^{-\theta_1}} & 0\\[2ex]
\displaystyle{-\frac{2y'_0(\theta_\infty-1)+\theta_\infty(\theta_1-1)}{2\tilde{k}_1^0\sqrt{\theta_\infty}(\theta_1+\theta_\infty)}z^{-\theta_1}} & \displaystyle{-\frac{2y'_0(\theta_\infty-1)-\theta_\infty(\theta_1+1)}{2\tilde{k}_2^0\sqrt{\theta_\infty}(\theta_1+\theta_\infty)}z^{-\theta_1}} & \displaystyle{\frac{\theta_1}{\theta_1+\theta_\infty} z^{-\theta_\infty}e^z}\\[2ex]
\displaystyle{\frac{2y'_0(\theta_\infty-1)+\theta_\infty(\theta_1-1)}{2\tilde{k}_1^0\sqrt{\theta_\infty}(\theta_1+\theta_\infty)}z^{-\theta_1}} &  \displaystyle{\frac{2y'_0(\theta_\infty-1)-\theta_\infty(\theta_1+1)}{2\tilde{k}_2^0\sqrt{\theta_\infty}(\theta_1+\theta_\infty)} z^{-\theta_1}} & \displaystyle{\frac{\theta_\infty}{\theta_1+\theta_\infty} z^{-\theta_\infty}e^z}
\end{pmatrix}
,
\end{align*}
where $\tilde{k}_1^0,\tilde{k}_2^0$ are given in~\eqref{par1}.

As in \cite{guz2006}, we consider  the case $\theta_1,\theta_\infty$ not integers.  For the fundamental matrix solutions of~\eqref{sysex1}, the general form of the elements of the first row is $\tilde{y}_1(z)=Cz^{-\theta_2}$ form \eqref{y1ex1}, while from~\eqref{y2ex1} the general form of the elements of the second row is
\[
\tilde{y}_2(z;a,b,m,n)=aM(ze^{2\pi i m};\theta_1,\theta_1+\theta_\infty)+bU(ze^{2\pi i n};\theta_1,\theta_1+\theta_\infty),\quad a,b\in\mathbb{C},n,m\in\mathbb{Z}.
\]
For the sake of the computation of the Stokes matrices, it is sufficient to compute the first two rows of the fundamental solutions $\widetilde{Y}_1,\widetilde{Y}_2,\widetilde{Y}_3$ with asymptotics $\widetilde{Y}_F$ in the Stokes sectors $\mathcal{S}_1,\mathcal{S}_2,\mathcal{S}_3$, respectively. The computation of the first row is immediate, by comparision with the first row of the leading term of $\widetilde{Y}_F$. The second row is obtained by comparison of the second row of the leading term of $\widetilde{Y}_F$ and the leading coefficients of the asymptotics of the confluent hypergeometric functions $M$ and $U$. For the fundamental solution $\widetilde{Y}_1$ we use formulas~\eqref{asymptoticsU} and~\eqref{asymptoticsM} with $\epsilon=-1$:
\[
\widetilde{Y}_1(z)=
\begin{pmatrix}
c_1\,z^{-\theta_2} & c_2\,z^{-\theta_2} & 0\\[1ex]
\tilde{y}_2(z;0,b_1,0,0) & \tilde{y}_2(z;0,b_2,0,0) & \tilde{y}_2(z;a,b_3,0,0)\\[1ex]
(\widetilde{Y}_1(z))_{31} & (\widetilde{Y}_1(z))_{32} & (\widetilde{Y}_1(z))_{33}
\end{pmatrix}
,
\]
where
\begin{align*}
&a= \frac{\theta_1}{\theta_1+\theta_\infty}\Gamma(\theta_1),\\[2ex]
&b_1=-\frac{2y'_0(\theta_\infty-1)+\theta_\infty(\theta_1-1)}{2\tilde{k}_1^0\sqrt{\theta_\infty}(\theta_1+\theta_\infty)},\quad b_2= -\frac{2y'_0(\theta_\infty-1)-\theta_\infty(\theta_1+1)}{2\tilde{k}_2^0\sqrt{\theta_\infty}(\theta_1+\theta_\infty)},\quad b_3=- \frac{\theta_1}{\theta_1+\theta_\infty}\frac{\Gamma(\theta_1)}{\Gamma(\theta_\infty)}e^{-i\pi\theta_1},\\[2ex]
&c_1=\frac{\tilde{k}_2^0}{\tilde{k}_1^0}\frac{2y'_0(\theta_\infty-1)+\theta_\infty(\theta_1-1)}{2\theta_1\theta_\infty},\quad c_2=\frac{2y'_0(\theta_\infty-1)+\theta_\infty(\theta_1-1)}{2\theta_1\theta_\infty}.
\end{align*}
For the fundamental solution $\widetilde{Y}_2$ we use formulas~\eqref{asymptoticsU} and~\eqref{asymptoticsM} with $\epsilon=1$:
\[
\widetilde{Y}_2(z)=
\begin{pmatrix}
c_1\,z^{-\theta_2} & c_2\,z^{-\theta_2} & 0\\[1ex]
\tilde{y}_2(z;0,b_1,0,0) & \tilde{y}_2(z;0,b_2,0,0) & \tilde{y}_2(z;a,b_3e^{2\pi i\theta_2},0,0)
\\[1ex]
(\widetilde{Y}_1(z))_{31} & (\widetilde{Y}_1(z))_{32} & (\widetilde{Y}_1(z))_{33}
\end{pmatrix}
.
\]
For the fundamental solution $\widetilde{Y}_3$ we have to use the cyclic relation~\eqref{cyclicU} with $n=-1$ and again formulas~\eqref{asymptoticsU} and~\eqref{asymptoticsM} with $\epsilon=-1$ to obtain the asymptotics of the function $U$ in the sector $\mathcal{S}_3$:
\[
\widetilde{Y}_3(z)=
\begin{pmatrix}
c_1\,z^{-\theta_2} &  c_2\,z^{-\theta_2} & 0\\[1ex]
\tilde{y}_2(z;0,b_1e^{-2\pi i\theta_2},0,-1) & \tilde{y}_2(z;0,b_2e^{-2\pi i \theta_2},0,-1) & (\widetilde{Y}_3(z))_{23}\\[1ex]
 (\widetilde{Y}_3(z))_{31} &  (\widetilde{Y}_3(z))_{32} &  (\widetilde{Y}_3(z))_{33}
\end{pmatrix}
\quad z\in\mathcal{S}_3.
\]
The non trivial entries  $(\mathbb{S}_1)_{13}$ and $(\mathbb{S}_1)_{23}$ can be computed from the entries $(1,3)$ and $(2,3)$ of the equation $\widetilde{Y}_2(z)=\widetilde{Y}_1(z)\mathbb{S}_1$, while the non trivial entries $(\mathbb{S}_2)_{31}$ and $(\mathbb{S}_2)_{32}$ can be computed from the entries $(2,1)$ and $(2,2)$ of the equation $\widetilde{Y}_3(z)=\widetilde{Y}_2(z)\mathbb{S}_2$, obtaining

 \begin{equation}
\label{s1}
{\small \mathbb{S}_1=
\begin{pmatrix}
1 & 0 & \displaystyle{2i\tilde{k}_1^0\frac{\sin(\pi\theta_1)}{\sqrt{\theta_\infty}}\frac{\Gamma(\theta_1)}{\Gamma(\theta_\infty)}}\\[2ex]
0 & 1 & \displaystyle{-2i\tilde{k}_2^0\frac{\sin(\pi\theta_1)}{\sqrt{\theta_\infty}}\frac{\Gamma(\theta_1)}{\Gamma(\theta_\infty)}}\\[2ex]
0 & 0 & 1
\end{pmatrix}
}
\end{equation}
and
\begin{equation}
\label{s2}
{\small 
\mathbb{S}_2=
\begin{pmatrix}
1 & 0 & 0
\\[1ex]
0 & 1 & 0
\\[2ex]
\displaystyle{\frac{i\pi e^{-i\pi(\theta_1-\theta_\infty)}}{\tilde{k}_1^0\sqrt{\theta_\infty}}\frac{2y'_0(\theta_\infty-1)+\theta_\infty(\theta_1-1)}{\Gamma(1+\theta_1)\Gamma(1-\theta_\infty)}} 
&
 \displaystyle{\frac{i\pi e^{-i\pi (\theta_1-\theta_\infty)}}{\tilde{k}_2^0\sqrt{\theta_\infty}}\frac{2y'_0(\theta_\infty-1)-\theta_\infty(\theta_1+1)}{\Gamma(1+\theta_1)\Gamma(1-\theta_\infty)}} & 1
\end{pmatrix}
}
.
\end{equation}

We are ready to apply the formulae \eqref{3febbraio2021-2}.  For $\tau=0$ and $\eta=3\pi/2$, the ordering relation is $1\prec 3$ and $2\prec 3$, while there is no ordering $1 \leftrightarrow 2$, because $u_1-u_2\to 0$. So we have $p_{12}=2$ and 
\begin{align*}
&p_{13}= 2\cos\pi(\theta_1-\theta_3)-e^{i\pi(\theta_1-\theta_3)} (\mathbb{S}_1)_{13}(\mathbb{S}_{2}^{-1})_{31}, \quad  1\prec 3,
\\ 
&p_{23}= 2\cos\pi(\theta_2-\theta_3)-e^{i\pi(\theta_2-\theta_3)} (\mathbb{S}_1)_{23}(\mathbb{S}_{2}^{-1})_{32}, \quad  2\prec 3.
\end{align*}
Substituting the entries of the Stokes matrices, we exactly obtain (and confirm!) the known result of theorem 3 of \cite{guz2006} or formulas in the footnote 4 at page 3258 of \cite{guz2012}, namely $p_{12}=2$ and 
\begin{align}
\label{3maggio2021-7} 
  & 
p_{13}=-\frac{4s\sin(\pi\theta_1) \sin(\pi\theta_\infty)}{\theta_1}+2\cos(\pi(\theta_1+\theta_\infty)),
\\
\label{3maggio2021-8}
& p_{23}=\frac{4s\sin(\pi\theta_1)\sin(\pi\theta_\infty)}{\theta_1}+2\cos(\pi(\theta_1-\theta_\infty)) ,
\end{align}
where the parameter $s$ is equivalent to $y'_0$ through
$$ 
y_0'= \frac{\theta_\infty(2s+\theta_1+1)}{2(\theta_\infty-1)}. 
$$

{\small Notice that we receive \eqref{3maggio2021-7} and \eqref{3maggio2021-8} also  if we choose $\tau=\pi$ and $\eta=\pi/2$. In this case,  the ordering relation is $3\prec 1$ and $3\prec 2$.  Then,  
$$
p_{13}= 2\cos\pi(\theta_1-\theta_3)-e^{i\pi(\theta_3-\theta_1)}(\mathbb{S}_2)_{31}(\mathbb{S}_{3}^{-1})_{13} \quad 1\succ 3,
$$
$$
p_{23}= 2\cos\pi(\theta_2-\theta_3)-e^{i\pi(\theta_3-\theta_2)}(\mathbb{S}_2)_{32}(\mathbb{S}_{3}^{-1})_{23} \quad 2\succ 3.
$$ 
We need in this case the Stokes matrices $\mathbb{S}_2$ and $\mathbb{S}_3$, the latter being obtainable  from $\mathbb{S}_1$ by the formulae 
\[
\mathbb{S}_{2\nu+1}=e^{2\pi i\nu\Theta}\mathbb{S}_1e^{-2\pi i\nu\Theta},\quad\mathbb{S}_{2\nu}=e^{2\pi i(\nu-1)\Theta}\mathbb{S}_1 e^{-2\pi i(\nu-1)\Theta},\quad\nu\in\mathbb{Z},
\] 
where $\Theta:=\operatorname{diag}(\theta_1,\theta_2,\theta_3)\equiv\operatorname{diag}(\theta_1,\theta_1,-\theta_\infty)$. 
}
\bre
{\rm 
The Stokes matrices~\eqref{s1} and~\eqref{s2} are the same as those for case (T1) of the table with $N=-1$, $\theta_1=\theta_2$, $\theta_3=-\theta_\infty$, since the system~\eqref{sysex1} that we have considered in this section is gauge equivalent to system~\eqref{sysex1bis} via the gauge transformation
\[
\widetilde{Y}\mapsto P\widetilde{Y},\quad P=\begin{pmatrix} 0 & 0 & 1\\ 0 & 1 & 0\\ 1 & 0 & 0 \end{pmatrix}.
\]
}
\ere

%%%%%%%%%%%%%%%

\subsection{Case (T3) of the table}
\label{StokesT3}

We consider the trascendent in case (T3) of the table
\be
\label{29luglio2021-5}
y(x)=y_0^\prime +\frac{y_0^\prime(y_0^\prime-1)(\theta_3^3-1-(\theta_\infty-1)^2)}{2}x^2+O(x^3),\quad \theta_1=\theta_2=0,
\ee
with  $\theta_3,\theta_\infty$ not integers (the  case $\theta_3=0$ and $\theta_\infty\in\mathbb{Z}\setminus\{0\}$ will be considered in Section \ref{dfmanifolds}).  Since $\theta_1-\theta_2\not\in\mathbb{Z}\backslash\{0\}$ there is no partial resonance, so that the formal solution of system~\eqref{syst3} is unique. It is computed following Section 4 of \cite{CDG}, receiving 
\[
\widetilde{Y}_F(z)=\left(I+O\left(\frac{1}{z}\right)\right)
\begin{pmatrix}
\displaystyle{-y'_0\frac{\theta_3-\theta_\infty}{2\tilde{k}_1^0\theta_\infty}} & \displaystyle{-(1-y'_0)\frac{\theta_3-\theta_\infty}{2\tilde{k}_2^0\theta_\infty}} & \displaystyle{-\frac{\theta_3-\theta_\infty}{2\theta_\infty}z^{-\theta_3}e^z}\\[2ex]
\displaystyle{-y_0'\frac{\tilde{k}_2^0}{\tilde{k}_1^0}} & y'_0 & 0\\[2ex]
\displaystyle{y'_0\frac{\theta_3-\theta_\infty}{2\tilde{k}_1^0\theta_\infty}} & \displaystyle{(1-y'_0)\frac{\theta_3-\theta_\infty}{2\tilde{k}_2^0\theta_\infty}} & \displaystyle{\frac{\theta_3+\theta_\infty}{2\theta_\infty}z^{-\theta_3}e^z}
\end{pmatrix},
\]
where $\tilde{k}_1^0,\tilde{k}_2^0$ are in \eqref{25giugno2021-3}.  The elements of the second row of the solutions of~\eqref{syst3} are constants, while from~\eqref{y1t3} the elements of the first row have the general form
\begin{align*}
&\tilde{y}_1(z;a,b,m,n)=\left[aM\left(ze^{2m\pi i};\frac{\theta_\infty-\theta_3}{2},\theta_\infty\right)+bU\left(ze^{2n\pi i};\frac{\theta_\infty-\theta_3}{2},\theta_\infty\right)\right]z^{(\theta_\infty-\theta_3)/2},
\\ & a,b\in\mathbb{C},n,m\in\mathbb{Z}.
\end{align*}
It is sufficient to compute the first two rows of the fundamental solutions $\widetilde{Y}_1,\widetilde{Y}_2,\widetilde{Y}_3$ with asymptotics $\widetilde{Y}_F$ in the Stokes sectors $\mathcal{S}_1,\mathcal{S}_2,\mathcal{S}_3$, respectively. The computation of the second row is immediate by comparison with the second row of the leading term of $\widetilde{Y}_F$. The first row is obtained by comparison of the first row of the leading term of $\widetilde{Y}_F$ and the leading coeffincients of the asymptotics of the confluent hypergeometric functions $M$ and $U$. For the fundamental solution $\widetilde{Y}_1$ we use formulas~\eqref{asymptoticsU} and~\eqref{asymptoticsM} with $\epsilon=-1$:
\[
\widetilde{Y}_1(z)=
\begin{pmatrix}
\tilde{y}_1(z;0,b_1,0,0) & \tilde{y}_1(z;0,b_2,0,0) & \tilde{y}_1(z;a_3,b_3,0,0)\\[1ex]
\displaystyle{-y_0'\frac{\tilde{k}_2^0}{\tilde{k}_1^0}} & y'_0 & 0\\[2ex]
(\widetilde{Y}_1)_{31} & (\widetilde{Y}_1)_{32} & (\widetilde{Y}_1)_{33} 
\end{pmatrix}
,
\]
where
\begin{align*}
&a_3=-\Gamma\left(\frac{\theta_\infty-\theta_3}{2}\right)\frac{\theta_3-\theta_\infty}{2\theta_\infty},\\[2ex]
&b_1=-y'_0\frac{\theta_3-\theta_\infty}{2\tilde{k}_1^0\theta_\infty},\quad b_2=-(1-y'_0)\frac{\theta_3-\theta_\infty}{2\tilde{k}_2^0\theta_\infty},\quad b_3=\frac{\Gamma((\theta_\infty-\theta_3)/2)}{\Gamma((\theta_\infty+\theta_3)/2)}\frac{\theta_3-\theta_\infty}{2\theta_\infty}e^{-i\pi(\theta_\infty-\theta_3)/2}.
\end{align*}
For the fundamental solution $\widetilde{Y}_2$ we use formulas~\eqref{asymptoticsU} and~\eqref{asymptoticsM} with $\epsilon=1$:
\[
\widetilde{Y}_2(z)=
\begin{pmatrix}
\tilde{y}_2(z;0,b_1,0,0) & \tilde{y}_1(z;0,b_2,0,0) & \tilde{y}_1(z;a,b_3e^{i\pi(\theta_\infty-\theta_3)/2},0,0)\\[1ex]
\displaystyle{-y_0'\frac{\tilde{k}_2^0}{\tilde{k}_1^0}} & y'_0 & 0\\[2ex]
(\widetilde{Y}_2)_{31} & (\widetilde{Y}_2)_{32} & (\widetilde{Y}_2)_{33} 
\end{pmatrix}
.
\]
For the fundamental solution $\widetilde{Y}_3$ we have to use the cyclic relation~\eqref{cyclicU} with $n=-1$ and again formulas~\eqref{asymptoticsU} and~\eqref{asymptoticsM} with $\epsilon=-1$ to obtain the asymptotics of the function $U$ in the sector $\mathcal{S}_3$:
\[
\widetilde{Y}_3(z)=
\begin{pmatrix}
\tilde{y}_1(z;a_1,\tilde{b}_1,0,0) & \tilde{y}_1(z;a_2,\tilde{b}_2,0,0) & \tilde{y}_1(z;a_3e^{-i\pi(\theta_\infty+\theta_3)},\tilde{b}_3,-1,-1)\\[1ex]
\displaystyle{-y_0'\frac{\tilde{k}_2^0}{\tilde{k}_1^0}} & y'_0 & 0\\[2ex]
(\widetilde{Y}_3(z))_{31} & (\widetilde{Y}_3(z))_{32} & (\widetilde{Y}_3(z))_{33}
\end{pmatrix}
,
\]
where
\begin{align*}
&a_1=2iy'_0\frac{\theta_3-\theta_\infty}{2\tilde{k}_1^0\theta_\infty}\Gamma\left(\frac{\theta_\infty+\theta_3}{2}\right)\sin\left[\frac{\pi}{2}(\theta_\infty+\theta_3)\right]e^{i\pi\theta_3},
\\
\noalign{\medskip}
&a_2=2i(1-y'_0)\frac{\theta_3-\theta_\infty}{2\tilde{k}_2^0\theta_\infty}\Gamma\left(\frac{\theta_\infty+\theta_3}{2}\right)\sin\left[\frac{\pi}{2}(\theta_\infty+\theta_3)\right]e^{i\pi \theta_3},
\\
\noalign{\medskip}
&\tilde{b}_1=b_1e^{i\pi(\theta_\infty+\theta_3)},\quad \tilde{b}_2=b_2e^{i\pi(\theta_\infty+\theta_3)},\quad \tilde{b}_3=b_3e^{-i\pi(\theta_\infty+\theta_3)}
\end{align*}
The non trivial entries  $(\mathbb{S}_1)_{13}$ and $(\mathbb{S}_1)_{23}$ can be computed from the entries $(1,3)$ and $(2,3)$ of the equation $\widetilde{Y}_2(z)=\widetilde{Y}_1(z)\mathbb{S}_1$, while the non trivial entries $(\mathbb{S}_2)_{31}$ and $(\mathbb{S}_2)_{32}$ can be computed from the entries $(1,1)$ and $(1,2)$ of the equation $\widetilde{Y}_3(z)=\widetilde{Y}_2(z)\mathbb{S}_2$, obtaining
\begin{equation}
\label{25giugno2021-5}
{\small
\mathbb{S}_1=
\begin{pmatrix}
1 & 0 & \displaystyle{-2i\tilde{k}_1^0\frac{\Gamma((\theta_\infty-\theta_3)/2)}{\Gamma((\theta_\infty+\theta_3)/2)}\sin\left[\frac{\pi}{2}(\theta_\infty-\theta_3)\right]}\\[2ex]
0 & 1 & \displaystyle{-2i\tilde{k}_2^0\frac{\Gamma((\theta_\infty-\theta_3)/2)}{\Gamma((\theta_\infty+\theta_3)/2)}\sin\left[\frac{\pi}{2}(\theta_\infty-\theta_3)\right]}\\[2ex]
0 & 0 & 1
\end{pmatrix}
}
\end{equation}
and
\begin{equation}
\label{25giugno2021-6}
{\small 
\mathbb{S}_2=
\begin{pmatrix}
1 & 0 & 0\\[2ex]
0 & 1 & 0\\[2ex]
\displaystyle{-\frac{2i\pi y_0'e^{i\pi\theta_3}}{\tilde{k}_1^0\Gamma(1-(\theta_\infty+\theta_3)/2)\Gamma((\theta_\infty-\theta_3)/2)}} & \displaystyle{-\frac{2i\pi (1-y_0')e^{i\pi\theta_3}}{\tilde{k}_2^0\Gamma(1-(\theta_\infty+\theta_3)/2)\Gamma((\theta_\infty-\theta_3)/2)}} & 1
\end{pmatrix}
}
.
\end{equation}
With the choice of the admissible direction $\tau=0$, corresponding to $\eta=3\pi/2$, formulas~\eqref{3febbraio2021-2} become
\begin{align*}
& p_{12}=2,
& p_{13}=2\cos(\pi\theta_3)-4y_0'\sin\left[\frac{\pi}{2}(\theta_\infty-\theta_3)\right]\sin\left[\frac{\pi}{2}(\theta_\infty+\theta_3)\right],
\\
&&p_{23}=2\cos(\pi\theta_3)-4(1-y_0')\sin\left[\frac{\pi}{2}(\theta_\infty-\theta_3)\right]\sin\left[\frac{\pi}{2}(\theta_\infty+\theta_3)\right],
\end{align*}
The above expressions confirm the  results of \cite{guz2012}, page 3260, Case (46).

%%%%%%%%%%%%%%%%%%%%%%%%%%%%%%%%%%%%%%%%%%%%%%

\subsection{Monodromy data of three dimensional semisimple Frobenius manifolds associated to transcendents (T3) of the Table}
\label{dfmanifolds}

A semisimple Frobenius manifold $M$ (of dimension $n$) is a complex analytic manifold whose tangent bundle is equipped with a Frobenius algebra structure,  semisimple on an open dense subset $M_{ss}\subset M$, and with a $z$-deformed flat connection (see \cite{Dub1,Dub2,CDG1,Her} for the details). In suitable coordinates $u=(u^1,\ldots,u^n)$, called \textit{canonical}, the flatness condition is equivalent to the Frobenius  integrability of  the $n$-dimensional Pfaffian system 
\[
dY=\omega(z,u)Y,\quad\omega(z,u)=\left(U(u)+\frac{V(u)}{z}\right)dz+\sum_{k=1}^n\left(zE_k+V_k(u)\right)du_k,
\]
where $U(u)=\operatorname{diag}(u^1,\ldots,u^n)$, $(E_k)_{ij}=\delta_{ik}\delta_{jk}$, $V(u)$ is skew-symmetric and
\[
V_k(u)=\frac{\partial\Psi(u)}{\partial u^k}\Psi^{-1}(u),
\]
where $\Psi(u)$ diagonalizes holomorphically $V(u)$ on $M_{ss}$. For  geometrical reasons, $\Psi^{\rm T}\cdot \Psi$ is constant, satisfying  a  normalization condition that in dimension 3 is given below in \eqref{8luglio2021-1}. 
 The monodromy data of the $z$-component of the above Pfaffian system   locally parametrize the manifold, whose structure can be locally reconstructed by solving a Riemann-Hilbert boundary value problem \cite{Dub1,Dub2}. This holds also at a semisimple coalescence point \cite{CDG1,sabbah,Cotti20}, where   $V(u)$ is holomorphic and  the condition 
$
V_{ij}(u)\to 0$ for $u_i-u_j\to 0$ holds. 

 For three dimensional semisimple Frobenius manifolds the Pfaffian system is a special case of the one we have studied in this paper, with 
 $$ 
 \theta_1=\theta_2=\theta_3=0,\quad\quad 
  \hbox{eigenvalues of }V = \mu,0,-\mu,\quad \quad   \mu:= \theta_\infty/2.
  $$
  The factorization \eqref{20gennaio2021-1} implies
  $$ 
  V(u)=\Omega(x).
  $$
  Moreover, the choice of  Remark \ref{11giugno2021-3} must be made, so that $V$ is skew-symmetric (see also \eqref{8luglio2021-2} below).  
  As explained in \cite{Dub2} and translated into explicit formulae in \cite{guzz2001}, the 
  structure of a 3-dimensional  Dubrovin-Frobenius manifold is expressed in terms of a PVI transcendent,  associated   to the equivalent class of the matrix  $V= \Omega$, as  in our Theorem \ref{20gennaio2021-2}.\footnote{  In   \cite{guzz2001}, in place of $\Omega(x)$, a matrix  
  $ 
  V(x)= \begin{pmatrix} 0 & -\Omega_3(x) & ~\Omega_2(x)
  \\
 ~ \Omega_3(x) & 0 & -\Omega_1(x) 
  \\
  -\Omega_2(x) & ~\Omega_1(x) & 0 
  \end{pmatrix}
  $
 appears,  while $U=\hbox{\rm diag}(0,1,x)$. The matrix $\Omega$ used here is $\Omega(x)=P^{-1}V(x) P$, where 
  $
  P=\begin{pmatrix} 1 & 0 & 0
  \\
  0 & 0 & 1
  \\
  0 & 1& 0 
  \end{pmatrix}
  $. In   \cite{guzz2001},  a parameterization of $V(x)$ in terms of $y(x)$ was already given, so anticipating for  the special case $\theta_1=\theta_2=\theta_3=0$ our Theorem \ref{20gennaio2021-2} and the analogous result of \cite{MazzoIrr}. 
  } 
  Every branch of the transcendent parameterizes a chamber of the manifold.  By the results in \cite{CDG1}, this fact extends at  semisimple coalescence points.  
     
     We consider here a Dubrovin-Frobenius manifold with a local chamber associated with a  branch    holomorphic at  $x=0$, with  the Taylor expansion  
\be
\label{29luglio2021-6}
     y(x)=y_0^\prime x+y_0^\prime(1-y_0^\prime)(2\mu^2-2\mu+1)x^2+O(x^3),\quad \quad y_0^\prime\neq 0,1.
    \ee
      We can compute the monodromy data of  the  local chamber above, namely the data associated with $y(x)$, using  system \eqref{mainsystem}. This is possible  because $y(x)$  is  in the class  (T3) of the table (with $\theta_1=\theta_2=\theta_3=0$, $\theta_\infty=2\mu$). 
    In \cite{CDG} an example of this computation is given for $\theta_\infty=-1/2$ and $y_0'=1/2$ (corresponding to a chamber of  the  Frobenius manifold on the space of orbits of the Coxeter group $A_3$), while here we perform the computation for all values of $\theta_\infty\ne 0$ and $y_0'\ne 0$. 
  
   The matrix $\Omega(0)=\Omega_0$   is \eqref{25giugno2021-3}  of Section \ref{t3}, skew-symmetric only for   (up to the sign)  
   \be
   \label{8luglio2021-2}
   \tilde{k}_1^0=i\sqrt{y_0^\prime}, \quad\tilde{k}_2^0=i\sqrt{1-y_0^\prime},
   \ee
namely
\[
\Omega_0=\begin{pmatrix}
0 & 0 &\Omega_{13}
\\[2ex]
0 & 0 & \Omega_{23}
\\[2ex]
-\Omega_{13} & -\Omega_{23} & 0
\end{pmatrix}
=
\begin{pmatrix}
0 & 0 & i\mu \sqrt{y_0^\prime}
\\[2ex]
0 & 0 & i\mu\sqrt{1-y_0^\prime}
\\[2ex]
-i\mu \sqrt{y_0^\prime} & -i\mu\sqrt{1-y_0^\prime} & 0
\end{pmatrix}
,
\]
The  diagonalizing matrix   $G_0$, given in \eqref{g0t3}, with the parameters $\tilde{k}_1^0,\tilde{k}_2^0$  fixed above, must be renormalized  by 
$$G_0\longmapsto \Psi:=G_0\cdot \hbox{\rm diag}(1/\sqrt{2},\Omega_{13}/(i\mu),1/\sqrt{2}) =G_0\cdot \hbox{\rm diag}(1/\sqrt{2},\sqrt{y_0^\prime},1/\sqrt{2}),
 $$ 
 so that   
 \be 
 \label{8luglio2021-1}
 \Psi^{\rm T} \cdot \Psi=\eta:= \begin{pmatrix}
 0 & 0 & 1
 \\
 0 & 1 & 0 
 \\
 1 & 0 & 0
 \end{pmatrix},\quad\quad \Psi^{-1}\Omega_0\Psi=\operatorname{diag}(\mu,0,-\mu)=:\hat{\mu},\quad\quad \mu\neq 0. 
 \ee
 This normalization does not affect the computation of Stokes matrices, but is important in order to compute   the central connection matrix, as will be done below,  in accordance with  \cite{Dub1,Dub2,CDG1}. 
Explicitly, using $(\Omega_{13})^2+(\Omega_{23})^2=-\mu^2$ we obtain
% \[
%\Psi=
%\begin{pmatrix}
%\displaystyle{-\frac{\Omega^0_3}{\mu}} & \displaystyle{\frac{\Omega^0_1}{\Omega^0_3}} & \displaystyle{\frac{\Omega^0_3}{\mu}} \\[2ex]
%\displaystyle{\frac{\Omega^0_1}{\mu}} & 1 & -\displaystyle{\frac{\Omega^0_1}{\mu}} \\[2ex]
%1 & 0 & 1
%\end{pmatrix}
%\mbox{,}
%\]
%with inverse
%\[
%\Psi^{-1}=
%\begin{pmatrix}
%\displaystyle{\frac{\Omega^0_3}{2\mu}} &-\displaystyle{\frac{ \Omega^0_1}{2\mu}} & \displaystyle{\frac{1}{2}} \\[2ex]
%-\displaystyle{\frac{\Omega^0_1\Omega^0_3}{\mu^2}} & -\displaystyle{\left(\frac{\Omega^0_3}{\mu}\right)^2} & 0\\[2ex]
%-\displaystyle{\frac{\Omega^0_3}{2\mu}} & \displaystyle{\frac{\Omega^0_1}{2\mu}} & \displaystyle{\frac{1}{2}}
%\end{pmatrix}
%.
%\]
{\small 
\begin{align*}
\Psi=\begin{pmatrix}
\dfrac{\Omega_{13}}{\sqrt{2}\mu}& \dfrac{i\Omega_{23}}{\mu} & -\dfrac{\Omega_{13}}{\sqrt{2}\mu}
\\
\noalign{\medskip}
\dfrac{\Omega_{23}}{\sqrt{2}\mu}& -\dfrac{i\Omega_{13}}{\mu} & -\dfrac{\Omega_{23}}{\sqrt{2}\mu}
\\
\noalign{\medskip}
\dfrac{1}{\sqrt{2}} & 0 & \dfrac{1}{\sqrt{2}}
\end{pmatrix}
 = \begin{pmatrix}
i\sqrt{y_0^\prime/2} & -\sqrt{1-y_0^\prime} & -i\sqrt{y_0^\prime/2}
\\
\noalign{\medskip}
i\sqrt{(1-y_0^\prime)/2} & \sqrt{y_0^\prime} & -i\sqrt{(1-y_0^\prime)/2}
\\
\noalign{\medskip}
1/\sqrt{2} & 0 & 1/\sqrt{2}
\end{pmatrix},
\\
\noalign{\medskip}
\Psi^{-1}=
\begin{pmatrix}
-\dfrac{\Omega_{13}}{\sqrt{2}\mu} & -\dfrac{\Omega_{23}}{\sqrt{2}\mu} & \dfrac{1}{\sqrt{2}}
\\
\noalign{\medskip}
\dfrac{i\Omega_{23}}{\mu} & -\dfrac{i \Omega_{13}}{\mu} & 0 
\\
\noalign{\medskip}
\dfrac{\Omega_{13}}{\sqrt{2}\mu} & \dfrac{\Omega_{23}}{\sqrt{2}\mu} & \dfrac{1}{\sqrt{2}}
\end{pmatrix}
= \begin{pmatrix}
-i\sqrt{y_0^\prime/2} & -i\sqrt{(1-y_0^\prime )/2} & 1/\sqrt{2}
\\
\noalign{\medskip}
-\sqrt{1-y_0^\prime} & \sqrt{y_0^\prime} & 0 
\\
\noalign{\medskip}
i\sqrt{y_0^\prime/2} & i \sqrt{(1-y_0^\prime)/2} & 1/\sqrt{2}
 \end{pmatrix}.
\end{align*}
}
In conclusion, system \eqref{syst3}  for $Y=\Psi \widetilde{Y}$ becomes
\begin{equation}
\label{normal}
\frac{d\widetilde{Y}}{dz}=\left[
\mathcal{U}_0
+
\frac{\mu}{z}
\begin{pmatrix}
1 & 0 & 0\\
0 & 0 & 0\\
0 & 0 & -1
\end{pmatrix}
\right]\widetilde{Y},\quad\quad 
\mathcal{U}_0:=\frac{1}{2}
\begin{pmatrix}
1 & 0 & 1\\
0 & 0 & 0\\
1 & 0 & 1
\end{pmatrix}
\mbox{.}
\end{equation}

\subsubsection{Stokes matrices and data $p_{jk}$}

\bpr
\label{29luglio2021-3}
The Stokes matrices  for the chamber of a  semisimple Dubrovin-Frobenius manifold  associated with the branch \eqref {29luglio2021-6}, with $\mu\neq 0$, are 
\be
\label{29luglio2021-1}
\mathbb{S}_1=
\begin{pmatrix}
1 & 0 & 2\sqrt{y_0'}\sin(\mu\pi)\\[2ex]
0 & 1 & 2\sqrt{1-y_0'}\sin(\mu\pi)\\[2ex]
0 &  0 & 1
\end{pmatrix}
,\quad
\mathbb{S}_2=
\begin{pmatrix}
1 & 0 &0\\[2ex]
0 & 1 &0\\[2ex]
-2\sqrt{y_0'}\sin(\mu\pi) &  -2\sqrt{1-y_0'}\sin(\mu\pi) & 1
\end{pmatrix}
.
\ee
Corresponding to them, 
$$p_{12}=2,\quad
\quad
p_{13}= 2-[(\mathbb{S}_1)_{13}]^2=2-4y_0^\prime \sin^2(\pi\mu) ,\quad \quad p_{23}=2-[(\mathbb{S}_1)_{23}]^2=2-4(1-y_0^\prime)\sin^2(\pi\mu).
$$

\epr 

\begin{proof}
If $2\mu\not\in\mathbb{Z}$, this is a particular case of the computations of Section \ref{StokesT3}, with $\theta_1=\theta_2=\theta_2=0$, $\theta_\infty=2\mu$ in the Stokes matrices \eqref{25giugno2021-5}-\eqref{25giugno2021-6}. 
If $2\mu\in\mathbb{Z}$, we need two facts. The first is that  the symmetries of the Pfaffian system imply \cite{Dub1}
\be
\label{3agosto2021-1}
\mathbb{S}_2=\mathbb{S}_1^{-T}.
\ee
 The second fact is that for a Dubrovin-Frobenius manifold, in order to compute $\mathbb{S}_1$, it suffices  to have $\hat{\mu}$, the nilpotent exponent  $R$ (defined below) and  the central connection matrix $\mathcal{C}^{(0)}$ (also defined below). Indeed,  the following  relations  hold (see \cite{Dub1,CDG1})
\be
\label{9luglio2021-1}
\mathbb{S}_1=(\mathcal{C}^{(0)})^{-1} e^{-i\pi R} e^{-i\pi \hat{\mu}} \eta^{-1} (\mathcal{C}^{(0)})^{-\rm{T}},\quad\hbox{ or }\quad \mathbb{S}_1^{\rm T}=(\mathcal{C}^{(0)})^{-1} e^{i\pi R}e^{i\pi \hat{\mu}} \eta^{-1} (\mathcal{C}^{(0)})^{-\rm{T}}.
\ee
 We will compute $\mathcal{C}^{(0)}$ for all values of $\mu\neq 0$ in Section \ref{29luglio2021-2} below, so  proving that $\mathbb{S}_1$ and $\mathbb{S}_2$ are as in \eqref{29luglio2021-1} in all cases.   
  The invariant $p_{jk}$ follow from Theorem \ref{16aprile2021-5}.
\end{proof}

The invariants $p_{12},p_{13},p_{23}$ in Proposition \ref{29luglio2021-3} coincide with  the data obtained by the direct analysis of the $2\times2$ Fuchsian system, at the bottom of page 1335 of \cite{guz2002} (see the formulae for $\sigma=0$ there, where there is a misprint, to be corrected by the replacement:  $x_1^2\mapsto x_1$, $x_\infty^2\mapsto x_\infty$). 

\subsubsection{A comment on the transcendent \eqref{29luglio2021-6} with integer $\mu$}
\label{13agosto2021-1}

 If  $\mu\in\mathbb{Z}\backslash\{0\}$, 
 $$ 
 \mathbb{S}_1=\mathbb{S}_2=I,
\quad\quad
p_{12}=p_{13}=p_{23}=2.
$$ 
The $p_{jk}=2$ correspond to a  degenerate case mentioned in the Introduction, when  the triple $\mathcal{M}_1,\mathcal{M}_2,\mathcal{M}_3$ and the the integration constant in the transcendent cannot be parametrized by  the $p_{jk}$.   For example, if $\mu=1$, the series  \eqref{29luglio2021-6} is the Taylor expansion at $x=0$ of the  one-parameter family of rational solutions   (see  \cite{Hit,MazRat})
\be
\label{29luglio2021-8}
y(x)= \frac{a x}{1-(1-a)x},\quad \quad a=y_0^\prime\in\mathbb{C}\backslash\{0\},\quad\quad \mu=1
.
\ee
The  Fuchsian system \eqref{20gennaio2021-21} or \eqref{23gennaio2021-8} has a fundamental solution\footnote{System \eqref{23gennaio2021-8} with matrix coefficients determined by \eqref{29luglio2021-8}  has the isomonodromic  fundamental matrix solutions  $\widetilde{\Phi}(\lambda,x) C$, $\det C\neq 0$, with 
$$ 
\widetilde{\Phi}(\lambda,x)= 
\begin{pmatrix}
1 & L(\lambda,x) 
\\
\noalign{\medskip}
r(\lambda,x) &\dfrac{(1 + (a-1)x)^2}{k0}+r(\lambda,x) L(\lambda,x)
\end{pmatrix}
,
\quad\quad 
k_0\neq 0,
$$
where 
$$L(\lambda,x)=-a\ln(\lambda)  + ( a-1)\ln(\lambda - x)+ \ln(\lambda - 1),\quad\quad  r(\lambda,x)= \frac{2((a-1)x + 1)\lambda + (1 - a)x^2 - 1}{2k_0}.
$$
 }  with monodromy matrices 
$$ 
\mathcal{M}_1= \begin{pmatrix} 1 & -2i\pi a \\ 0 & 1 \end{pmatrix},
\quad
\mathcal{M}_2= \begin{pmatrix} 1 & 2i\pi (a-1) \\ 0 & 1 \end{pmatrix},
\quad
\mathcal{M}_2= \begin{pmatrix} 1 & 2\pi i  \\ 0 & 1 \end{pmatrix},
$$
generating  a reducible monodromy group. Notice that $a$ appears in $y(x)$ and  the matrices, but not in the traces $p_{jk}$. 
For other integer values of  $\mu$, we obtain the corresponding solutions with  expansion   \eqref{29luglio2021-6}, applying to \eqref{29luglio2021-8}  the symmetry of PVI given in Lemma 1.7 of \cite{DM}, which transforms $y(x)$ for a PVI with given $\mu$, to a solution $\tilde{y}(x)$ of a PVI with $-\mu$ or equivalently $1+\mu$. In this way, all values $\mu+N$, $N\in\mathbb{Z}$, are obtained. For example, from \eqref{29luglio2021-8}, we obtain
$$ 
\tilde{y}(x)= \frac{ax(1 + (a-1)x^2)^2}{(1 + x(a-1))(1 + x(a-1)(ax^3 - x^3 + 4x^2 - 6x + 4))},\quad \quad \mu=-1,2.
$$

{\small
We also mention that for $\mu$ integer,  the system 
$$ 
\frac{dY}{dz}=\Bigl(\hbox{\rm diag}(0,x,1)+\frac{\Omega(x)}{z}\Bigr) Y
$$ 
is solvable in terms of elementary functions and has fundamental matrix solutions without Stokes phenomenon, so that $ \mathbb{S}_1=\mathbb{S}_2=I$. For example, for $\mu=1$, corresponding to the transcendent \eqref{29luglio2021-8}, we have $k_1(x)=k_1^0$, $k_2(x)=k_2^0$. For the choice $k_1^0=i\sqrt{a}$, $k_2^0=i\sqrt{1-a}$ (the choice is fixed up to signs) we obtain  a skew-symmetric   
$$\Omega(x)
=\frac{1}{1+(a-1)x}\begin{pmatrix}
0 & -ix\sqrt{a(a-1)} & -i\sqrt{a}
\\
ix\sqrt{a(a-1)} &  0 & (x-1) \sqrt{a-1}
\\
i\sqrt{a} & -(x-1) \sqrt{a-1} & 0 
\end{pmatrix}.
$$
There is an elementary  fundamental matrix solution with canonical behaviour at $z=\infty$ and trivial Stokes matrices: $$ 
Y(z,x)=\left[I+\frac{1}{(ax-x + 1)}
\begin{pmatrix}
a & -i\sqrt{a(a - 1)} & -I\sqrt{a}
\\
 -i\sqrt{a(a - 1)} & 1-a & -\sqrt{a-1}
 \\
-i\sqrt{a} & -\sqrt{a-1} & -1 
\end{pmatrix}\cdot \frac{1}{z}\right]\begin{pmatrix}
1 & 0 & 0 
\\
0 & e^{xz} & 0 
\\
0 & 0 & e^z
\end{pmatrix}.
$$
Another explanation  for  $\mathbb{S}_1=\mathbb{S}_2=I$ is  that
$$
\lim_{x\to 0}\Omega_{12}(x)=\lim_{x\to 0}\Omega_{21}(x)=0\quad \Longrightarrow \quad (\mathbb{S}_j)_{12}= (\mathbb{S}_j)_{21}=0,
$$
$$
\lim_{x\to 1}\Omega_{23}(x)=\lim_{x\to 1}\Omega_{32}(x)=0\quad \Longrightarrow \quad (\mathbb{S}_j)_{23}= (\mathbb{S}_j)_{32}=0,
$$
$$
\lim_{x\to \infty}\Omega_{13}(x)=\lim_{x\to \infty}\Omega_{31}(x)=0\quad \Longrightarrow \quad (\mathbb{S}_j)_{13}= (\mathbb{S}_j)_{31}=0.
$$
The implications ``$\Longrightarrow$''  above follow from theorem 1.1 of \cite{CDG}, described in the Introduction. 
}

%%%%%%%%%%%%%
\subsubsection{Fundametal solution at the Fuchsian singularity $z=0$}
\label{29luglio2021-2}

From the general theory of Fuchsian singularities, system~\eqref{normal} has a fundamental solution at $z=0$ of the form
\be
\label{8luglio2021-4}
\widetilde{Y}^{(0)}(z)=\Bigl(I+\sum_{k=0}^\infty G_kz^k \Bigr)z^{\hat{\mu}}z^R,
\ee
where the series is convergent and  $R$ and the $G_k$'s are constructed recursively as follows:
\begin{itemize}
\item For $k=1$, if $\mu_i-\mu_j\ne 1$ then
\[
\begin{cases}
\mbox{choose }(R_1)_{ij}=0\\[1ex]
\displaystyle{(G_1)_{ij}=\frac{(\mathcal{U}_0)_{ij}}{\mu_j-\mu_i+1}},
\end{cases}
\]
else if $\mu_i-\mu_j=1$ then
\[
\begin{cases}
\mbox{necessarily }(R_1)_{ij}=(\mathcal{U}_0)_{ij}\\[1ex]
(G_1)_{ij}\mbox{ is arbitrary};
\end{cases}
\]
\item For $k\ge 2$, if $\mu_i-\mu_j\ne k$, then
\[
\begin{cases}
\mbox{choose }(R_k)_{ij}=0\\[1ex]
\displaystyle{(G_k)_{ij}=\frac{1}{\mu_j-\mu_i+k}\left(\mathcal{U}_0G_{k-1}-\sum_{p=1}^{k-1}G_pR_{k-p}\right)_{ij}},
\end{cases}
\]
else if $\mu_i-\mu_j=k$ then
\[
\begin{cases}
\mbox{necessarily }(R_k)_{ij}=\left(\mathcal{U}_0G_{k-1}-\sum_{p=1}^{k-1}G_pR_{k-p}\right)_{ij}\\[1ex]
(G_k)_{ij}\mbox{ is arbitrary};
\end{cases}
\]
\end{itemize}
here $\mu_1=\mu$, $\mu_2=0$, $\mu_3=-\mu$ are the diagonal elements of $\hat{\mu}$. The nilpotent matrix $R$ is
\[
R=\sum_{k=1}^\infty R_k \quad \hbox{ finite sum}.
\]

\bpr
If  system \eqref{normal} is non-resonant (i.e. $2\mu\not\in\mathbb{Z}\backslash \{0\}$), then it has a fundamental matrix solution at $z=0$ of the form \eqref{8luglio2021-4} with 
\begin{align}
&R=0
\\
\label{8luglio2021-5}
&(G_1)_{ij}=\frac{1}{\mu_j-\mu_i+1}(\mathcal{U}_0)_{ij},
\\
\noalign{\medskip}
\label{gk}
&(G_k)_{ij}=\frac{1}{\mu_j-\mu_i+k}\sum_{l_{k-1},\ldots,l_1=1}^{3}\frac{(\mathcal{U}_0)_{il_1}\ldots(\mathcal{U}_0)_{l_{k-1}j}}{(\mu_j-\mu_{l_1}+k-1)\ldots(\mu_j-\mu_{l_{k-1}}+1)},\quad k\ge 2,
\end{align}
where $i,j=1,2, 3$ and $\mu_1=\mu$, $\mu_2=0$, $\mu_3=-\mu$.
\epr
\begin{proof}
From the general theory sketched above we can choose $R=0$ and $G_1$ is  \eqref{8luglio2021-5}. To prove formula \eqref{gk} we proceed by induction: with the choice $R=0$, the recursive relations for $(G_k)_{ij}$, $k\ge 2$, reduce to
\[
(G_k)_{ij}=\frac{1}{\mu_j-\mu_i+k}(\mathcal{U}_0G_{k-1})_{ij}\mbox{.}
\]
Then, \eqref{gk} is easily verified for $k=2$. Let us suppose \eqref{gk} holds for $k-1$, then
\begin{align*}
(\mathcal{U}_0G_{k-1})_{ij}&=\sum_{l_{1}=1}^3(\mathcal{U}_0)_{il_{1}}(G_{k-1})_{l_{1}j}\\[2ex]&
 =\sum_{l_{1}=1}^3\frac{(\mathcal{U}_0)_{il_{1}}}{\mu_j-\mu_{l_{1}}+k-1}\sum_{l_{k-1},\ldots,l_2=1}^{3}\frac{(\mathcal{U}_0)_{l_1l_2}\ldots(\mathcal{U}_0)_{l_{k-2}j}}{(\mu_j-\mu_{l_2}+k-2)\ldots(\mu_j-\mu_{l_{k-2}}+1)}\mbox{,}
\end{align*}
proving \eqref{gk}.
\end{proof}

\bre
{\rm Due to the particular form of $\mathcal{U}_0$, the only non zero entries of $G_k$, $k\ge1$, are  those at positions $(1,1)$, $(1,3)$, $(3,1)$ and $(3,3)$.
}
\ere

Consider now the resonant case, $2\mu\in\mathbb{Z}\setminus\{0\}$. For $\mu>0$, the resonance $\mu_i-\mu_j\in\mathbb{Z}\setminus\{0\}$ occurs only for $i=1$ and $j=3$ at step $2\mu$ of the recursive construction if $\mu$ is half-integer, and also for $i=1$, $j=2$ and $i=2$, $j=3$ at step $\mu$ when $\mu$ is an integer. For  $\mu<0$, the above applies with $i,j$ exchanged.

\bpr
\label{3agosto2021-2}
If system \eqref{normal} with $\mu>0$ is resonant, then it has a fundamental matrix solution at $z=0$ of the form \eqref{8luglio2021-4}. 
For $\mu$ half-integer we can choose the matrix $R$ in such a way that the only (possibly) non zero entry is
\begin{align}
\label{r12}
R_{13}&=(\mathcal{U}_0)_{13},\quad \hbox{ for }\mu=\frac{1}{2},
\\
\label{rnorm}
R_{13}&=\sum_{m=1}^3\frac{(\mathcal{U}_0)_{1m}}{\mu_3-\mu_m+2\mu-1}\sum_{l_{2\mu-2},\ldots,l_1=1}^{3}\frac{(\mathcal{U}_0)_{ml_1}\ldots(\mathcal{U}_0)_{l_{2\mu-2}3}}{(\mu_3-\mu_{l_1}+2\mu-2)\ldots(\mu_3-\mu_{l_{2\mu-2}}+1)},\quad \mu\geq \frac{3}{2},
\end{align}
For $\mu$  positive integer
$$R=0.$$
 If $\mu$ is half-integer, the matrix coefficient $G_{2\mu}$ has entries $(1,1)$, $(3,1)$, $(3,3)$ fully determined and the entry $(1,3)$ is an arbitrary parameter; if $\mu$ is an integer the matrix coefficient $G_{\mu}$ has free parameters at entries $(1,2)$ and $(2,3)$ and another free parameter occurs in $G_{2\mu}$ at position $(1,3)$.
 
For $\mu<0$, the above results apply exchanging the entry $(1,3)$ with the entry $(3,1)$ for $\mu$ half-integer or exchanging the entries $(1,2)$, $(2,3)$, $(1,3)$ with the entries $(2,1)$, $(3,2)$, $(3,1)$, respectively, for $\mu$ integer.
\epr
\begin{proof}
From the general theory of Fuchsian singularities, we can choose $R$ such that  only $R_{13}$ is  possibly non zero. Let us start considering half-integer values of $\mu$. For the case $\mu=1/2$ we have
$$ (R_1)_{ij}=(\mathcal{U}_0)_{ij}\quad \hbox{ and } \quad 
(G_1)_{ij}\mbox{ arbitrary}.
$$
Let $\mu\geq 3/2$, then $(R)_{13}=(\mathcal{U}_0G_{2\mu-1})_{13}$. Now, for $k=1,\dots,2\mu-1$, formula \eqref{gk} holds, hence we easily obtain the sought expression of $(R)_{13}$. Again, the fact that $(G_{2\mu})_{13}$ is a free parameter is just a consequence of the general theory. Consider now the case of integer $\mu$. For $k=1,\ldots,\mu-1$ formula \eqref{gk} holds and we have $(R_\mu)_{12}=(\mathcal{U}_0G_{\mu-1})_{12}$, but $(\mathcal{U}_0)_{i2}=(\mathcal{U}_0)_{2j}=0$, for each $i$, $j=1$, $2$, $3$, hence $(R_\mu)_{12}=0$. Similarly we get $(R_\mu)_{23}=0$. From the general theory, $(G_\mu)_{12}=g_1$ and $(G_\mu)_{23}=g_2$ are free parameters.  Let us use the inductive definition to compute 
\[
(G_{\mu+1})_{ij}=\frac{1}{\mu_j-\mu_i+\mu+1}(\mathcal{U}_0G_{\mu})_{ij}\mbox{,}
\]
where the product $\mathcal{U}_0G_{\mu}$ has structure
\[
2\mathcal{U}_0G_\mu=
\begin{pmatrix}
B & g_1 & C\\
0 & 0 & 0\\
B & g_1 & C
\end{pmatrix}
\mbox{,}
\]
with $B=(G_\mu)_{11}+(G_\mu)_{31}$ and $C=(G_\mu)_{13}+(G_\mu)_{33}$. Now, from \eqref{gk} we see that $B=C=0$, thus $(\mathcal{U}_0G_k)_{11}=(\mathcal{U}_0G_k)_{31}=(\mathcal{U}_0G_k)_{13}=(\mathcal{U}_0G_k)_{33}=0$ for all $k\ge\mu$ and therefore $(R_{2\mu})_{13}=0$.

The last statement for $\mu<0$ is obvious. 
\end{proof}

Summing up, for the non-resonant case a fundamental solution at $z=0$ can be written in the form
\begin{equation}
\label{gggg}
\widetilde{Y}^{(0)}(z)=\mathcal{G}(z)~\hbox{\rm diag}(z^\mu,1,z^{-\mu})\mbox{,}
\quad
\quad
\mathcal{G}(z)=I+\frac{z}{2}
\begin{pmatrix}
1 & 0 &1/(1-2\mu)\\
0 & 0 & 0\\
1/(1+2\mu) & 0 & 1
\end{pmatrix}
+O(z^2)\mbox{,}
\end{equation}
while for the resonant case we can write a fundamental solution of the form
\[
\widetilde{Y}^{(0)}(z)=\mathcal{G}(z)~\hbox{\rm diag}(z^\mu,1,z^{-\mu})~z^R\mbox{,}
\]
where
\[
z^R=
{\small \begin{pmatrix}
1 & 0 & R_{13}\log(z)\\
0 & 1 & 0\\
0 & 0 & 1
\end{pmatrix}
}
,\quad \mu\ge\frac{1}{2};
\quad
\quad
\quad
z^R=
{\small \begin{pmatrix}
1 & 0 & 0\\
0 & 1 & 0\\
R_{13}\log(z) & 0 & 1
\end{pmatrix}
}
,\quad \mu\le-\frac{1}{2}\mbox{.}
\]
We recall that $R_{13}$ and $R_{31}$ depend on $\mu$. If $\mu\ne\pm 1/2,\pm 1$, the expansion of $\mathcal{G}(z)$ up to the first order is the same as in \eqref{gggg}, whereas if $\mu=\pm 1/2,\pm 1$ then $\mathcal{G}(z)$ has the form
$$
\mathcal{G}(z)=I+\frac{z}{2}
\begin{pmatrix}
1 & 0 & g\\
0 & 0 & 0\\
1/2 & 0 & 1
\end{pmatrix}
+O(z^2)\mbox{,}
\quad 
\mu=\frac{1}{2}\mbox{;}
\quad\quad
\mathcal{G}(z)=I+\frac{z}{2}
\begin{pmatrix}
1 & 0 & 1/2\\
0 & 0 & 0\\
g & 0 & 1
\end{pmatrix}
+O(z^2)\mbox{,}\quad \mu=-\frac{1}{2}\mbox{,}
$$
$$
\mathcal{G}(z)=I+\frac{z}{2}
\begin{pmatrix}
1 & g_1 & -1\\
0 & 0 & g_2\\
1/3 & 0 & 1
\end{pmatrix}
+O(z^2)\mbox{,}\quad \mu=1\mbox{;}
\quad\quad
\mathcal{G}(z)=I+\frac{z}{2}
\begin{pmatrix}
1 & 0 & 1/3\\
g_1 & 0 & 0 \\
-1 & g_2 & 1
\end{pmatrix}
+O(z^2)\mbox{,}\quad \mu=-1\mbox{,}
$$
being $g,g_1,g_2\in\mathbb{C}$ a free parameter.

\subsubsection{Central connection matrix}
The full monodromy data of the chamber of the  Dubrovin-Frobenius manifold associated to the branch \eqref{29luglio2021-6} include the constant {\bf central connection matrix} $\mathcal{C}^{(0)}$, defined by
\begin{equation}
\label{connectionmatrix}
\widetilde{Y}_1(z)=\widetilde{Y}^{(0)}(z)\mathcal{C}^{(0)}\mbox{,}
\end{equation}
where $\widetilde{Y}_1(z)$ is the fundamental solution in $\mathcal{S}_1$, related to the one of  Section \ref{StokesT3} by a change of normalization 
$$ 
\widetilde{Y}_1(z) = \hbox{\rm diag}\Bigl(\sqrt{2},~1/\sqrt{y_0^\prime},~\sqrt{2}\Bigr) \cdot \widetilde{Y}_1(z)^{\rm Sect.~ \ref{StokesT3}},
$$
with $\tilde{k}_1^0=i\sqrt{y_0^\prime}$, $\tilde{k}_2^0=i\sqrt{1-y_0^\prime}$, 
corresponding to  $\Psi=G_0^{\rm Sect.~ \ref{StokesT3}}\cdot\hbox{\rm diag}(1/\sqrt{2},\sqrt{y_0^\prime},1/\sqrt{2})$.

    The parameters of the confluent hypergeometric functions $U(z;a,b),M(z;a,b)$ are $a=\mu$ and $b=2\mu$, so we can  express $\widetilde{Y}_1(z)$ in terms of the {\bf Hankel functions} $H_\nu^{(1)}(z),H_\nu^{(2)}(z)$, with 
   $$
   \nu=\mu-\frac{1}{2}  
   $$
    (see Remark \ref{hankel} of Appendix \ref{specfunc}). We notice immediately that the second row of the connection matrix is
\[
\mathcal{C}^{(0)}_{21}=\frac{(\Omega_0)_{13}(\Omega_0)_{23}}{\mu^2\sqrt{y_0^\prime}}\equiv -\sqrt{1-y'_0},
\quad
\mathcal{C}^{(0)}_{22}=-\frac{1}{\sqrt{y_0^\prime}}\left(\frac{(\Omega_0)_{13}}{\mu}\right)^2 \equiv \sqrt{y_0^\prime},\quad\mathcal{C}^{(0)}_{23}=0\mbox{.}
\]
In order to compute the other rows, we have to distinguish again between the non-resonant and resonant cases.

Let us start with the non-resonant case. We need the series expansion of the Hankel functions in a neighborhood of $z=0$ when $\nu\not\in\mathbb{Z}$:
\begin{equation}
\label{h1}
H^{(1)}_\nu(z)=\sum_{k=0}^\infty\left(a_k^{(1)}(\nu)\left(\frac{z}{2}\right)^\nu-b_k^{(1)}(\nu)\left(\frac{z}{2}\right)^{-\nu}\right)z^{2k}\mbox{,}
\end{equation}
where
\[
a_k^{(1)}(\nu)=i\csc(\nu\pi)\frac{(-1)^k}{4^kk!}\frac{e^{-i\nu\pi}}{\Gamma(\nu+k+1)},\quad b_k^{(1)}(\nu)=i\csc(\nu\pi)\frac{(-1)^k}{4^kk!}\frac{1}{\Gamma(-\nu+k+1)}\mbox{,}
\]
and
\begin{equation}
\label{h2}
H^{(2)}_\nu(z)=\sum_{k=0}^\infty\left( a_k^{(2)}(\nu)\left(\frac{z}{2}\right)^{-\nu}-b_k^{(2)}(\nu)\left(\frac{z}{2}\right)^\nu\right)z^{2k}\mbox{,}
\end{equation}
where
\[
a_{k}^{(2)}(\nu)=i\csc(\nu\pi)\frac{(-1)^k}{4^kk!}\frac{1}{\Gamma(-\nu+k+1)},\quad b_{k}^{(2)}(\nu)=i\csc(\nu\pi)\frac{(-1)^k}{4^kk!}\frac{e^{i\nu\pi}}{\Gamma(\nu+k+1)}\mbox{.}
\]
Let us consider the entry $(1,1)$ of \eqref{connectionmatrix}: 
\[
\frac{-i\sqrt{\pi}(\Omega_0)_{13}e^{i\nu\pi/2}}{\mu}\left[\frac{i^{\nu} a_0^{(1)}(\nu)}{2^{\nu+\mu+1}}z^\mu-\frac{i^{-\nu}b_0^{(1)}(\nu)}{2^{-\nu-\mu+1}}\frac{z}{2}z^{-\mu}\right](1+O(z))=\left(\left(1+\frac{z}{2}\right)z^\mu\mathcal{C}^{(0)}_{11}+\frac{z^{-\mu}}{1-2\mu}\frac{z}{2}\mathcal{C}^{(0)}_{31}\right)(1+O(z^2))
\]
thus we can read the entries $(1,1)$ and $(3,1)$ of the connection matrix:
\[
\mathcal{C}^{(0)}_{11}=\displaystyle{-i\frac{\sqrt{2\pi y'_0}}{2^{1+2\mu}}\frac{\sec(\mu\pi)}{\Gamma(1/2+\mu)}},
\quad
\mathcal{C}^{(0)}_{31}=\displaystyle{i(1-2\mu)\frac{\sqrt{2\pi y^\prime_0}}{2^{2(1-\mu)}}\frac{\sec(\mu\pi)}{\Gamma(3/2-\mu)}}\mbox{.}
\]
The entries $(1,2)$ and $(3,2)$ have the same form with the substitution $y'_0\to 1-y'_0$:
\[
\mathcal{C}^{(0)}_{12}=-i\frac{\sqrt{2\pi(1-y'_0)}}{2^{1+2\mu}}\frac{\sec(\mu\pi)}{\Gamma(1/2+\mu)},
\quad
\mathcal{C}^{(0)}_{32}=i(1-2\mu)\frac{\sqrt{2\pi(1-y'_0)}}{2^{2(1-\mu)}}\frac{\sec(\mu\pi)}{\Gamma(3/2-\mu)}\mbox{.}
\]
The computation of the entries $(1,3)$ and $(3,3)$ is carried out with the same procedure using the expansion of $H_\nu^{(2)}(z)$. We obtain
\[
\mathcal{C}^{(0)}_{13}=\displaystyle{\frac{\sqrt{2\pi}}{2^{1+2\mu}}\frac{\sec(\mu\pi)}{\Gamma(1/2+\mu)}e^{i\mu\pi}},
\quad
\mathcal{C}^{(0)}_{33}=\displaystyle{(1-2\mu)\frac{\sqrt{2\pi}}{2^{2(1-\mu)}}\frac{\sec(\mu\pi)}{\Gamma(3/2-\mu)}e^{-i\mu\pi}}\mbox{.}
\]
The above $\mathcal{C}^{(0)}$, through the the relations \eqref{9luglio2021-1} with $R=0$, yields the Stokes matrices \eqref{29luglio2021-1}. 

\vskip 0.2 cm 
Next, we consider the resonant case $2\mu\in\mathbb{Z}\backslash\{0\}$. 
 
For integer $\mu$ the computations are the same as those of the previous non-resonant case. Through the the relations \eqref{9luglio2021-1} with $R=0$,  we obtain  the Stokes matrices \eqref{29luglio2021-1}, and more specifically 
$ 
\mathbb{S}_1=\mathbb{S}_2=I
$.

 We consider then half-integer values of $\mu$. In this case
$$
\nu=\mu-1/2=n\in\mathbb{Z},
 $$
  hence the local representations of the Hankel functions in a neighborhood of $z=0$ are different for the ones used in the precending section:
\begin{equation}
\label{h1n}
\begin{aligned}
H_n^{(1)}(z)=\left(\frac{z}{2}\right)^n\sum_{k=0}^\infty f^{(1)}_k(z)\left(\frac{z}{2}\right)^{2k}-\frac{i}{\pi}\left(\frac{z}{2}\right)^{-n}\sum_{k=0}^{n-1}\frac{(n-k-1)!}{k!}\left(\frac{z}{2}\right)^{2k}\mbox{,}
\end{aligned}
\end{equation}
where
\[
\begin{aligned}
f^{(1)}_k(z)=\frac{(-1)^k}{k!}\left[\frac{1+(2i/\pi)\log(z/2)}{\Gamma(n+k+1)}-\frac{i}{\pi}\frac{\psi(k+1)+\psi(n+k+1)}{(n+k)!}\right]\mbox{,}
\end{aligned}
\]
being $\psi(x)=\Gamma'(x)/\Gamma(x)$. For the second Hankel function we have:
\begin{equation}
\label{h2n}
\begin{aligned}
H_n^{(2)}(z)=\left(\frac{z}{2}\right)^n\sum_{k=0}^\infty f^{(2)}_k(z)\left(\frac{z}{2}\right)^{2k}+\frac{i}{\pi}\left(\frac{z}{2}\right)^{-n}\sum_{k=0}^{n-1}\frac{(n-k-1)!}{k!}\left(\frac{z}{2}\right)^{2k}\mbox{,}
\end{aligned}
\end{equation}
where
\[
\begin{aligned}
f^{(2)}_k(z)=\frac{(-1)^k}{k!}\left[\frac{1-(2i/\pi)\log(z/2)}{\Gamma(n+k+1)}+\frac{i}{\pi}\frac{\psi(k+1)+\psi(n+k+1)}{(n+k)!}\right]\mbox{.}
\end{aligned}
\]
Formulae \eqref{h1n} and \eqref{h2n} are valid for $n\in\mathbb{N}$. Nevertheless, we can use them also for negative integer $n$ recalling that
\[
H_{-n}^{(1)}(z)=(-1)^nH_{n}^{(1)}(z),
\quad 
H_{-n}^{(2)}(z)=(-1)^nH_{n}^{(2)}(z)
, \quad\quad n\in\mathbb{N}.
\]

 For brevity, we restrict to the case when $\mu>0$, being the case of negative $\mu$ analogous.  In this case $R_{13}\neq 0$.  We start studying the entry $(1,1)$ of \eqref{connectionmatrix}: the left hand side is
\begin{align*}
\frac{i\sqrt{\pi}(\Omega_0)_{13}}{\mu}\frac{e^{in\pi}}{2^{2\mu+1/2}(\mu-1/2)!}\left(2\log2+\frac{i}{\pi}\left(\psi\left(\mu+1/2\right)-\gamma\right)\right)(1+O(z))z^\mu+
\\[2ex]
\frac{\sqrt{\pi}(\Omega_0)_{13}}{\mu}\frac{e^{in\pi}}{2^{2\mu-1/2}\pi(\mu-1/2)!}(1+O(z))z^\mu\log z+O(z^2)(1+O(\log z))z^\mu+O(z)z^{-\mu}\mbox{,}
\end{align*}
while the right hand side is
\[
(1+O(z))(\mathcal{C}^{(0)})_{11}z^\mu+(\mathcal{C}^{(0)})_{31}(R)_{13}z^\mu\log z\mbox{.}
\]
Equating the two last relations and exploiting the dependence on $\mu$ and $y'_0$ we receive 
\[
\mathcal{C}^{(0)}_{11}=-\frac{e^{i\pi\mu}\sqrt{2y_0'}}{2^{2\mu+1}\sqrt{\pi}(\mu-1/2)!}\left(4\log 2+\psi(\mu+1/2)-\gamma\right),
\quad
\mathcal{C}^{(0)}_{31}=\frac{e^{i\pi\mu}\sqrt{2y_0'}}{2^{2\mu}(R)_{13}\sqrt{\pi}(\mu-1/2)!}\mbox{.}
\]
Similarly, we compute 
\[
\mathcal{C}^{(0)}_{12}=-\frac{e^{i\pi\mu}\sqrt{2(1-y_0')}}{2^{2\mu+1}\sqrt{\pi}(\mu-1/2)!}\left(4\log 2+\psi(\mu+1/2)-\gamma\right),
\quad
\mathcal{C}^{(0)}_{32}=\frac{e^{i\pi\mu}\sqrt{2(1-y_0')}}{2^{2\mu}(R)_{13}\sqrt{\pi}(\mu-1/2)!}\mbox{.}
\]
The entries $(1,3)$ and $(3,3)$ are computed in a similar way:
\[
\mathcal{C}^{(0)}_{13}=\frac{\sqrt{2\pi}}{2^{2\mu}(\mu-1/2)!}+\frac{i\sqrt{2}}{2^{2\mu+1}\sqrt{\pi}(\mu-1/2)!}(4\log 2+\psi(\mu+1/2)-\gamma),\quad\mathcal{C}^{(0)}_{33}=-\frac{i\sqrt{2}}{2^{2\mu}(R)_{13}\sqrt{\pi}(\mu-1/2)!}
\]

Substituting the above $\mathcal{C}^{(0)}$ into the r.h.s. of   \eqref{9luglio2021-1}, and also using \eqref{3agosto2021-1}, we receive the Stokes matrices and prove that the expressions \eqref{29luglio2021-1} hold also for half-integer $\mu$. Notice that for this computation, one also has to use the explicit value of $R_{13}$ as in Proposition \ref{3agosto2021-2}.

%{\color{blue} Check 
%\[
%\mathcal{C}^{(0)}\mathbb{S}_1(\mathcal{C}^{0})^T=e^{-i\pi R}e^{-i\pi \hat{\mu}}\eta^{-1}=\begin{pmatrix} -i\pi R_{13}e^{i\pi\mu} & 0 & e^{-i\pi\mu} \\ 0 & 1 & 0 \\ e^{i\pi\mu} & 0 & 0   \end{pmatrix}
%\]
%Gli elementi della seconda riga, della seconda colonna e l'elemento $(3,3)$ di $\mathcal{C}^{(0)}\mathbb{S}_1(\mathcal{C}^{(0)})^T$ non da $R$ e danno il risultato corretto. Gli elementi $(1,1), (1,3), (3,1)$ dipendono da $R$:
%\[
%(\mathcal{C}^{(0)}\mathbb{S}_1(\mathcal{C}^{(0)})^T)_{11}=\frac{2\pi}{2^{4\mu}(\mu-1/2)!^2 },
%\]
%per $\mu=1/2$, $R_{13}=1/2$ \`{e} $(\pi/2)$, ok, per $\mu=3/2$, $R_{13}=-1/32$ \`{e} $(\pi/32)$, ok;
%\[
%(\mathcal{C}^{(0)}\mathbb{S}_1(\mathcal{C}^{(0)})^T)_{13}=-\frac{2i}{2^{4\mu}(\mu-1/2)!^2 R_{13}},
%\]
%per $\mu=1/2$, $R_{13}=1/2$ \`{e} $(-i)$, ok, per $\mu=3/2$, $R_{13}=-1/32$ \`{e} $(i)$, ok;
%\[
%(\mathcal{C}^{(0)}\mathbb{S}_1(\mathcal{C}^{(0)})^T)_{31}=\frac{2i}{2^{4\mu}(\mu-1/2)!^2 R_{13}}=-(\mathcal{C}^{(0)}\mathbb{S}_1(\mathcal{C}^{(0)})^T)_{13},\mbox{ ok.}
%\]
%}

%%%%%%%%%%%%%%%%%%%%%%%%%%%%%%%%%%%%%%%%%%%%

\appendix

 %%%%%%%%%%%%%%%%%%%%%%%%%%%%%%%%%%%%%%%%%%%%.

\section{Appendix}
\label{specfunc}
It is an exercise to prove the following 
 \ble
 \label{21gennaio2021-3}
Any solution  of the system of PDEs
$$ 
\sum_{k=1}^3 \frac{\partial f}{\partial u_k}=0, 
\quad \quad 
\sum_{k=1}^3 u_k\frac{\partial f}{\partial u_k}=\alpha f, \quad\quad  \alpha\in\mathbb{C},
$$
has structure
$$ 
f=(u_2-u_1)^\alpha\mathcal{F}\left(\frac{u_2-u_1}{u_3-u_1}\right) ,\quad\hbox{ or }\quad f=(u_3-u_1)^\alpha \mathcal{G}\left(\frac{u_2-u_1}{u_3-u_1}\right),
$$
where $\mathcal{F}$ and $\mathcal{G}$ are some functions of their argument.

\ele

\subsection{Confluent hypergeometric functions}
\label{chf}
The confluent hypergeometric equation is
\[
zw''+(b-z)w'+aw=0,\quad a,b\in\mathbb{C}.
\]
Two linearly independent solutions are
\[
M(z;a,b):=\sum_{s=0}^\infty\frac{(a)_s}{s!\Gamma(b+s)}z^s,
\]
where $(a)_s$ denotes the Pochhammer symbol, and the function $U(z;a,b)$, which is uniquely determined by the asymptotic condition
\begin{equation}
\label{asymptoticsU}
U(z;a,b)\sim z^{-a},\quad z\to\infty,\quad -\frac{3}{2}\pi<\arg(z)<\frac{3}{2}\pi.
\end{equation}
The function $M(z;a,b)$ is an entire functions of $z,a,b$, while $U(z;a,b)$ has a branch point at $z=0$, all its braches being entire in $a,b$. The analytic continuation of $U(z;a,b)$ is given by the cyclic relation
\begin{equation}
\label{cyclicU}
U(ze^{2\pi i n})=\frac{2\pi i e^{-\pi i b n}\sin(\pi b n)}{\Gamma(1+a-b)\sin(\pi b)}M(z;a,b)+e^{-2\pi i b n}U(z;a,b),\quad-\pi<\arg(z)<\pi,\quad n\in\mathbb{Z}.
\end{equation}
The function $M(z;a,b)$ admits an asymptotic expansion as $z\to \infty$, $-\pi/2<\epsilon\arg(z)<3\pi/2$ given by
\begin{equation}
\label{asymptoticsM}
M(z;a,b)\sim\frac{e^z z^{a-b}}{\Gamma(a)}\sum_{s=0}^\infty\frac{(1-a)_s(b-a)_s}{s!}z^{-s}+\frac{e^{\epsilon\pi i a}z^{-a}}{\Gamma(b-a)}\sum_{s=0}^\infty\frac{(a)_s(a-b+1)_s}{s!}(-z)^{-s},
\end{equation}
where $\epsilon=-1,1$ and $a,b-a$ are not zero or a negative integer.

\bre
\label{hankel}
{\rm 
If $b=2a$, we can write the general solution of the hypergeometric equatio in terms of the Hankel functions $H_\nu^{(1)}(z),H_\nu^{(2)}(z)$, with $\nu=(b-1)/2$, through the following relations:
\begin{align*}
&U\left(-2iz;\nu+\frac{1}{2},2\nu+1\right)=\frac{i\sqrt{\pi}}{2}e^{i\pi\nu}(2z)^{-\nu}e^{-iz}H_\nu^{(1)}(z),\\
&U\left(2iz;\nu+\frac{1}{2},2\nu+1\right)=-\frac{i\sqrt{\pi}}{2}e^{-i\pi\nu}(2z)^{-\nu}e^{iz}H_\nu^{(2)}(z).
\end{align*}
The Hankel functions have the asymptotics
\[
H^{(1)}_\nu(z)\sim\sqrt{\frac{2}{\pi z}}e^{i\left(z-\nu\pi/2-\pi/4\right)}\mbox{,}\quad z\to\infty\mbox{,}\quad -\pi<\arg(z)<2\pi \mbox{,}
\]
\[
H^{(2)}_\nu(z)\sim\sqrt{\frac{2}{\pi z}}e^{-i\left(z-\nu\pi/2-\pi/4\right)}\mbox{,}\quad z\to\infty\mbox{,}\quad-2\pi<\arg(z)<\pi\mbox{.}
\]
}
\ere

\subsection{$(2,2)$ - Generalized hypergeometric functions}
The generalized hypergeometric equation of kind $(2,2)$ is
\[
z^2w'''+z(b_2+a_2z)w''+(b_1+a_1z)w+a_0w=0.
\]
If $b_1,b_2$ are not negative integers and $b_1-b_2$ is not an integer, then a fundamental set of solutions is
\begin{align*}
&w_0(z;\mathbf{a},\mathbf{b})={}_2F_2\left(\left.\begin{matrix} a_1,\,a_2\\ b_1,\, b_2 \end{matrix}\right| z \right), &w_1(z;\mathbf{a},\mathbf{b})=z^{1-b_1}{}_2F_2\left(\left.\begin{matrix} 1+a_1-b_1,\,1+a_2-b_1\\2- b_1,\,1+ b_2-b_1 \end{matrix}\right|z\right),\\
&&w_2(z;\mathbf{a},\mathbf{b})=z^{1-b_2}{}_2F_2\left(\left.\begin{matrix} 1+a_1-b_2,\,1+a_2-b_2\\1+ b_1-b_2,\,2- b_2 \end{matrix}\right|z\right),
\end{align*}
where
\[
{}_2F_2\left(\left.\begin{matrix} a_1,\,a_2\\ b_1,\, b_2 \end{matrix}\right| z \right)=\sum_{k=0}^\infty\frac{(a_1)_k(a_2)_k}{\Gamma(b_1+k)\Gamma(b_2+k)}\frac{z^k}{k!}
\]
is an entire function of $z$ and of the parameters $a_1,a_2,b_1,b_2$. The following asymptotics hold:
\[
{}_2F_2\left(\left.\begin{matrix} a_1,\,a_2\\ b_1,\, b_2 \end{matrix}\right| z \right)\sim\frac{1}{\Gamma(a_1)\Gamma(a_2)}\left[K_{2,2}(z)+L_{2,2}(ze^{i\epsilon\pi})\right],\quad z\to\infty,\quad -(2+\epsilon)\frac{\pi}{2}<\arg(z)<(2-\epsilon)\frac{\pi}{2},
\]
where $\epsilon=-1,1$,
\[
K_{2,2}(z)=e^z z^\gamma\sum_{k=0}^\infty d_kz^{-k},\quad d_0=1,\quad \gamma=\sum_{h=1,2}(a_h-b_h)
\]
(the recursive formulas for $d_k$, $k\ge 1$ can be found in \cite{Luke1}, formula (6) of section 5.11.3) and
\[
L_{2,2}(z)=\sum_{m=1,2}z^{-a_m}\sum_{k=0}^\infty c_{m,k}\frac{(-1)^kz^{-k}}{k!},\quad c_{m,k}=\Gamma(a_m+k)\frac{\Gamma(a_l-a_m-k)}{\prod_{n=1,2}\Gamma(b_n-a_m-k)},\,l\ne m.
\]

\section*{}  
{\bf Acknowledgements.} We thank M. Bertola, P. Boalch and M. Mazzocco for drawing our attention to some references. D. Guzzetti  is a member of  the European Union's H2020 research and innovation programme under the Marie Sk\l{l}odowska-Curie grant No. 778010 {\it IPaDEGAN}.

%\end{spacing}

\end{document}